\documentclass[aop]{imsart}

\usepackage[english]{babel}
\usepackage[T1]{fontenc}
\usepackage[applemac]{inputenc}
\usepackage{microtype}
\usepackage{bbm}
\usepackage{tikz}
\usepackage{graphicx}
\usepackage{subfloat}
\usepackage{amsmath}
\usepackage{amssymb}
\usepackage{amsthm}
\usepackage{latexsym}
\usepackage{color}
\usepackage{dsfont}
\usepackage{appendix}
\newcommand{\argmin}{\operatornamewithlimits{argmin}}
\usepackage{hyperref}
\usepackage{float}
\usepackage{enumitem}
\usepackage{mathabx}
\usepackage{wrapfig}
\usepackage{graphicx} 
\graphicspath{{figures/}}

\DeclareSymbolFont{calletters}{OMS}{cmsy}{m}{n}
\DeclareSymbolFontAlphabet{\mathcal}{calletters}

\def\be{\begin{eqnarray}}
\def\ee{\end{eqnarray}}

\def\b*{\begin{eqnarray*}}
\def\e*{\end{eqnarray*}}

\def \E{\mathbb{E}}

\def \L{\mathbb{L}}

\def \N{\mathbb{N}}
\def \P{{\mathbb P}}
\def \Q{{\mathbb Q}}
\def \R{\mathbb{R}}

\def \eps{\varepsilon}
\def \S{\Sigma}
\def \xb{{\bf{x}}}

\def\Ac{{\cal A}}
\def\Bc{{\cal B}}

\def\Dc{{\cal D}}

\def\Fc{{\cal F}}

\def\Kc{{\cal K}}
\def\Lc{{\cal L}}
\def\Mc{{\cal M}}
\def\Nc{{\cal N}}

\def\Pc{{\cal P}}

\def\Sc{{\cal S}}
\def\Tc{{\cal T}}
\def\Uc{{\cal U}}

\def\Xc{{\cal X}}
\def\Yc{{\cal Y}}

\def\Cfrak{{\mathfrak C}}

\def\Pb{\overline{\P}}

\def\xb{{\bar x}}
\def\yb{{\bar y}}

\def \ri{{\rm ri\hspace{0cm}}}
\def \rf{{\rm rf\hspace{0cm}}}
\def \cl{{\rm cl\hspace{0.05cm}}}
\def \interior{{\rm int\hspace{0cm}}}
\def \dom{{\rm dom}}
\def \supp{{\rm supp}}

\def \csupp{{\wideparen{\supp}}}

\def \conv{{\rm conv}}

\def \aff{{\rm aff}}

\def \Leb{{\L}}
\def \Ctn{{\rm{C}}}

\def \Ibf{{\mathbf I}}
\def \Sbf{{\mathbf S}}
\def \Tbf{{\mathbf T}}

\def \leb{{\lambda}}

\def \dist{{\rm dist}}

\def \muxpw{\mu{\otimes}{\rm pw}}

\def \Kcirc{{\wideparen{\Kc}}}

\def\no{\noindent}
\def\x{\times}

\def\={\;=\;}
\def\.{\;.}

\def\eps{\varepsilon}

\def \1{{\bf 1}}
\def \ep{\hbox{ }\hfill{ ${\cal t}$~\hspace{-5.1mm}~${\cal u}$   } }

\def \proof{{\noindent \bf Proof. }}
\def \ep{\hbox{ }\hfill$\Box$}

 \def\normeL2#1{\left\|{#1}\right\|_{L^2}}

\begin{document}

\title{Irreducible Convex Paving for Decomposition of Multi-dimensional Martingale Transport Plans
}

\author{Hadrien De March\thanks{hadrien.de-march@polytechnique.org.} \and Nizar Touzi\thanks{nizar.touzi@polytechnique.edu.}
        }

\begin{aug}

\affiliation{CMAP, \'Ecole Polytechnique, France\thanks{The authors gratefully acknowledge the financial support of the ERC 321111 Rofirm, and the Chairs Financial Risks (Risk Foundation, sponsored by Soci\'et\'e G\'en\'erale) and Finance and Sustainable Development (IEF sponsored by EDF and CA).}}

\end{aug}


\newtheorem{Theorem}{Theorem}[section]
\newtheorem{Lemma}[Theorem]{Lemma}
\newtheorem{Corollary}[Theorem]{Corollary}
\newtheorem{Proposition}[Theorem]{Proposition}
\newtheorem{Remark}[Theorem]{Remark}
\newtheorem{Example}[Theorem]{Example}
\newtheorem{Definition}[Theorem]{Definition}
\newtheorem{Assumption}[Theorem]{Assumption}

\begin{abstract}
: Martingale transport plans on the line are known from Beiglb\"ock \& Juillet \cite{beiglboeck2016problem} to have an irreducible decomposition on a (at most) countable union of intervals. We provide an extension of this decomposition for martingale transport plans in $\R^d$, $d\ge 1$. Our decomposition is a partition of $\R^d$ consisting of a possibly uncountable family of relatively open convex components, with the required measurability so that the disintegration is well-defined. We justify the relevance of our decomposition by proving the existence of a martingale transport plan filling these components. We also deduce from this decomposition a characterization of the structure of polar sets with respect to all martingale transport plans.
\end{abstract}

\maketitle


\section{Introduction}

The problem of martingale optimal transport was introduced as the dual of the problem of robust (model-free) superhedging of exotic derivatives in financial mathematics, see Beiglb\"ock, Henry-Labord\`ere \& Penkner \cite{beiglbock2013model} in discrete time, and Galichon, Henry-Labord\`ere \& Touzi \cite{galichon2014stochastic} in continuous-time. The robust superhedging problem was introduced  by Hobson \cite{hobson2011skorokhod}, and was addressing specific examples of exotic derivatives by means of corresponding solutions of the Skorohod embedding problem, see \cite{cox2011robust,hobson2015robust,hobson2012robust}, and the survey \cite{hobson2011skorokhod}. 

Given two probability measures $\mu,\nu$ on $\R^d$, with finite first order moment, martingale optimal transport differs from standard optimal transport in that the set of all coupling probability measures $\Pc(\mu,\nu)$ on the product space is reduced to the subset $\Mc(\mu,\nu)$ restricted by the martingale condition. We recall from Strassen \cite{strassen1965existence} that $\Mc(\mu,\nu)\neq\emptyset$ if and only if $\mu\preceq\nu$ in the convex order, i.e. $\mu(f)\le\nu(f)$ for all convex functions $f$. Notice that the inequality $\mu(f)\le\nu(f)$ is a direct consequence of the Jensen inequality, the reverse implication follows from the Hahn-Banach theorem.

This paper focuses on the critical observation by Beiglb\"ock \& Juillet \cite{beiglboeck2016problem} that, in the one-dimensional setting $d=1$, any such martingale interpolating probability measure $\P$ has a canonical decomposition $\P=\sum_{k\ge 0}\P_k$, where $\P_k\in\Mc(\mu_k,\nu_k)$ and $\mu_k$ is the restriction of $\mu$ to the so-called irreducible components $I_k$, and $\nu_k := \int_{x\in I_k}\P(dx,\cdot)$, supported in $J_k$, $k\ge 0$, is independent of the choice of $\P_k$. Here, $(I_k)_{k\ge 1}$ are open intervals, $I_0:=\R\setminus(\cup_{k\ge 1} I_k)$, and $J_k$ is an augmentation of $I_k$ by the inclusion of either one of the endpoints of $I_k$, depending on whether they are charged by the distribution $\P_k$. Remarkably, the irreducible components $(I_k,J_k)_{k\ge 0}$ are independent of the choice of $\P\in\Mc(\mu,\nu)$. To understand this decomposition, notice that convex functions in one dimension are generated by the family $f_{x_0}(x):=|x-x_0|,$ $x_0\in\R$, $x_0\in\R$. Then, in terms of the potential functions $U^\mu(x_0):=\mu(f_{x_0})$, and $U^\nu(x_0):=\nu(f_{x_0})$, $x_0\in\R$, we have $\mu\preceq\nu$ if and only if $U^\mu\le U^\nu$ and $\mu,\nu$ have same mean. Then, at any contact points $x_0$, of the potential functions, $U^\mu(x_0)=U^\nu(x_0)$, we have equality in the underlying Jensen equality, which means that the singularity $x_0$ of the underlying function $f_{x_0}$ is not seen by the measure. In other words, the point $x_0$ acts as a barrier for the mass transfer in the sense that martingale transport maps do not cross the barrier $x_0$. Such contact points are precisely the endpoints of the intervals $I_k$, $k\ge 1$.

The decomposition in irreducible components plays a crucial role for the quasi-sure formulation introduced by Beiglb\"ock, Nutz, and Touzi \cite{beiglbock2015complete}, and represents an important difference between martingale transport and standard transport. Indeed, while the martingale transport problem is affected by the quasi-sure formulation, the standard optimal transport problem is not changed. We also refer to Ekren \& Soner \cite{ekren2016constrained} for further functional analytic aspects of this duality.

Our objective in this paper is to extend the last decomposition to an arbitrary $d-$dimensional setting, $d\ge 1$. The main difficulty is that convex functions do not have anymore such a simple generating family. Therefore, all of our analysis is based on the set of convex functions. A first extension of the last decomposition to the multi-dimensional case was achieved by Ghoussoub, Kim \& Lim \cite{ghoussoub2015structure}. Motivated by the martingale monotonicity principle of Beiglb\"ock \& Juillet \cite{beiglboeck2016problem} (see also Zaev \cite{zaev2015monge} for higher dimension and general linear constraints), their strategy is to find a monotone set $\Gamma\subset\R^d\times\R^d$, where the robust superhedging holds with equality, as a support of the optimal martingale transport in $\Mc(\mu,\nu)$. Denoting $\Gamma_x:=\{y:(x,y)\in\Gamma\}$, this naturally induces the relation $x \,\mbox{Rel}\, x'$ if $x\in \ri\,\conv(\Gamma_{x'})$, which is then completed to an equivalence relation $\sim$. The corresponding equivalence classes define their notion of irreducible components.

Our subsequent results differ from \cite{ghoussoub2015structure} from two perspectives. First, unlike \cite{ghoussoub2015structure}, our decomposition is universal in the sense that it is not relative to any particular martingale measure in $\Mc(\mu,\nu)$ (see example \ref{example:Compodepend}). Second, our construction of the irreducible convex paving allows to prove the required measurability property, thus justifying completely the existence of a disintegration of martingale plans. 

Finally, during the final stage of writing the present paper, we learned about
the parallel work by Jan Ob{\l}{\'o}j and Pietro Siorpaes \cite{OblojSiopraes2017}.
Although the results are close, our approach is different from theirs.
We are grateful to them for pointing to us the notions of "convex
face" and "Wijsmann topology" and the relative references, which allowed us to streamline our presentation.
In an earlier version of this work we used instead a topology that we called the compacted Hausdorff distance, defined as the topology generated by the countable restrictions of the space to the closed balls centered in the origin with integer radii; the two are in our case the same topologies, as the Wijsman topology is locally equivalent to the Hausdorff topology in a locally compact set. We also owe Jan and Pietro special thanks for their useful remarks and comments on a first draft of this paper privately exchanged with them.

The paper is organized as follows. Section \ref{sect:mainresults} contains the main results of the paper, namely our decomposition in irreducible convex paving, and shows the identity with the Beiglb\"ock \& Juillet \cite{beiglboeck2016problem} notion in the one-dimensional setting. Section \ref{sect:preliminaries} collects the main technical ingredients needed for the statement of our main results, and gives the structure of polar sets. In particular, we introduce the new notions of relative face and tangent convex functions, together with the required topology on the set of such functions.  The remaining sections contain the proofs of these results. In particular, the measurability of our irreducible convex paving is proved in Section \ref{sect:measurability}.
\\

\no {\bf Notation}\quad We denote by $\overline\R:=\R\cup\{-\infty,\infty\}$ the completed real line, and similarly denote $\overline{\R}_+:=\R_+\cup \{\infty\}$. We fix an integer $d\ge 1$. For $x\in\R^d$ and $r\ge 0$, we denote $B_r(x)$ the closed ball for the Euclidean distance, centered in $x$ with radius $r$. We denote for simplicity $B_r := B_r(0)$. If $x\in\Xc$, and $A\subset \Xc$, where $(\Xc,{\rm d})$ is a metric space, $\dist(x,A):=\inf_{a\in A}{\rm d}(x,a)$. In all this paper, $\R^d$ is endowed with the Euclidean distance.

If $V$ is a topological affine space and $A\subset V$ is a subset of $V$, $\interior A$ is the interior of $A$, $\cl A$ is the closure of $A$, $\aff A$ is the smallest affine subspace of $V$ containing $A$, $\conv A$ is the convex hull of $A$, $\dim(A):=\dim(\aff A)$, and $\ri A$ is the relative interior of $A$, which is the interior of $A$ in the topology of $\aff A$ induced by the topology of $V$. We also denote by $\partial A:=\cl A\setminus\ri A$ the relative boundary of $A$, and by $\leb_A$ the Lebesgue measure of $\aff A$.

The set $\Kc$ of all closed subsets of $\R^d$ is a Polish space when endowed with the Wijsman topology\footnote{The Wijsman topology on the collection of all closed subsets of a metric space $(\Xc,{\rm d})$ is the weak topology generated by $\{\dist(x,\cdot):x\in\Xc\}$.} (see Beer \cite{beer1991polish}). As $\R^d$ is separable, it follows from a theorem of Hess \cite{hess1986contribution} that a function $F:\R^d\longrightarrow\Kc$ is Borel measurable with respect to the Wijsman topology if and only if its associated multifunction is Borel measurable, i.e.
\b*
F^-(V):=\{x\in\R^d:F(x)\cap V\neq\emptyset\}&\mbox{is Borel measurable}\\
&\mbox{for all open subset}~V\subset\R^d.
\e*
The subset $\Kcirc\subset\Kc$ of all the convex closed subsets of $\R^d$ is closed in $\Kc$ for the Wijsman topology, and therefore inherits its Polish structure. Clearly, $\Kcirc$ is isomorphic to $\ri\,\Kcirc := \{\ri K : K\in\Kcirc\}$ (with reciprocal isomorphism $\rm{cl}$). We shall identify these two isomorphic sets in the rest of this text, when there is no possible confusion.

We denote $\Omega:=\R^d\times\R^d$ and define the two canonical maps
\b*
 X :(x,y)\in\Omega
 \longmapsto x\in\R^d
 &\mbox{and}&
 Y :(x,y)\in\Omega
 \longmapsto y\in\R^d.
 \e*
For $\varphi,\psi:\R^d\longrightarrow\bar\R$, and $h:\R^d\longrightarrow\R^d$, we denote 
\b*
\varphi\oplus\psi
:=
\varphi(X)+\psi(Y),
&\mbox{and}&
h^\otimes := h(X)\cdot(Y-X),
\e*
with the convention $\infty-\infty = \infty$.

For a Polish space $\Xc$, we denote by $\Bc(\Xc)$ the collection of Borel subsets of $\Xc$, and $\Pc(\Xc)$ the set of all probability measures on $\big(\Xc,\Bc(\Xc)\big)$. For $\P\in\Pc(\Xc)$, we denote by $\Nc_\P$ the collection of all $\P-$null sets, $\supp\,\P$ the smallest closed support of $\P$, and $\csupp\P:=\cl\,\conv\,\supp\P$ the smallest convex closed support of $\P$. For a measurable function $f:\Xc\longrightarrow\R$, we denote $\dom\,f:=\{|f|<\infty\}$, and we use again the convention $\infty-\infty = \infty$ to define its integral, and denote
 \b*
 \P[f]:=\E^\P[f] 
 = \int_\Xc f d\P 
 = \int_\Xc f(x) \P(dx)
 &\mbox{for all}&
 \P\in\Pc(\Xc).
 \e*
Let $\Yc$ be another Polish space, and $\P\in\Pc(\Xc\x\Yc)$. The corresponding conditional kernel\footnote{The usual definition of a kernel requires that the map $x\mapsto \P_x[B]$ is Borel measurable for all Borel set $B\in\Bc(\R^d)$. In this paper, we only require this map to be analytically measurable.} $\P_x$ is defined $\mu-$a.e. by:
$$\P(dx,dy) = \mu(dx)\otimes \P_x(dy),\text{ where }\mu:=\P\circ X^{-1}.
$$
We denote by $\Leb^0(\Xc,\Yc)$ the set of Borel measurable maps from $\Xc$ to $\Yc$. We denote for simplicity $\Leb^0(\Xc):=\Leb^0(\Xc,\bar\R)$ and $\Leb^0_+(\Xc):=\Leb^0(\Xc,\bar\R_+)$. Let $\Ac$ be a $\sigma-$algebra of $\Xc$, we denote by $\Leb^\Ac(\Xc,\Yc)$ the set of $\Ac-$measurable maps from $\Xc$ to $\Yc$. For a measure $m$ on $\Xc$, we denote $\Leb^1(\Xc,m):=\{f\in\Leb^0(\Xc):m[|f|]<\infty\}$. We also denote simply $\Leb^1(m):=\Leb^1(\bar\R,m)$ and $\Leb^1_+(m):=\Leb^1_+(\bar\R_+,m)$.

We denote by $\Cfrak$ the collection of all finite convex functions $f:\R^d\longrightarrow\R$. We denote by $\partial f(x)$ the corresponding subgradient at any point $x\in\R^d$. We also introduce the collection of all measurable selections in the subgradient, which is nonempty by Lemma \ref{lemma:partialnonempty},
$$
\partial f:=\big\{p\in\L^0(\R^d,\R^d): p(x)\in\partial f(x)\text{ for all }x\in\R^d\big\}.
$$
We finally denote $\underline{f}_\infty := \liminf_{n\to\infty}f_n$, for any sequence $(f_n)_{n\ge 1}$ of real numbers, or of real-valued functions.

\section{Main results}
\label{sect:mainresults}
\setcounter{equation}{0}

Throughout this paper, we consider two probability measures $\mu$ and $\nu$ on $\R^d$ with finite first order moment, and $\mu \preceq \nu$ in the convex order, i.e. $\nu(f)\ge \mu(f)$ for all $f\in\Cfrak$. Using the convention $\infty-\infty=\infty$, we may define $(\nu-\mu)(f)\in[0,\infty]$ for all $f\in\Cfrak$. 

We denote by $\Mc(\mu,\nu)$ the collection of all probability measures on $\R^d\times\R^d$ with marginals $\P\circ X^{-1}=\mu$ and $\P\circ Y^{-1}=\nu$. Notice that $\Mc(\mu,\nu)\neq\emptyset$ by Strassen \cite{strassen1965existence}.

An $\Mc(\mu,\nu)-$polar set is an element of $\cap_{\P\in\Mc(\mu,\nu)}\Nc_\P$. A property is said to hold $\Mc(\mu,\nu)-$quasi surely (abbreviated as q.s.) if it holds on the complement of an $\Mc(\mu,\nu)-$polar set.

\subsection{The irreducible convex paving}

The next first result shows the existence of a maximum support martingale transport plan, i.e. a martingale interpolating measure $\widehat\P$ whose disintegration $\widehat\P_x$ has a maximum convex hull of supports among all measures in $\Mc(\mu,\nu)$.

\begin{Theorem}\label{thm:icp}
There exists $\widehat\P\in\Mc(\mu,\nu)$ such that
\be\label{eq:maximal_proba}
\mbox{for all }\P\in\Mc(\mu,\nu),&
 \csupp\,\P_X
 \subset
 \csupp\,\widehat\P_X, &\mu-\mbox{a.s.}
 \ee
 
 Furthermore $\csupp\,\widehat\P_X$ is $\mu-$a.s. unique, and we may choose this kernel so that
 \\
{\rm (i)} $x\longmapsto \csupp\,\widehat\P_x$ is analytically measurable\footnote{Analytically measurable means measurable with respect to the smallest $\sigma-$algebra containing the analytic sets. All Borel sets are analytic and all analytic sets are universally measurable, i.e. measurable with respect to all Borel measures (see Proposition 7.41 and Corollary 7.42.1 in \cite{bertsekas1978stochastic}).} $\R^d\longrightarrow \Kcirc,$ 
\\
{\rm (ii)}  $x\in I(x):=\ri\,\csupp\,\widehat\P_x$, for all $x\in\R^d$, and $\big\{I(x),x\in\R^d\big\}$ is a partition of $\R^d$.
\end{Theorem}

This Theorem will be proved in Subsection \ref{subsect:probamax}. The ($\mu-$a.s. unique) set valued map $I(X)$ paves $\R^d$ by its image by (ii) of Theorem \ref{thm:icp}. By \eqref{eq:maximal_proba}, this paving is stable by all $\P\in\Mc(\mu,\nu)$:
\be\label{eq:stability}
Y\in\cl I(X),&\Mc(\mu,\nu)-\mbox{q.s.}
\ee
Finally, the measurability of the map $I$ in the Polish space $\Kcirc$ allows to see it as a random variable, which allows to condition probabilistic events to $X\in I$, even when these components are all $\mu-$negligible when considered apart from the others. Under the conditions of Theorem \ref{thm:icp}, we call such $I(X)$ the irreducible convex paving associated to $(\mu,\nu)$.

Now we provide an important counterexample proving that for some $(\mu,\nu)$ in dimension larger than $1$, particular couplings in $\Mc(\mu,\nu)$ may define different pavings.

\begin{Example}\label{example:Compodepend}{\rm
In $\R^2$, we introduce $x_0:=(0,0)$, $x_1:=(1,0)$, $y_0:=x_0$, $y_{-1}:=(0,-1)$, $y_1:=(0,1)$, and $y_2:=(2,0)$. Then we set $\mu := \frac12(\delta_{x_0}+\delta_{x_1})$ and $\nu := \frac18(4\delta_{y_0}+\delta_{y_{-1}}+\delta_{y_1}+2\delta_{y_2})$. We can show easily that $\Mc(\mu,\nu)$ is the nonempty convex hull of $\P_1$ and $\P_2$ where
$$\P_1:=\frac{1}{8}
\big(4\delta_{x_0,y_0}+2\delta_{x_1,y_2}+\delta_{x_1,y_1}+\delta_{x_1,y_{-1}}\big)$$
and
$$\P_2:=\frac{1}{8}
\big(2\delta_{x_0,y_0}+\delta_{x_0,y_1}+\delta_{x_0,y_{-1}}+2\delta_{x_1,y_0}+2\delta_{x_1,y_2}\big)$$
\begin{figure}[h]
\centering
 \includegraphics[width=0.8\linewidth]{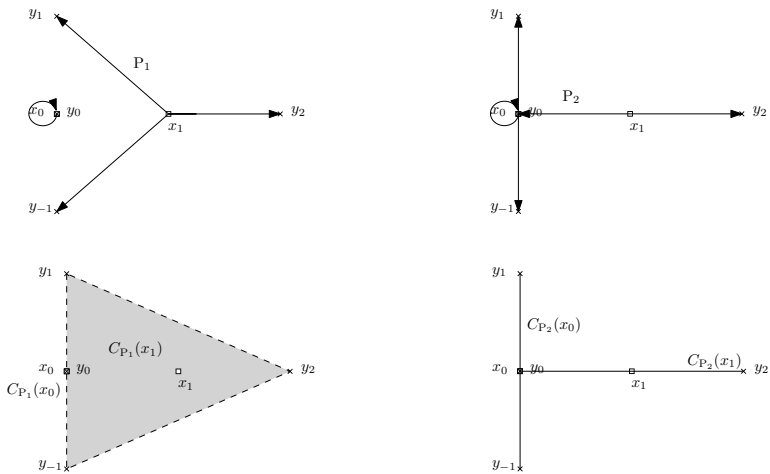}
    \caption{\label{fig:compGKL} The extreme probabilities and associated irreducible paving.}
\end{figure}

(i) {\it The Ghoussoub-Kim-Lim \cite{ghoussoub2015structure} (GKL, hereafter) irreducible convex paving.} Let $c_1 = \mathbf{1}_{\{X=Y\}}$, $c_2 = 1-c_1 = \mathbf{1}_{\{X\neq Y\}}$, and notice that $\P_i$ is the unique optimal martingale transport plan for $c_i$, $i=1,2$. Then, it follows that the corresponding $\P_i-$irreducible convex paving according to the definition of \cite{ghoussoub2015structure} are given by 
 $$
 \begin{array}{c}
 C_{\P_1}(x_0)=\{x_0\},
 ~C_{\P_1}(x_1)=\ri\,\conv\{y_1,y_{-1},y_2\},
 \\
 \mbox{and}~
 C_{\P_2}(x_0)=\ri\,\conv\{y_1,y_{-1}\},
 ~C_{\P_2}(x_1)=\ri\,\conv\{y_0,y_2\}.
 \end{array}
 $$
Figure \ref{fig:compGKL} shows the extreme probabilities $\P_1$ and $\P_2$, and their associated irreducible convex pavings map $C_{\P_1}$ and $C_{\P_2}$.

(ii) {\it Our irreducible convex paving.} The irreducible components are given by 
 \b*
 I(x_0)=\ri\,\conv(y_1,y_{-1})
 &\mbox{and}&
 I(x_1)=\ri\,\conv(y_1,y_{-1},y_2).
 \e*
To see this, we use the characterization of Proposition \ref{prop:CharacProb}. Indeed, as $\Mc(\mu,\nu) = \conv(\P_1,\P_2)$, for any $\P\in\Mc(\mu,\nu)$, $\P \ll \widehat\P:=\frac{\P_1+\P_2}{2}$, and $ {\rm supp}\,\P_x
 \subset
 {\rm conv}\big({\rm supp}\,\widehat\P_x\big)$ for $x=x_0,x_1$. Then $I(x) = \ri {\rm conv}\big({\rm supp}\,\widehat\P_x\big)$ for $x=x_0,x_1$ (i.e. $\mu$-a.s.) by Proposition \ref{prop:CharacProb}.}
\end{Example}

\begin{Remark}
In the one dimensional case, a convex paving which is invariant with respect to some $\P \in \Mc(\mu, \nu)$ is automatically invariant with respect to all $\P \in \Mc(\mu, \nu)$. Given a particular coupling $\P \in \Mc(\mu, \nu)$, the finest convex paving which is $\P-$invariant roughly corresponds to the GKL convex paving constructed in \cite{ghoussoub2015structure}. Then Example \ref{example:Compodepend} shows that this does not hold any more in dimension greater than two.

Furthermore, in dimension one the "restriction" $\nu_I:=\int_I\P(dx,\cdot)$ does not depend on the choice of the coupling $\P\in\Mc(\mu,\nu)$. Once again Example \ref{example:Compodepend} shows that it does not hold in higher dimension. Conditions guaranteeing that this property still holds in higher dimension will be investigated in \cite{HDM2017}.
\end{Remark}

\subsection{Behavior on the boundary of the components}

For a probability measure $\P$ on a topological space, and a Borel subset $A$, $\P|_A:=\P[\cdot\cap A]$ denotes its restriction to $A$.

\begin{Proposition}\label{prop:CharacProb}
We may choose $\widehat\P\in\Mc(\mu,\nu)$ in Theorem \ref{thm:icp} so that for all $\P\in\Mc(\mu,\nu)$ and $y\in\R^d$,
 \b*
&\mu\big[ \P_X[\{y\}]>0\big]
 \le 
 \mu\big[ \widehat \P_X[\{y\}]>0\big],\\
 \mbox{and}&
 \supp\,\P_X|_{\partial I(X)}
 \subset
 \csupp\,\widehat\P_X|_{\partial I(X)},~\mu-\mbox{a.s.}
 \e*

\no{\rm (i)} The set-valued maps $\underline{J}(X):= I(X)\cup \big\{y\in\R^d:\nu[{y}]>0,\mbox{ and }\widehat\P_X\big[\{y\}\big]>0\big\}$, and $\bar{J}(X):=I(X)\cup \csupp\,\widehat\P_X |_{\partial I(X)}$ are unique $\mu-$a.s, and $Y\in \bar{J}(X)$, $\Mc(\mu,\nu)-$q.s.

\no{\rm (ii)} We may chose the kernel $\widehat\P_X$ so that the map $\bar J$ is convex valued, $I\subset\underline{J}\subset \bar J\subset\cl I$, and both $\underline{J}$ and $\bar J$ are constant on $I(x)$, for all $x\in\R^d$.
\end{Proposition}

The proof is reported in Subsection \ref{subsect:probamax}.

\subsection{Structure of polar sets}

Here we state the structure of polar sets that is a direct consequence, and will be made more precise by Theorem \ref{thm:CharacPolar}.

\begin{Theorem}\label{thm:CharacPolarnotheta}
A Borel set $N\in\Bc(\Omega)$ is $\Mc(\mu,\nu)-$polar if and only if
$$
N \subset
\{X\in N_\mu\}\cup\{Y\in N_\nu\}\cup\{Y\notin J(X)\},
$$
 for some $(N_\mu,N_\nu)\in\Nc_\mu\x\Nc_\nu$ and a set valued map $J$ such that $\underline{J}\subset J\subset \bar J$, the map $J$ is constant on $I(x)$ for all $x\in\R^d$, $I(X)\subset \conv(J(X)\setminus N_\nu')$, $\mu-$a.s. for all $N_\nu'\in\Nc_\nu$, and $Y\in J(X)$, $\Mc(\mu,\nu)-$q.s.
\end{Theorem}

\subsection{The one-dimensional setting}

In the one-dimensional case, the decomposition in irreducible components and the structure of $\Mc(\mu,\nu)-$polar sets were introduced in Beiglb\"ock \& Juillet \cite{beiglboeck2016problem} and Beiglb\"ock, Nutz \& Touzi \cite{beiglbock2015complete}, respectively. 

Let us see how the results of this paper reduce to the known concepts in the one dimensional case. First, in the one-dimensional setting, $I(x)$ consists of open intervals (at most countable number) or single points. Following \cite{beiglboeck2016problem} Proposition 2.3, we denote the full dimension  components $(I_k)_{k\ge 1}$. 

We also have $\underline{J}=\bar J$ (see Proposition \ref{prop:IvsIBJ} below) therefore, Theorem \ref{thm:CharacPolarnotheta} is equivalent to Theorem 3.2 in \cite{beiglbock2015complete}. Similar to $(I_k)_{k\ge 1}$, we introduce the corresponding sequence $(J_k)_{k\ge 1}$, as defined in \cite{beiglbock2015complete}. Similar to \cite{beiglboeck2016problem}, we denote by $\mu_k$ the restriction of $\mu$ to $I_k$, and $\nu_k:=\int_{x\in I_k}\P[dx,\cdot]$ is independent of the choice of $\P\in\Mc(\mu,\nu)$. We define the Beiglb\"ock \& Juillet (BJ)-irreducible components
\b*
\left(I^{BJ},J^{BJ}\right):x\mapsto
\begin{cases}
(I_k,J_k)
&\mbox{if } x\in I_k, \mbox{ for some }k\ge 1, \\ 
\big(\{x\},\{x\}\big)
& \mbox{if } x\notin \cup_k I_k.
\end{cases}
\e*

\begin{Proposition}\label{prop:IvsIBJ}
Let $d=1$. Then $I=I^{BJ}$, and $\bar J=\underline{J} = J^{BJ}$, $\mu-a.s.$
\end{Proposition}

\proof
By Proposition \ref{prop:CharacProb} (i)-(ii), we may find $\widehat\P\in\Mc(\mu,\nu)$ such that
$\csupp\,\widehat\P_X
 =
 \cl I(X)$, and $\csupp\,\widehat\P_X|_{\partial I(X)} = \bar J\setminus I(X)$, $\mu-$a.s. Notice that as $\bar J\setminus I(\R^d)$ only consists of a countable set of points, we have $\underline{J} = \bar J$. By Theorem 3.2 in \cite{beiglbock2015complete}, we have $Y\in J^{BJ}(X)$, $\Mc(\mu,\nu)-$q.s. Therefore, $Y\in J^{BJ}(X)$, $\widehat \P-$a.s. and we have $\bar J(X)\subset J^{BJ}(X)$, $\mu-$a.s.

On the other hand, let $k\ge 1$. By the fact that $u_\nu-u_\mu>0$ on $I_k$, together with the fact that $J_k\setminus I_k$ is constituted with atoms of $\nu$, for any $N_\nu\in\Nc_\nu$, $J_k\subset\conv(J_k\setminus N_\nu)$. As $\mu = \nu$ outside of the components,
\be\label{eq:inclusionBJ}
J^{BJ}(X)\subset\conv(J^{BJ}(X)\setminus N_\nu), &\mu-\mbox{a.s.}
\ee
Then by Theorem 3.2 in \cite{beiglbock2015complete}, as $\{Y\notin \bar J(X)\}$ is polar, we may find $N_\nu\in\Nc_\nu$ such that $J^{BJ}(X)\setminus N_\nu\subset \bar J(X)$, $\mu-$a.s. The convex hull of this inclusion, together with \eqref{eq:inclusionBJ} gives the remaining inclusion $J^{BJ}(X)\subset \bar J(X)$, $\mu-$a.s.

The equality $I(X)=I^{BJ}(X)$, $\mu-$a.s. follows from the relative interior taken on the previous equality.
\ep
\section{Preliminaries}
\label{sect:preliminaries}
\setcounter{equation}{0}

The proof of these results needs some preparation involving convex analysis tools.

\subsection{Relative face of a set}

For a subset $A\subset\R^d$ and $a\in\R^d$, we introduce the face of $A$ relative to $a$ (also denoted $a-$relative face of $A$):
 \begin{equation}\label{notation:ri}
 \rf_a A  
 :=
 \big\{y\in A: (a-\eps(y-a),y+\eps(y-a)) \subset A,\text{ for some }\eps>0\big\}.
 \end{equation}
\begin{figure}[h]
\centering
 \includegraphics[width=0.8\linewidth]{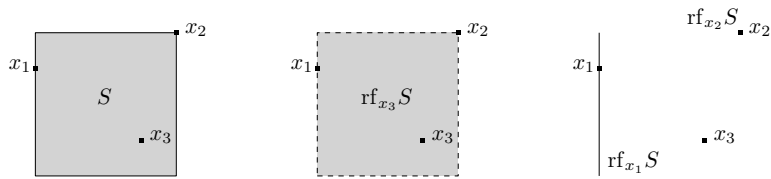}
    \caption{\label{fig:examplerf} Examples of relative faces.}
\end{figure}
\no Figure \ref{fig:examplerf} illustrates examples of relative faces of a square $S$, relative to some points. 
For later use, we list some properties whose proofs are reported in Section \ref{sect:convexanalysis}.
\footnote{
$\rf_aA$ is equal to the only relative interior of face of $A$ containing $a$, where we extend the notion of face to non-convex sets. A face $F$ of $A$ is a nonempty subset of $A$ such that for all $[a,b]\subset A$, with $(a,b)\cap F\neq \emptyset$, we have $[a,b]\subset F$. It is discussed in Hiriart-Urruty-Lemar\'echal \cite{hiriart2013convex} as an extension of Proposition 2.3.7 that when $A$ is convex, the relative interior of the faces of $A$ form a partition of $A$, see also Theorem 18.2 in Rockafellar \cite{rockafellar2015convex}.
}

\begin{Proposition}\label{prop:ri}
{\rm (i)} For $A,A'\subset\R^d$, we have $\rf_a(A\cap A') = \rf_a(A)\cap \rf_a(A')$, and $\rf_a A \subset \rf_a A'$ whenever $A\subset A'$. Moreover, $\rf_a A \neq\emptyset$ iff $a\in \rf_a A$ iff $a\in A$.
\\
{\rm (ii)} For a convex $A$, $\rf_a A  = \ri A \neq\emptyset$ iff $a\in\ri A$. Moreover, $\rf_a A$ is convex relatively open, $A\setminus \cl \rf_a A$ is convex, and if $x_0\in A\setminus \cl \rf_a A$ and $y_0\in A$, then $[x_0,y_0)\subset A\setminus \cl\rf_a A$. Furthermore, if $a\in A$, then $\dim(\rf_a\cl A) =\dim(A)$ if and only if $a\in \ri A$. In this case, we have $\cl\rf_a\cl A = \cl\ri\,\cl A = \cl A = \cl\rf_aA$.
\end{Proposition}

\subsection{Tangent Convex functions}

Recall the notation \eqref{notation:ri}, and denote for all $\theta:\Omega\to\bar\R$: 
 \b*
 \dom_x\theta
 &:=&
 \rf_x\conv\,\dom\,\theta(x,\cdot).
 \e*
For $\theta_1,\theta_2:\Omega\longrightarrow\R$, we say that $\theta_1=\theta_2$, $\muxpw$, if 
 \b*
 \dom_X\theta_1 = \dom_X\theta_2,
 &\mbox{and}&
 \theta_1(X,\cdot)=\theta_2(X,\cdot)
 ~\mbox{on}~\dom_X\theta_1,~\mu-\mbox{a.s.}
 \e*

The crucial ingredient for our main result is the following.

\begin{Definition}
A measurable function $\theta:\Omega\to\overline{\R}_+$ is a tangent convex function if 
 \b*
 &\theta(x,\cdot)
  ~\mbox{is convex, and}~
  \theta(x,x)=0,
 ~\mbox{for all}~
 x\in\R^d.
 \e*
We denote by $\Theta$ the set of tangent convex functions, and we define
 \b*
 \Theta_\mu
 &:=&
 \big\{\theta\in\L^0(\Omega,\overline{\R}_+):
       \theta = \theta',~\muxpw,
       \mbox{and}~\theta\ge\theta',~\mbox{for some}~
       \theta'\in\Theta   
 \big\}.
 \e*
\end{Definition}

In order to introduce our main example of such functions, let
\b*
\Tbf_pf(x,y)
:=
f(y)-f(x)-p^\otimes(x,y)\ge 0 ,
&\mbox{for all}&
f\in\Cfrak,~\mbox{and}~p\in\partial f.
\e*
Then, $
 \Tbf (\Cfrak):=\{\Tbf_p f:f\in\Cfrak,p\in\partial f\}\subset\Theta\subset\Theta_\mu
 $.

\begin{Example}\label{expl:Thetamununonconvex}
The second inclusion is strict. Indeed, let $d=1$, and consider the convex function $f:=\infty\mathbf{1}_{(-\infty,0)}$. Then $\theta' := f(Y-X)\in\Theta$. Now let $\theta = \theta'+\sqrt{|Y-X|}$. Notice that since $\dom_X\theta' = \dom_X\theta = \{X\}$, we have $\theta' = \theta$, $\muxpw$ for any measure $\mu$, and $\theta\ge \theta'$. Therefore $\theta\in\Theta_\mu$. However, for all $x\in\R^d$, $\theta(x,\cdot)$ is not convex, and therefore $\theta\notin\Theta$.

In higher dimension we may even have $X\in\ri\,\dom\theta(X,\cdot)$, and $\theta(X,\cdot)$ is not convex. Indeed, for $d=2$, let $f:(y_1,y_2)\longmapsto \infty(\mathbf{1}_{\{|y_1|> 1\}} + \mathbf{1}_{\{|y_2|> 1\}})$, so that $\theta := f(Y-X)\in\Theta$. Let $x_0:=(1,0)$ and $\theta := \theta' + \mathbf{1}_{\{Y = X+x_0\}}$. Then, $\theta = \theta'$, $\muxpw$ for any measure $\mu$, and $\theta\ge \theta'$. Therefore $\theta\in\Theta_\mu$. However, $\theta\notin\Theta$ as $\theta(x,\cdot)$ is not convex for all $x\in\R^d$.
\end{Example}

\begin{Proposition}\label{prop:Thetamu}
{\rm (i)} Let $\theta\in\Theta_\mu$, then $\dom_X\theta=\rf_X\dom\theta(X,\cdot)\subset\dom\theta(X,\cdot)$, $\mu-$a.s.\\
{\rm (ii)} Let $\theta_1,\theta_2\in\Theta_\mu$, then $\dom_X(\theta_1+\theta_2) = \dom_X\theta_1\cap\dom_X\theta_2$, $\mu-$a.s.\\
{\rm (iii)} $\Theta_\mu$ is a convex cone.
\end{Proposition}
\proof
\no\underline{\rm (i)} It follows immediately from the fact that $\theta(X,\cdot)$ is convex and finite on $\dom_X\theta$, $\mu-$a.s. by definition of $\Theta_\mu$. Then $\dom_X\theta\subset\rf_X\dom\theta(X,\cdot)$. On the other hand, as $\dom\theta(X,\cdot)\subset \conv\,\dom\theta(X,\cdot)$, the monotony of $\rf_x$ gives the other inclusion: $\rf_X\dom\theta(X,\cdot)\subset\dom_X\theta$.
\\
\no\underline{\rm (ii)} As $\theta_1,\theta_2\!\ge\! 0$, $\dom(\theta_1+\theta_2) \!=\! \dom\theta_1\cap\dom\theta_2$. Then, for $x\in\R^d$, $\conv\,\dom(\theta_1(x,\cdot)\!+\!\theta_2(x,\cdot)) \!\subset\! \conv\,\dom\theta_1(x,\cdot)\cap\conv\,\dom\theta_2(x,\cdot)$. By Proposition \ref{prop:ri} (i),
 \b*
 \dom_x(\theta_1+\theta_2) 
 \subset 
 \dom_x\theta_1\cap\dom_x\theta_2,
 &\mbox{for all}&
 x\in\R^d.
 \e*
As for the reverse inclusion, notice that (i) implies that $\dom_X\theta_1\cap\dom_X\theta_2\subset \dom\theta_1(X,\cdot)\cap\dom\theta_2(X,\cdot) = \dom\big(\theta_1(X,\cdot)+\theta_2(X,\cdot)\big)\subset\conv\,\dom\big(\theta_1(X,\cdot)+\theta_2(X,\cdot)\big)$, $\mu-$a.s. Observe that $\dom_x\theta_1\cap\dom_x\theta_2$ is convex, relatively open, and contains $x$. Then, 
 \b*
 \dom_X\theta_1\cap\dom_X\theta_2
 &=&
 \rf_X\big(\dom_X\theta_1\cap\dom_X\theta_2\big)\\
 &\subset&
 \rf_X\Big(\conv\,\dom\big(\theta_1(X,\cdot)+\theta_2(X,\cdot)\big)\Big)\\
 &=&
 \dom_X(\theta_1+\theta_2)
 ~~\mu-\mbox{a.s.}
 \e*
\no\underline{\rm (iii)} Given (ii), this follows from direct verification.
\ep

\begin{Definition}
A sequence $(\theta_n)_{n\ge 1}\subset\L^0(\Omega)$ converges $\muxpw$ to some $\theta\in\L^0(\Omega)$ if 
 \b*
 \dom_X\left(\underline\theta_\infty\right) =  \dom_X\theta &\mbox{and}&\theta_n(X,\cdot) \longrightarrow  \theta(X,\cdot),\,\mbox{pointwise on}~\dom_X\theta,
 ~\mu-\mbox{a.s.}
 \e*
\end{Definition}

Notice that the $\muxpw$-limit is $\muxpw$ unique. In particular, if $\theta_n$ converges to $\theta$, $\muxpw$, it converges as well to $\underline\theta_\infty$.

\begin{Proposition}\label{prop:convergenceprop}
Let $(\theta_n)_{n\ge 1}\subset \Theta_\mu$, and $\theta:\Omega\longrightarrow \bar\R_+$, such that $\theta_n\underset{n\to\infty}{\longrightarrow} \theta$, $\muxpw$,\\
{\rm (i)} $\dom_X\theta\subset \liminf_{n\to\infty}\dom_X\theta_n$, $\mu-$a.s.\\
{\rm (ii)} If $\theta_n' = \theta_n$, $\muxpw$, and $\theta_n'\ge \theta_n$, then $\theta_n'\underset{n\to\infty}{\longrightarrow} \theta$, $\muxpw$;\\
{\rm (iii)} $\underline{\theta}_\infty\in\Theta_\mu$.
\end{Proposition}
\proof
 \no\underline{\rm (i)} Let $x\in\R^d$, such that $\theta_n(x,\cdot)$ converges on $\dom_x\theta$ to $\theta(x,\cdot)$. Let $y\in\dom_x\theta$, let $y'\in\dom_x\theta$ such that $y' = x-\epsilon(y-x)$, for some $\epsilon>0$. As $\theta_n(x,y)\underset{n\to\infty}{\longrightarrow}\theta(x,y)$, and $\theta_n(x,y')\underset{n\to\infty}{\longrightarrow}\theta(x,y')$, then for $n$ large enough, both are finite, and $y\in \dom_x\theta_n$. $y\in\liminf_{n\to\infty}\dom_x\theta_n$, and $\dom_x\theta\subset\liminf_{n\to\infty}\dom_x\theta_n$. The inclusion is true for $\mu-$a.e. $x\in\R^d$, which gives the result.
 
 \no\underline{\rm (ii)} By (i), we have $\dom_X\theta\subset\liminf_{n\to\infty}\dom_X\theta_n = \liminf_{n\to\infty}\dom_X\theta_n'$, $\mu-$a.s. As $\theta_n\le\theta_n'$, $\dom_X\underline{\theta'}_\infty\subset\dom_X\underline{\theta}_\infty\subset\liminf_{n\to\infty}\dom_X\theta_n$, $\mu-$a.s. We denote $N_\mu\in\Nc_\mu$, the set on which $\theta_n(X,\cdot)$ does not converge to $\theta(X,\cdot)$ on $\dom_X\theta(X,\cdot)$. For $x\notin N_\mu$, for $y\in\dom_x\theta$, $\theta_n(x,y)=\theta_n'(x,y)$, for $n$ large enough, and $\theta_n'(x,y)\underset{n\to\infty}{\longrightarrow}\theta(x,y)<\infty$. Then $\dom_X\theta = \dom_X\underline{\theta'}_\infty$, and $\theta_n'(X,\cdot)$ converges to $\theta(X,\cdot)$, on $\dom_X\theta$, $\mu-$a.s. We proved that $\theta_n'\underset{n\to\infty}{\longrightarrow}\theta$, $\muxpw$.
 
\no\underline{\rm (iii)} has its proof reported in Subsection \ref{subsect:Komlos} due to its length and technicality.
\ep\\

The next result shows the relevance of this notion of convergence for our setting.

\begin{Proposition}\label{prop:Komlos}
Let $(\theta_n)_{n\ge 1}\subset\Theta_\mu$. Then, we may find a sequence $\widehat\theta_n\in \conv(\theta_k,k\ge n)$, and $\widehat\theta_\infty\in\Theta_\mu$ such that $\widehat\theta_n \longrightarrow\widehat\theta_\infty$, $\muxpw$ as $n\to\infty$.
\end{Proposition}

The proof is reported in Subsection \ref{subsect:Komlos}.

\begin{Definition}
{\rm (i)} A subset $\Tc\subset\Theta_\mu$ is $\muxpw$-Fatou closed if $\underline\theta_\infty\in\Tc$ for all $(\theta_n)_{n\ge 1}\subset\Tc$ converging $\muxpw$ (in particular, $\Theta_\mu$ is $\muxpw-$Fatou closed by Proposition \ref{prop:convergenceprop} {\rm (iii)}).
\\
{\rm (ii)} The $\muxpw-$Fatou closure of a subset $A\subset\Theta_\mu$ is the smallest $\muxpw-$Fatou closed set containing $A$:
 $$
 \widehat A
 :=
 \bigcap\big\{\Tc\subset\Theta_\mu:~
              A\subset\Tc,
              ~\text{and}~
              \Tc~\mbox{$\muxpw$-Fatou closed}\,
         \big\}.
 $$
\end{Definition}

We next introduce for $a\ge 0$ the set $\Cfrak_a:=\big\{f\in\Cfrak:(\nu-\mu)(f)\le a\big\}$, and
$$\widehat{\Tc}(\mu,\nu):=\underset{a\ge 0}{\bigcup}\,\widehat{\Tc}_a,
\text{ where }
\widehat{\Tc}_a:=\widehat{\Tbf(\Cfrak_a)}, ~\mbox{and}~
 \Tbf \big(\Cfrak_a\big)
 :=
 \big\{ \Tbf_p f: f\in\Cfrak_a,p\in\partial f\big\}.
 $$

\begin{Proposition}\label{prop:TmunuConvexCone}
$\widehat{\Tc}(\mu,\nu)$ is a convex cone.
\end{Proposition}
\proof
We first prove that $\widehat{\Tc}(\mu,\nu)$ is a cone. We consider $\lambda,a>0$, as we have $\lambda\Cfrak_a=\Cfrak_{\lambda a}$, and as convex combinations and inferior limit commute with the multiplication by $\lambda$, we have $\lambda\widehat\Tc_a=\widehat\Tc_{\lambda a}$. Then $\widehat{\Tc}(\mu,\nu) = {\rm cone}(\widehat{\Tc}_1)$, and therefore it is a cone.

We next prove that $\widehat\Tc_a$ is convex for all $a\ge 0$, which induces the required convexity of $\widehat{\Tc}(\mu,\nu)$ by the non-decrease of the family $\{\widehat\Tc_a,a\ge 0\}$. Fix $0\le\lambda\le 1$, $a\ge 0$, $\theta_0\in\widehat{\Tc}_a$, and denote $\Tc(\theta_0):=\big\{\theta\in\widehat{\Tc}_a:\lambda\theta_0+(1-\lambda)\theta\in\widehat{\Tc}_a\big\}$. In order to complete the proof, we now verify that $\Tc(\theta_0)\supset\Tbf \big(\Cfrak_a\big)$ and is $\muxpw-$Fatou closed, so that $\Tc(\theta_0)=\widehat\Tc_a$. 

To see that $\Tc(\theta_0)$ is Fatou-closed, let $(\theta_n)_{n\ge 1}\subset\Tc(\theta_0)$, converging $\muxpw$. By definition of $\Tc(\theta_0)$, we have $\lambda\theta_0+(1-\lambda)\theta_n\in\widehat{\Tc}_a$ for all $n$. Then, $\lambda\theta_0+(1-\lambda)\theta_n\longrightarrow \liminf_{n\to\infty}\lambda\theta_0+(1-\lambda)\underline\theta_n$, $\muxpw$, and therefore $\lambda\theta_0+(1-\lambda)\underline\theta_\infty\in \widehat{\Tc}_a$, which shows that $\underline\theta_\infty\in\Tc(\theta_0)$. 

We finally verify that $\Tc(\theta_0)\supset\Tbf \big(\Cfrak_a\big)$. First, for $\theta_0\in\Tbf \big(\Cfrak_a\big)$, this inclusion follows directly from the convexity of $\Tbf \big(\Cfrak_a\big)$, implying that $\Tc(\theta_0)=\widehat\Tc_a$ in this case. For general $\theta_0\in\widehat{\Tc}_a$, the last equality implies that $\Tbf \big(\Cfrak_a\big)\subset\Tc(\theta_0)$, thus completing the proof. 
\ep\\

Notice that even though $\Tbf(\Cfrak_a)\subset\Theta$, the functions in $\widehat{\Tc}(\mu,\nu)$ may not be in $\Theta$ as they may not be convex in $y$ on $(\dom_x\theta)^c$ for some $x\in\R^d$ (see Example \ref{expl:Thetamununonconvex}). The following result shows that some convexity is still preserved.

\begin{Proposition}\label{prop:ConvexityInvariance}
For all $\theta\in\widehat{\Tc}(\mu,\nu)$, we may find $N_\mu\in\Nc_\mu$ such that for $x_1,x_2\notin N_\mu$, $y_1,y_2\in\R^d$, and $\lambda\in[0,1]$ with $\yb:= \lambda y_1 + (1-\lambda)y_2\in\dom_{x_1}\theta\cap\dom_{x_2}\theta$, we have:
 \b*
 \lambda \theta(x_1,y_1)+(1-\lambda)\theta(x_1,y_2)-\theta(x_1,\yb)
 &=&
\lambda \theta(x_2,y_1)+(1-\lambda)\theta(x_2,y_2)\\
&&-\theta(x_2,\yb)
 \ge 0.
 \e*
\end{Proposition}

The proof of this claim is reported in Subsection \ref{subsect:convexityinvariance}. We observe that the statement also holds true for a finite number of points $y_1,...,y_k$.\footnote{
This is not a direct consequence of Proposition \ref{prop:ConvexityInvariance}, as the barycentre $\yb$ has to be in $\dom_{x_1}\theta\cap\dom_{x_2}\theta$.
}

\subsection{Extended integral}

We now introduce the extended $(\nu-\mu)-$integral:
 \b*
 \nu\widehat{\ominus}\mu[\theta]
 :=
 \inf\big\{a\ge 0 :\theta\in \widehat{\Tc}_a\big\} 
 &\mbox{for}&
 \theta\in\widehat{\Tc}(\mu,\nu).
 \e*

\begin{Proposition}\label{prop:numoinsmu}
{\rm (i)} $\P[\theta]\leq \nu\widehat{\ominus}\mu[\theta]<\infty$ for all $\theta\in\widehat{\Tc}(\mu,\nu)$ and $\P\in\Mc(\mu,\nu)$.\\
{\rm (ii)}  $\nu\widehat{\ominus}\mu[\Tbf_{p} f] = (\nu-\mu)[f]$ for $f\in\Cfrak\cap\L^1(\nu)$ and $p\in\partial f$.\\
{\rm (iii)} $\nu\widehat{\ominus}\mu$ is homogeneous and convex.
\end{Proposition}
\proof
\no\underline{\rm (i)} For $a>\nu\widehat{\ominus}\mu[\theta]$, set $S^a:=\big\{F\in\Theta_\mu:\P[F]\leq a~\mbox{for all}~\P\in\Mc(\mu,\nu)\big\}$. Notice that $S^a$ is $\muxpw-$Fatou closed by Fatou's lemma, and contains $\Tbf(\Cfrak_a)$, as for $f\in\Cfrak\cap\Leb^1(\nu)$ and $p\in\partial f$, $\P[T_pf]= (\nu-\mu)[f]$ for all $\P\in\Mc(\mu,\nu)$. Then $S^a$ contains $\widehat{\Tc}_a$ as well, which contains $\theta$. Hence, $\theta\in S^a$ and $\P[\theta]\leq a$ for all $\P\in\Mc(\mu,\nu)$. The required result follows from the arbitrariness of $a>\nu\widehat{\ominus}\mu[\theta]$.

\no\underline{\rm (ii)} Let $\P\in\Mc(\mu,\nu)$. For $p\in\partial f$, notice that $T_{p} f\in \Tbf(\Cfrak_a)\subset\widehat{\Tc_{a}}$ for some $a=(\nu-\mu)[f]$, and therefore $(\nu-\mu)[f] \ge \nu\widehat{\ominus}\mu[T_{p} f]$. Then, the result follows from the inequality $(\nu-\mu)[f] = \P[T_{p} f] \leq \nu\widehat{\ominus}\mu[T_{p} f]$.

\no\underline{\rm (iii)} Similarly to the proof of Proposition \ref{prop:TmunuConvexCone}, we have $\lambda\widehat\Tc_a=\widehat\Tc_{\lambda a}$, for all  $\lambda,a>0$. Then with the definition of $\nu\widehat{\ominus}\mu$ we have easily the homogeneity.

To see that the convexity holds, let $0<\lambda<1$, and $\theta,\theta'\in\widehat{\Tc}(\mu,\nu)$ with $a>\nu\widehat{\ominus}\mu[\theta]$, $a'>\nu\widehat{\ominus}\mu[\theta']$, for some $a,a'>0$. By homogeneity and convexity of $\widehat\Tc_1$, $\lambda \theta +(1-\lambda)\theta'\in\widehat\Tc_{\lambda a + (1-\lambda)a'}$, so that $\nu\widehat{\ominus}\mu[\lambda \theta +(1-\lambda)\theta']\le\lambda a + (1-\lambda)a'$. The required convexity property now follows from arbitrariness of $a>\nu\widehat{\ominus}\mu[\theta]$ and $a'>\nu\widehat{\ominus}\mu[\theta']$.
\ep \\

The following compacteness result plays a crucial role.

\begin{Lemma}\label{lemma:KomlosLiminf}
Let $(\theta_n)_{n\ge 1}\subset\widehat{\Tc}(\mu,\nu)$ be such that
$\sup_{n\ge 1}\,\nu\widehat{\ominus}\mu(\theta_n)<\infty$.
Then we can find a sequence $\widehat\theta_n\in \conv(\theta_k,k\ge n)$ such that
\b*
\underline{\widehat\theta}_\infty\in\widehat{\Tc}(\mu,\nu),&
\widehat\theta_n \longrightarrow\underline{\widehat\theta}_\infty,
~\muxpw,
&\mbox{and}~~
\nu\widehat{\ominus}\mu(\underline{\widehat\theta}_\infty)
\leq
\liminf_{n\to\infty}\nu\widehat{\ominus}\mu(\theta_n).
\e*
\end{Lemma}
\proof
By possibly passing to a subsequence, we may assume that $\lim_{n\to\infty}(\nu\widehat{\ominus}\mu)(\theta_n)$ exists. The boundedness of $\nu\widehat{\ominus}\mu(\theta_n)$ ensures that this limit is finite. We next introduce the sequence $\widehat\theta_n$ of Proposition \ref{prop:Komlos}. Then $\widehat\theta_n \longrightarrow\widehat\theta_\infty$, $\mu\otimes{\rm pw}$, and therefore $\underline{\widehat\theta}_\infty\in\widehat{\Tc}(\mu,\nu)$, because of the convergence $\widehat\theta_n \longrightarrow\underline{\widehat\theta}_\infty$, $\muxpw$. As $(\nu\widehat{\ominus}\mu)(\widehat\theta_n)\le \sup_{k\ge n}(\nu\widehat{\ominus}\mu)(\theta_k)$ by Proposition \ref{prop:numoinsmu} (iii), we have $\infty>\lim_{n\to\infty}(\nu\widehat{\ominus}\mu)(\theta_n)=\lim_{n\to\infty}\sup_{k\ge n}(\nu\widehat{\ominus}\mu)(\theta_k)\ge\limsup_{n\to\infty}(\nu\widehat{\ominus}\mu)(\widehat\theta_n)$. Set $l:= \limsup_{n\to\infty}\,\nu\widehat{\ominus}\mu(\widehat\theta_n)$. For $\epsilon >0$, we consider $n_0\in\N$ such that $\sup_{k\ge n_0} \nu\widehat{\ominus}\mu(\widehat\theta_k)\leq l+\epsilon$. Then for $k\geq n_0$, $\widehat\theta_k\in\widehat{\Tc}_{l+2\epsilon}(\mu,\nu)$, and therefore $\underline{\widehat\theta}_\infty = \liminf_{k\geq n_0}\widehat\theta_k\in \widehat{\Tc}_{l+2\epsilon}(\mu,\nu)$, implying $\nu\widehat{\ominus}\mu(\widehat\theta)\leq l+2\epsilon\longrightarrow l$, as $\epsilon\to 0$. Finally, $\liminf_{n\to\infty}(\nu\widehat{\ominus}\mu)(\theta_n)\ge \nu\widehat{\ominus}\mu(\underline{\widehat\theta}_\infty)$.
\ep

\subsection{The dual irreducible convex paving}

Our final ingredient is the following measurement of subsets $K\subset\R^d$: 
 \begin{equation}\label{eq:defG}
 G(K) := \dim(K) + g_{K}(K)
 ~\mbox{where}~
 g_K(dx) 
 :=
  \frac{e^{-\frac12|x|^2}}{(2\pi)^{\frac12\dim K}}
  \,\lambda_K(dx),
 \end{equation}
Notice that $0 \leq G\leq d+1$ and, for any convex subsets $C_1\subset C_2$ of $\R^d$, we have 
\be\label{eq:G}
G(C_1) = G(C_2)
&\mbox{iff}&
\ri C_1 =\ri C_2
~~\mbox{iff}~~~
\cl C_1 = \cl C_2.
\ee
For $\theta\in\Leb_+^0(\Omega),A\in\Bc(\R^d)$, we introduce the following map from $\R^d$ to the set $\Kcirc$ of all relatively open convex subsets of $\R^d$:
 \be\label{KthetaA}
 K_{\theta,A}(x)
 :=
 \rf_x\conv(\dom\theta(x,\cdot)\setminus A)
 =\dom_X(\theta+\infty 1_{\R^d\x A}),
 &\text{for all}&
 x\in\R^d.
 \ee
We recall that a function is universally measurable if it is measurable with respect to every complete probability measure that measures all Borel subsets.

\begin{Lemma}\label{lemma:measurabilityI}
For $\theta\in\Leb_+^0(\Omega)$ and $A\in\Bc(\R^d)$, we have:
\\
{\rm (i)} $\cl\,\conv\,\dom\theta(X,\cdot):\R^d\longmapsto \Kcirc$, $\dom_X\theta:\R^d\longmapsto \ri\Kcirc$, and $K_{\theta,A}:\R^d\longmapsto \ri\Kcirc$ are universally measurable;
\\
{\rm (ii)} $G:\Kcirc\longrightarrow\R$ is Borel measurable;
\\
{\rm (iii)} if $A\in\Nc_\nu$, and $\theta\in\widehat{\Tc}(\mu,\nu)$, then up to a modification on a $\mu-$null set, $K_{\theta,A}(\R^d)\subset {\ri\,\Kcirc}$ is a partition of $\R^d$ with $x\in K_{\theta,A}(x)$ for all $x\in\R^d$.
\end{Lemma}

The proof is reported in Subsections \ref{subsect:partition} for (iii), \ref{subsect:measurabilityG} for (ii), and \ref{subsect:measurabilityI} for (i). The following property is a key-ingredient for our dual decomposition in irreducible convex paving.

\begin{Proposition}\label{prop:decompo0}
For all $(\theta,N_\nu)\in\widehat{\Tc}(\mu,\nu)\times\Nc_\nu$, we have the inclusion $Y\in \cl K_{\theta,N_\nu}(X)$, $\Mc(\mu,\nu)-$q.s.
\end{Proposition}

\proof
For an arbitrary $\P\in \Mc(\mu,\nu)$, we have by Proposition \ref{prop:numoinsmu} that $\P[\theta]<\infty$. Then, $\P\left[\dom\theta\setminus(\R^d\times N_\nu)\right] = 1$ i.e. $\P[Y\in D_X] = 1$ where $D_x := \conv(\dom\theta(x,\cdot)\setminus N_\nu)$. By the martingale property of $\P$, we deduce that
 \b*
 X &=& \E^\P[Y\mathbf{1}_{Y\in D_X}|X]
 \;=\;
 (1-\Lambda)E_{K} + \Lambda E_D, 
 ~~\mu-\mbox{a.s.}
 \e*
Where $\Lambda := \P_X[Y\in D_X\setminus\cl {K_{\theta,N_\nu}}(X)]$, $E_D := \E^{P_X}[Y|Y\in D_X\setminus\cl {K_{\theta,N_\nu}}(X)]$, $E_{K} := \E^{\P_X}[Y|Y\in \cl {K_{\theta,N_\nu}}(X)]$, and $\P_X$ is the conditional kernel to $X$ of $\P$. We have $E_{K}\in \cl  \rf_XD_X\subset D_X$ and $E_D\in D_X\setminus \cl \rf_XD_X$ because of the convexity of $D_X\setminus \cl \rf_XD_X$ given by Proposition \ref{prop:ri} (ii) ($D_X$ is convex). The lemma also gives that if $\Lambda\neq 0$, then $\E^\P[Y|X]=\Lambda E_D + (1-\Lambda) E_{K}\in D_X\setminus \cl  {K_{\theta,N_\nu}}(X)$. This implies that
 \b*
 \{\Lambda\neq 0\}
 \;\subset\;
 \{\E^\P[Y|X]\in D_X\setminus \cl  {K_{\theta,N_\nu}}(X)\}
 &\subset& 
 \{\E^\P[Y|X]\notin {K_{\theta,N_\nu}}(X)\}\\
 &\subset& 
 \{\E^\P[Y|X]\neq X\}.
 \e*
Then $\P[\Lambda \neq 0] = 0$, and therefore $\P\left[Y\in D_X\setminus\cl {K_{\theta,N_\nu}}(X)\right] = 0$. Since $\P[Y\in D_X] = 1$, this shows that $\P[Y\in \cl {K_{\theta,N_\nu}}(X)] = 1.$
\ep\\

In view of Proposition \ref{prop:decompo0} and Lemma \ref{lemma:measurabilityI} (iii), we introduce the following optimization problem which will generate our irreducible convex paving decomposition: 
 \be\label{eq:irreducible-optim}
 \underset{(\theta,N_\nu)\in\widehat{\Tc}(\mu,\nu)\x\Nc_\nu}{\inf}\mu[G(K_{\theta,N_\nu})].
 \ee
The following result gives another possible definition for the irreducible paving.

\begin{Proposition}\label{prop:icpdual}
{\rm (i)} We may find a $\mu$-a.s. unique universally measurable minimizer $\widehat{K}:=K_{\widehat\theta,\widehat N_\nu}:\R^d\to \Kcirc$ of \eqref{eq:irreducible-optim}, for some $(\widehat\theta,\widehat N_\nu)\in\widehat{\Tc}(\mu,\nu)\x\Nc_\nu$; 
\\
{\rm (ii)} for all $\theta\in\widehat{\Tc}(\mu,\nu)$ and $N_\nu\in\Nc_\nu$, we have $\widehat{K}(X)\subset K_{\theta,N_\nu}(X)$, $\mu$-a.s;
\\
{\rm (iii)} we have the equality $\widehat{K}(X) = I(X)$, $\mu-$a.s.
\end{Proposition}

In item (i), the measurability of $I$ is induced by Lemma \ref{lemma:measurabilityI} (i). Existence and uniqueness, together with (ii), are proved in Subsection \ref{subsect:minsupport}. finally, the proof of (iii) is reported in Subsection \ref{subsect:probamax}, and is a consequence of Theorem \ref{thm:CharacPolar} below. Proposition \ref{prop:icpdual} provides a characterization of the irreducible convex paving by means of an optimality criterion on $\big(\widehat{\Tc}(\mu,\nu),\Nc_\nu\big)$.

\begin{Remark}\label{rem:Compodependdual}{\rm
We illustrate how to get the components from optimization Problem \eqref{eq:irreducible-optim} in the case of Example \ref{example:Compodepend}. A $\widehat{\Tc}(\mu,\nu)$ function minimizing this problem (with $N_\nu:=\emptyset\in\Nc_\nu$) is $\widehat\theta :=\liminf_{n \to \infty}\Tbf_{p_n}f_n$, where $f_n := n f$, $p_n:=n p$ for some $p\in\partial f$, and
$$f(x):=\dist\big(x,\aff(y_1,y_{-1})\big)+\dist\big(x,\aff(y_1,y_2)\big)+\dist\big(x,\aff(y_2,y_{-1})\big).$$
One can easily check that $\mu[f] = \nu[f]$ for any $n\ge 1$: $f,f_n\in\Cfrak_0$. These functions separate $I(x_0)$, $I(x_1)$ and $\big(I(x_0)\cup I(x_1)\big)^c$.

Notice that in this example, we may as well take $\theta := 0$, and $N_\nu := \{y_{-1},y_0,y_1,y_2\}^c$, which minimizes the optimization problem as well.}
\end{Remark}

\subsection{Structure of polar sets}

Let $\theta\in\widehat{\Tc}(\mu,\nu)$, we denote the set valued map $J_\theta(X):=\dom\,\theta(X,\cdot)\cap \bar J(X)$, where $\bar J$ is introduced in Proposition \ref{prop:CharacProb}.

\begin{Remark}\label{rem:underlineJ}
{\rm
Let $\theta\in\widehat{\Tc}(\mu,\nu)$, up to a modification on a $\mu-$null set, we have
 \be\label{eq:Jtheta}
Y\in J_\theta(X),~\Mc(\mu,\nu)-\mbox{q.s,}~~
 \underline{J} 
 \subset 
 J_\theta\subset \bar J,\\
 \mbox{and}~~J_\theta~\mbox{constant on}~I(x),
 ~~\mbox{for all}~~
 x\in\R^d.
 \ee
These claims are a consequence of Proposition \ref{prop:Jqs} together with Lemma \ref{lemma:JisK}.}
\end{Remark}

Our second main result shows the importance of these set-valued maps:

\begin{Theorem}\label{thm:CharacPolar}
A Borel set $N\in\Bc(\Omega)$ is $\Mc(\mu,\nu)-$polar if and only if
$$
N \subset
 \{X\in N_\mu\}\cup\{Y\in N_\nu\}\cup\{Y\notin J_\theta(X)\},$$
for some $(N_\mu,N_\nu)\in\Nc_\mu\x\Nc_\nu$ and $\theta\in\widehat{\Tc}(\mu,\nu)$.
\end{Theorem}

The proof is reported in Section \ref{subsect:probamax}. This Theorem is an extension of the one-dimensional characterization of polar sets given by Theorem 3.2 in \cite{beiglbock2015complete}, indeed in dimension one $\underline{J} = J_\theta = \bar J$ by Proposition \ref{prop:IvsIBJ}, together with the inclusion in Remark \ref{rem:underlineJ}.

We conclude this section by reporting a duality result which will be used for the proof of Theorem \ref{thm:CharacPolar}. We emphasize that the primal objective of the accompanying paper De March \cite{HDM2017} is to push further this duality result so as to be suitable for the robust superhedging problem in financial mathematics. 

Let $c:\R^d\times\R^d\longrightarrow\R_+$, and consider the martingale optimal transport problem:
 \be
 \Sbf_{\mu,\nu}(c)
 &:=&
 \sup_{\P\in\Mc(\mu,\nu)}
 \P[c].
 \ee

Notice from Proposition \ref{prop:numoinsmu} (i) that $\Sbf_{\mu,\nu}(\theta) \le \nu\widehat{\ominus}\mu(\theta)$ for all $\theta\in\widehat{\Tc}$. We denote by $\Dc_{\mu,\nu}^{mod}(c)$ the collection of all $(\varphi,\psi,h,\theta)$ in $\L^1_+(\mu)\x\L^1_+(\nu)\x\L^0(\R^d,\R^d)\x\widehat{\Tc}(\mu,\nu)$ such that
\b*
\Sbf_{\mu,\nu}(\theta) = \nu\widehat{\ominus}\mu(\theta),
&\mbox{and}&
\varphi\oplus\psi+h^\otimes+\theta
\;\ge\; c, 
~\mbox{on}~
\{Y\in \aff K_{\theta,\{\psi = \infty\}}(X)\}.
\e*
The last inequality is an instance of the so-called robust superhedging property. The dual problem is defined by:
\b*
\Ibf_{\mu,\nu}^{mod}(c)
&:=&
\underset{(\varphi,\psi,h,\theta)\in\Dc_{\mu,\nu}^{mod}(c)}{\inf}\;\mu[\varphi]+\nu[\psi]+\nu\widehat{\ominus}\mu(\theta).
\e*
Notice that for any measurable function $c:\Omega\longrightarrow\R_+$, any $\P\in\Mc(\mu,\nu)$, and any $(\varphi,\psi,h,\theta)\in \Dc_{\mu,\nu}^{mod}(c)$, we have $\P[c]\le \mu[\varphi]+\nu[\psi]+\P[\theta]\le \mu[\varphi]+\nu[\psi]+\Sbf_{\mu,\nu}(\theta)$, as a consequence of the above robust superhedging inequality, together with the fact that $Y\in \aff K_{\theta,\{\psi = \infty\}}(X)$, $\Mc(\mu,\nu)$-q.s. by Proposition \ref{prop:decompo0} This provides the weak duality:
 \be\label{eq:weakduality}
 \Sbf_{\mu,\nu}(c)&\le &\Ibf^{mod}_{\mu,\nu}(c).
 \ee
The following result states that the strong duality holds for upper semianalytic functions. We recall that a function $f:\R^d\to\R$ is upper semianalytic if $\{f\ge a\}$ is an analytic set for any $a\in\R$. In particular, a Borel function is upper semianalytic.

\begin{Theorem}\label{thm:probduality}
Let $c:\Omega\to\overline{\R}_+$ be upper semianalytic. Then we have
\\
{\rm (i)} $\Sbf_{\mu,\nu}(c) = \Ibf_{\mu,\nu}^{mod}(c)$;
\\
{\rm (ii)} If in addition $\Sbf_{\mu,\nu}(c)<\infty$, then existence holds for the dual problem $\Ibf_{\mu,\nu}^{mod}(c)$.
\end{Theorem}

\begin{Remark}
By allowing $h$ to be infinite in some directions, orthogonal to $\aff K_{\theta,\{\psi = \infty\}}(X)$, together with the convention $\infty-\infty = \infty$, we may reformulate the robust superhedging inequality in the dual set as $\varphi\oplus\psi+h^\otimes+\theta\ge c$ pointwise.
\end{Remark}

\subsection{One-dimensional tangent convex functions}

For an interval $J\subset \R$, we denote $\Cfrak(K)$ the set of convex functions on $K$. 

\begin{Proposition}\label{prop:Setsindimone}
Let $d=1$, then
 \b*
 \widehat{\Tc}(\mu,\nu) 
 = 
 \Big\{\sum_k\mathbf{1}_{\{X\in I_k\}}\Tbf_{p_k} f_k:
        f_k\in\Cfrak(J_k),~p_k\in\partial f_k,~\sum_k(\nu_k-\mu_k)[f_k]<  \infty\Big\},
 \e*
 $\Mc(\mu,\nu)-$q.s. Furthermore, for all such $\theta\in\widehat{\Tc}(\mu,\nu)$ and its corresponding $(f_k)_k$, we have
 $
 \nu\widehat{\ominus}\mu(\theta) 
 = 
 \sum_k(\nu_k-\mu_k)[f_k]
 $.
\end{Proposition}
\proof
As all functions we consider are null on the diagonal, equality on $\cup_kI_k\x J_k$ implies $\Mc(\mu,\nu)-$q.s. equality by Theorem 3.2 in \cite{beiglbock2015complete}. Let $\Lc$ be the set on the right hand side. 

\no\underline{\rm Step 1:} We first show $\subset$, for $a\ge 0$, we denote $\Lc_a:=\{\theta\in\Lc:\sum_k(\nu_k-\mu_k)[f_k]\le a\}$. Notice that $\Lc_a$ contains $\Tbf(\Cfrak_a)$ modulo $\Mc(\mu,\nu)-$q.s. equality. We intend to prove that $\Lc_a$ is $\muxpw-$Fatou closed, so as to conclude that $\widehat\Tc_a\subset\Lc_a$, and therefore $\widehat\Tc(\mu,\nu)\subset\Lc$ by the arbitrariness of $a\ge 0$. 

Let $\theta_n=\sum_k\mathbf{1}_{\{X\in I_k\}}\Tbf_{p_{k^n}} f^n_k\in\Lc_a$ converging $\muxpw$. By Proposition \ref{prop:convergenceprop}, $\theta_n\longrightarrow\theta:=\underline\theta_\infty$, $\muxpw$. For $k\ge 1$, let $x_k\in I_k$ be such that $\theta_n(x_k,\cdot)\longrightarrow\theta(x_k,\cdot)$ on $\dom_{x_k}\theta$, and set $f_k:=\theta(x_k,\cdot)$. By Proposition 5.5 in \cite{beiglbock2015complete}, $f_k$ is convex on $I_k$, finite on $J_k$, and we may find $p_k\in\partial f_k$ such that for $x\in I_k$, $\theta(x,\cdot)=\Tbf_{p_k}f_k(x,\cdot)$. Hence, $\theta=\sum_k\mathbf{1}_{\{X\in I_k\}}\Tbf_{p_k} f_k$, and $\sum_k(\nu_k-\mu_k)[f_k]\le a$ by Fatou's Lemma, implying that $\theta\in \Lc_a$, as required.

\no\underline{\rm Step 2:} To prove the reverse inclusion $\supset$, let $\theta=\sum_k\mathbf{1}_{\{X\in I_k\}}\Tbf_{p_{k}} f_k\in\Lc$. Let $f_k^\epsilon$ be a convex function defined by $f_k^\epsilon:=f_k$ on $J_k^\epsilon = J_k\cap\left\{x\in  J_k:\dist\left(x, J^c_k\right)\ge\epsilon\right\}$, and $f_k^\epsilon$ affine on $\R\setminus J_k^\epsilon$. Set $\epsilon_n:=n^{-1}$, $\bar f_n = \sum_{k=1}^nf_k^{\epsilon_n}$, and define the corresponding subgradient in $\partial\bar f_n$:
 \b*
 \bar p_n
 :=
 p_k+\nabla(\bar f_n-f^{\eps_n}_k)
 ~~\mbox{on}~~J_k^{\eps_n},~k\ge 1,
 &\mbox{and}&
 \bar p_n := \nabla \bar f_n
 ~~\mbox{on}~~\R\setminus \big(\cup_k J_k^{\eps_n}\big).
 \e*
We have $(\nu-\mu)[\bar f_n]=\sum_{k=1}^n(\nu_k-\mu_k)[f_k^{\epsilon_n}]\leq \sum_k(\nu_k-\mu_k)[f_k]<\infty$. By definition, we see that $\Tbf_{\bar p_n}\bar f_n$ converges to $\theta$ pointwise on $\cup_k(I_k)^2$ and to $\theta_*(x,y):=\liminf_{\yb\to y}\theta(x,\yb)$ on $\cup_kI_k\times\cl I_k$ where, using the convention $\infty-\infty =\infty$, $\theta':=\theta-\theta_*\ge 0$, and $\theta'=0$ on $\cup_k(I_k)^2$. For $k\ge 1$, set $\Delta^l_k:=\theta'(x_k,l_k)$, and $\Delta^r_k:=\theta'(x_k,l_k)$ where $I_k = (l_k,r_k)$, and we fix some $x_k\in I_k$. For positive $\epsilon<\frac{r_k-l_k}{2}$, and $M\ge 0$, consider the piecewise affine function $g_k^{\epsilon,M}$ with break points $l_k+\epsilon$ and $r_k-\epsilon$, and:
 $$
 g_k^{\epsilon,M}(l_k) = M\wedge\Delta_k^l,
 ~~
 g_k^{\epsilon,M}(r_k) = M\wedge\Delta_k^r,  ~~~
 g_k^{\epsilon,M}(l_k+\epsilon) = 0,
 ~\mbox{and}~
 g_k^{\epsilon,M}(r_k-\epsilon) = 0.
 $$
Notice that $g_k^{\epsilon,M}$ is convex, and converges pointwise to $g_k^M:=M\wedge\theta'\left(\frac{l_k+r_k}{2},\cdot\right)$ on $J_k$, as $\epsilon\to 0$, with
 \b*
 (\nu_k-\mu_k)(g_k^M) 
 &=& 
 \nu_k\left[\{l_k\}\right](M\wedge\Delta_k^l)
 +\nu_k[r_k](M\wedge\Delta_k^r)
 \\
 &\le&
(\nu_k-\mu_k)[f_k]-(\nu_k-\mu_k)[(f_k)_*]
 \;\le\;
 (\nu_k-\mu_k)[f_k],
 \e*
where $(f_k)_*$ is the lower semi-continuous envelop of $f_k$. Then by the dominated convergence theorem, we may find positive $\epsilon_k^{n,M}<\frac{r_k-l_k}{2n}$ such that 
\b*
(\nu_k-\mu_k)(g_k^{\epsilon_k^{n,M},M}) \le (\nu_k-\mu_k)[f_k] +2^{-k}/n.
\e*
Now let $\bar g_n = \sum_{k=1}^ng_k^{\epsilon_k^{n,n},n}$, and $\bar p_n'\in\partial\bar g_n$. Notice that $\Tbf_{\bar p_n'}g_n\longrightarrow\theta'$ pointwise on $\cup_kI_k\times J_k$, furthermore, $(\nu-\mu)(\bar g_n)\le \sum_k(\nu_k-\mu_k)[f_k]+1/n\le \sum_k(\nu_k-\mu_k)[f_k]+1<\infty$.

Then we have $\theta_n:=\Tbf_{\bar p_n}\bar f_n+\Tbf_{\bar p_n'}\bar g_n$ converges to $\theta$ pointwise on $\cup_kI_k\times J_k$, and therefore $\Mc(\mu,\nu)-$q.s. by Theorem 3.2 in \cite{beiglbock2015complete}. Since $(\nu-\mu)(\bar f_n+\bar g_n)$ is bounded, we see that $(\theta_n)_{n\ge 1}\subset \Tbf(\Cfrak_a)$ for some $a\ge 0$. Notice that $\theta_n$ may fail to converge $\muxpw$. However, we may use Proposition \ref{prop:Komlos} to get a sequence $\widehat\theta_n\in \conv(\theta_k,k\ge n)$, and $\widehat\theta_\infty\in\Theta_\mu$ such that $\widehat\theta_n \longrightarrow\widehat\theta_\infty$, $\muxpw$ as $n\to\infty$, and satisfies the same $\Mc(\mu,\nu)-$q.s. convergence properties than $\theta_n$. Then $\underline{\widehat\theta}_\infty\in\widehat{\Tc}(\mu,\nu)$, and $\underline{\widehat\theta}_\infty = \theta$, $\Mc(\mu,\nu)-$q.s. 
\ep


\section{The irreducible convex paving}
\setcounter{equation}{0}

\subsection{Existence and uniqueness}\label{subsect:minsupport}

\no{\bf Proof of Theorem \ref{thm:icp}}
(i) The measurability follows from Lemma \ref{lemma:measurabilityI}. We first prove the existence of a minimizer for the problem \eqref{eq:irreducible-optim}. Let $m$ denote the infimum in \eqref{eq:irreducible-optim}, and consider a minimizing sequence $(\theta_n,N_\nu^n)_ {n\in\N}\subset\widehat{\Tc}(\mu,\nu)\x \Nc_\nu$ with $\mu[G(K_{\theta_n,N^n_\nu})]\leq m+1/n$. By possibly normalizing the functions $\theta_n$, we may assume that $\nu\widehat{\ominus}\mu(\theta_n)\leq 1$. Set
 \b*
 &\widehat\theta 
 \;:=\; 
 \sum_{n\ge 1} 2^{-n}\theta_n
 \quad\text{ and }\quad 
 \widehat N_\nu :=\cup_{n\ge 1} N_\nu^n\in\Nc_\nu.&
 \e*
Notice that $\widehat\theta$ is well-defined as the pointwise limit of a sequence of the nonnegative functions $\widehat\theta_N:=\sum_{n\le N}2^{-n}\theta_n$. Since $\nu\widehat{\ominus}\mu\big[\widehat\theta_N\big]\le \sum_{n\ge 1}2^{-n}<\infty$ by convexity of $\nu\widehat{\ominus}\mu$, $\widehat\theta_N\longrightarrow\widehat\theta$, pointwise, and $\widehat\theta\in\widehat{\Tc}(\mu,\nu)$ by Lemma \ref{lemma:KomlosLiminf}, since any convex extraction of $(\theta_n)_{n\ge 1}$ still converges to $\widehat\theta$. Since $\theta^{-1}_n(\{\infty\})\subset \widehat\theta^{-1}(\{\infty\})$, it follows from the definition of $\widehat N_\nu$ that
 $
 m+1/n
 \ge
 \mu[G(K_{\theta_n,N^n_\nu})]
 \;\ge\;
 \mu[G(K_{\widehat\theta,\widehat N_\nu})]
 $, hence $\mu[G(K_{\widehat\theta,\widehat N_\nu})]= m$ as $\widehat\theta\in\widehat{\Tc}(\mu,\nu)$, $\widehat N_\nu\in\Nc_\nu$.

\no(ii) For an arbitrary $(\theta,N_\nu)\in \widehat{\Tc}(\mu,\nu)\times\Nc_\nu$, we define $\bar\theta:=\theta+\widehat\theta\in\widehat{\Tc}(\mu,\nu)$ and $\bar N_\nu:=\widehat N_\nu\cup N_\nu$, so that $K_{\bar\theta,\bar N_\nu}\subset K_{\widehat\theta,\widehat N_\nu}$. By the non-negativity of $\theta$ and $\widehat\theta$, we have $m\leq\mu[G\big(K_{\bar\theta,\bar N_\nu}\big)]\leq\mu[G\big(K_{\widehat\theta,\widehat N_\nu}\big)] = m$. Then $G\big(K_{\bar\theta,\bar N_\nu}\big)=G\big(K_{\widehat\theta,\widehat N_\nu}\big)$, $\mu$-a.s. By \eqref{eq:G}, we see that, $\mu-$a.s. $K_{\bar\theta,\bar N_\nu}=K_{\widehat\theta,\widehat N_\nu}$ and $K_{\bar\theta,\bar N_\nu}=K_{\widehat\theta,\widehat N_\nu}=I$. This shows that $I\subset K_{\theta,N_\nu}$, $\mu$-a.s.
\ep

\subsection{Partition of the space in convex components}
\label{subsect:partition}

This section is dedicated to the proof of Lemma \ref{lemma:measurabilityI} (iii), which is an immediate consequence of Proposition \ref{prop:partition} (ii).

\begin{Proposition}\label{prop:partition}
Let $\theta\in\widehat{\Tc}(\mu,\nu)$, and $A\in\Bc(\R^d)$. We may find $N_\mu\in\Nc_\mu$ such that:
\\
{\rm (i)} for all $x_1,x_2\notin N_\mu$ with $K_{\theta,A}(x_1)\cap K_{\theta,A}(x_2)\neq\emptyset$, we have $K_{\theta,A}(x_1)=K_{\theta,A}(x_2);$
\\
{\rm (ii)} if $A\in\Nc_\nu$, then $x\in  K_{\theta,A}(x)$ for $x\notin N_\mu$, and up to a modification of $K_{\theta,A}$ on $N_\mu$, $K_{\theta,A}(\R^d)$ is a partition of $\R^d$ such that $x\in K_{\theta,A}(x)$ for all $x\in\R^d$.
\end{Proposition}

\proof
\no\underline{\rm (i)} Let $N_\mu$ be the $\mu-$null set given by Proposition \ref{prop:ConvexityInvariance} for $\theta$. For $x_1,x_2\notin N_\mu$, we suppose that we may find $\yb\in K_{\theta,A}(x_1)\cap K_{\theta,A}(x_2)$. Consider $y\in \cl K_{\theta,A}(x_1)$. As $K_{\theta,A}(x_1)$ is open in its affine span, $y':= \yb + \frac{\epsilon}{1-\epsilon}(\yb-y)\in K_{\theta,A}(x_1)$ for $0<\epsilon<1$ small enough. Then $\yb = \epsilon y+(1-\epsilon)y'$, and by Proposition \ref{prop:ConvexityInvariance}, we get
$$\epsilon\theta(x_1,y)+(1-\epsilon)\theta(x_1,y')-\theta(x_1,\yb)=\epsilon\theta(x_2,y)+(1-\epsilon)\theta(x_2,y')-\theta(x_2,\yb)$$
By convexity of $\dom_{x_i}\theta$, $K_{\theta,A}(x_i)\subset\dom_{x_i}\theta\subset \dom\theta(x_i,\cdot)$. Then $\theta(x_1,y')$, $\theta(x_1,\yb)$, $\theta(x_2,y')$, and $\theta(x_2,\yb)$ are finite and
 \b*
 \theta(x_1,y)<\infty 
 &\mbox{if and only if}&
 \theta(x_2,y)<\infty.
 \e*
Therefore $\cl K_{\theta,A}(x_1)\cap\dom\theta(x_1,\cdot)\subset \dom\theta(x_2,\cdot)$. We have obviously $\cl K_{\theta,A}(x_2)\cap\dom\theta(x_2,\cdot)\subset \dom\theta(x_2,\cdot)$ as well. Subtracting $A$, we get
$$\big(\cl K_{\theta,A}(x_1)\cap\dom\theta(x_1,\cdot)\setminus A\big)\cup\big(\cl K_{\theta,A}(x_2)\cap\dom\theta(x_2,\cdot)\setminus A\big)\subset \dom\theta(x_2,\cdot)\setminus A.$$
Taking the convex hull and using the fact that the relative face of a set is included in itself, we see that $\conv\big(K_{\theta,A}(x_1)\cup K_{\theta,A}(x_2)\big)\subset \conv\big(\dom\theta(x_2,\cdot)\setminus A\big)$. Notice that, as $K_{\theta,A}(x_2)$ is defined as the $x_2-$relative face of some set, either $x_2\in \ri K_{\theta,A}(x)$ or $K_{\theta,A}(x)=\emptyset$ by the properties of $\rf_{x_2}$. The second case is excluded as we assumed that $K_{\theta,A}(x_1)\cap K_{\theta,A}(x_2)\neq\emptyset$. Therefore, as $K_{\theta,A}(x_1)$ and $K_{\theta,A}(x_2)$ are convex sets intersecting in relative interior points and $x_2\in\ri K_{\theta,A}(x_2)$, it follows from Lemma \ref{lemma:intersectinterior} that $x_2\in\ri\,\conv\big(K_{\theta,A}(x_1)\cup K_{\theta,A}(x_2)\big)$. Then by Proposition \ref{prop:ri} (ii),
 \b*
 \rf_{x_2}\conv\big(K_{\theta,A}(x_1)\cup K_{\theta,A}(x_2)\big)
 &=&
 \ri\,\conv\big(K_{\theta,A}(x_1)\cup K_{\theta,A}(x_2)\big) 
 \\
 &=&
 \conv\big(K_{\theta,A}(x_1)\cup K_{\theta,A}(x_2)\big).
 \e*
Then, we have $\conv\big(K_{\theta,A}(x_1)\cup K_{\theta,A}(x_2)\big)\subset \rf_{x_2}\conv\big(\dom\theta(x_2,\cdot)\setminus A\big) = K_{\theta,A}(x_2)$, as $\rf_{x_2}$ is increasing. Therefore $K_{\theta,A}(x_1)\subset K_{\theta,A}(x_2)$ and by symmetry between $x_1$ and $x_2$, $K_{\theta,A}(x_1)= K_{\theta,A}(x_2)$.

\no\underline{\rm (ii)} We suppose that $A\in\Nc_\nu$. First, notice that, as $K_{\theta,A}(X)$ is defined as the $X-$relative face of some set, either $x\in K_{\theta,A}(x)$ or $K_{\theta,A}(x)=\emptyset$ for $x\in\R^d$ by the properties of $\rf_x$. Consider $\P\in\Mc(\mu,\nu)$. By Proposition \ref{prop:decompo0}, $\P[Y\in \cl K_{\theta,A}(X)]=1$. As $\supp(\P_X)\subset \cl K_{\theta,A}(X)$, $\mu$-a.s., $K_{\theta,A}(X)$ is non-empty, which implies that $x\in K_{\theta,A}(x)$. Hence, $\{X\in K_{\theta,A}(X)\}$ holds outside the set $N_\mu^{0}:=\{\supp(\P_X)\not\subset \cl I(X)\}\in\Nc_\mu$. Then we just need to have this property to replace $N_\mu$ by $N_\mu\cup N_\mu^0\in\Nc_\mu$.

Finally, to get a partition of $\R^d$, we just need to redefine $K_{\theta,A}$ on $N_\mu$. If $x\in\underset{x'\notin N_\mu}{\bigcup}K_{\theta,A}(x')$ then by definition of $N_\mu$, the set $K_{\theta,A}(x')$ is independent of the choice of $x'\notin N_\mu$ such that $x\in K_{\theta,A}(x')$: indeed, if $x_1',x_2'\notin N_\mu$ satisfy $x\in K_{\theta,A}(x_1')\cap K_{\theta,A}(x_2')$, then in particular $K_{\theta,A}(x_1')\cap K_{\theta,A}(x_2')$ is non-empty, and therefore $K_{\theta,A}(x_1')= K_{\theta,A}(x_2')$ by (i). We set $K_{\theta,A}(x):=K_{\theta,A}(x')$. Otherwise, if $x\notin\underset{x'\notin N_\mu}{\bigcup}K_{\theta,A}(x')$, we set $K_{\theta,A}(x):=\{x\}$ which is trivially convex and relatively open. With this definition, $K_{\theta,A}(\R^d)$ is a partition of $\R^d$.
\ep

\section{Proof of the duality}
\setcounter{equation}{0}

For simplicity, we denote $\mbox{\rm Val}(\xi):=\mu[\varphi]+\nu[\psi]+\nu\widehat{\ominus}\mu(\theta)$, for $\xi:=(\varphi,\psi,h,\theta)\in\Dc^{mod}_{\mu,\nu}(c)$.

\subsection{Existence of a dual optimizer}

\begin{Lemma}\label{lemma:compactnessprob}
Let $c,c_n:\Omega\longrightarrow \overline{\R_+}$, and $\xi_n\in \Dc^{mod}_{\mu,\nu}(c_n)$, $n\in\N$, be such that
$$
c_n \longrightarrow c,~\mbox{pointwise, and}~
\mbox{\rm Val}(\xi_n)\longrightarrow \Sbf_{\mu,\nu}(c)<\infty
~\mbox{as}~n\to\infty.
$$
Then there exists $\xi\in \Dc^{mod}_{\mu,\nu}(c)$ such that
$\mbox{\rm Val}(\xi_n)\longrightarrow\mbox{\rm Val}(\xi)$ as $n\to\infty$.
\end{Lemma}
\proof
Denote
$\xi_n := (\varphi_n,\psi_n, h_n,\theta_n)$, and observe that the convergence of $\mbox{\rm Val}(\xi_n)$ implies that the sequence $\big(\mu(\varphi_n),\nu(\psi_n),\nu\widehat{\ominus}\mu(\theta_n)\big)_n$ is bounded, by the non-negativity of $\varphi_n,\psi_n$ and $\nu\widehat{\ominus}\mu(\theta_n)$. We also recall the robust superhedging inequality
 \be\label{eq:superhedging-n}
 \varphi_n\oplus\psi_n + h_n^\otimes + \theta_n \ge c_n,
 &\text{on }\{Y\in \aff K_{\theta_n,\{\psi_n = \infty\}}(X)\},&n\ge 1.
 \ee

\no \underline{\it Step 1.} By Koml\`os Lemma together with Lemma \ref{lemma:KomlosLiminf}, we may find a sequence $(\widehat\varphi_n,\widehat\psi_n,\widehat\theta_n)\in\conv\{(\varphi_k,\psi_k,\theta_k),k\ge n\}$ such that 
 $$
 \begin{array}{c}
 \widehat\varphi_n\longrightarrow\varphi:=\underline{\widehat\varphi}_\infty,
 ~\mu-\mbox{a.s., }
 \widehat\psi_n\longrightarrow\psi:=\underline{\widehat\psi}_\infty,~\nu-\mbox{a.s., and}
 \\
 \widehat\theta_n\longrightarrow \widetilde\theta
 :=
 \underline{\widehat\theta}_\infty\in\widehat{\Tc}(\mu,\nu),~\mu\otimes{\rm pw}.
 \end{array}
 $$
Set $\varphi:=\infty$ and $\psi:=\infty$ on the corresponding non-convergence sets, and observe that $\mu[\varphi]+\nu[\psi]<\infty$, by the Fatou Lemma, and therefore $N_\mu:=\{\varphi=\infty\}\in\Nc_\mu$ and $N_\nu:=\{\psi=\infty\}\in\Nc_\nu$. We denote by $(\widehat h_n,\widehat c_n)$ the same convex extractions from $\{(h_k,c_k),k\ge n\}$, so that the sequence $\widehat\xi_n:=(\widehat\varphi_n,\widehat\psi_n,\widehat h_n,\widehat\theta_n)$ inherits from \eqref{eq:superhedging-n} the robust superhedging property, as for $\theta_1,\theta_2\in\widehat{\Tc}(\mu,\nu)$, $\psi_1,\psi_2\in\Leb_+^1(\R^d)$, and $0<\lambda<1$, we have $\aff K_{\lambda\theta_1+(1-\lambda)\theta_2,\{\lambda\psi_1+(1-\lambda)\psi_2=\infty\}}\subset \aff K_{\theta_1,\{\psi_1=\infty\}}\cap \aff K_{\theta_2,\{\psi_2=\infty\}}$: 
\be\label{eq:ineqfonda} 
\widehat\varphi_n\oplus\widehat\psi_n+\widehat\theta_n+ \widehat h_n^\otimes\ge \widehat c_n\ge 0,&\text{pointwise on }\aff K_{\widehat\theta_n,\{\widehat\psi_n = \infty\}}(X).
\ee

\no \underline{\it Step 2.} Next, notice that $l_n:=\left(\widehat h_n^\otimes\right)^-:=\max\left(-\widehat h_n^\otimes,0\right)\in\Theta$ for all $n\in\N$. By the convergence Proposition \ref{prop:Komlos}, we may find convex combinations $\widehat {l}_n:=\sum_{k\geq n}\lambda_k^nl_k\longrightarrow l:= \underline{\widehat l}_\infty$, $\mu\otimes{\rm pw}$. Updating the definition of $\varphi$ by setting $\varphi:=\infty$ on the zero $\mu-$measure set on which the last convergence does not hold on $(\partial^x\dom l)^c$, it follows from \eqref{eq:ineqfonda}, and the fact that $\aff K_{\bar\theta,\{\psi=\infty\}}\subset \liminf_{n\to\infty} \aff K_{\widehat\theta_n,\{\widehat\psi_n = \infty\}}$, that 
 $$
 l
 =
 \underline{\widehat l}_\infty
 \le 
 \liminf_n \sum_{k\geq n}\lambda_k^n
                         \big(\widehat\varphi_k\oplus\widehat\psi_k+\widehat\theta_k
                         \big)
 \le
 \varphi\oplus\psi+\bar\theta,
 ~\mbox{pointwise on }\left\{Y\in \aff K_{\bar\theta,\{\psi=\infty\}}(X)\right\}.
 $$
 where $\bar\theta:=\liminf_n\sum_{k\geq n}\lambda_k^n\widehat\theta_k\in\widehat{\Tc}(\mu,\nu)$. As $\{\varphi = \infty\}\in\Nc_\mu$, by possibly enlarging $N_\mu$, we assume without loss of generality that $\{\varphi = \infty\}\subset N_\mu$, we see that $\dom\,l\supset (N_\mu^c\times N_\nu^c) \cap\dom\bar\theta\cap \left\{Y\in \aff K_{\bar\theta,\{\psi=\infty\}}(X)\right\}$, and therefore
  \begin{equation}\label{doml}
  K_{\bar\theta,\{\psi=\infty\}}(X)\subset\dom_X l'\subset \dom\, l'(X,\cdot),~~
  \mu\mbox{-a.s.}
  \end{equation}
\no \underline{\it Step 3.}
Let $\widehat{\widehat{h}}_n:=\sum_{k\geq n}\lambda_k^n \widehat h_k$. Then $b_n:=\widehat{\widehat{h}}^\otimes_n+\widehat l_n =\sum_{k\geq n}\lambda_k^n \left(\widehat h_k^\otimes\right)^+$ defines a non-negative sequence in $\Theta$. By Proposition \ref{prop:Komlos}, we may find a sequence $\widehat b_n=:\widetilde h_n^\otimes+\widetilde l_n\in\conv(b_k,k\ge n)$ such that $\widehat b_n\longrightarrow b:=\underline{\widehat b}_\infty$, $\mu\otimes{\rm pw}$, where $b$ takes values in $[0,\infty]$. $\widehat b_n(X,\cdot)\longrightarrow b(X,\cdot)$ pointwise on $\dom_X b$, $\mu-$a.s. Combining with \eqref{doml}, this shows that
 \b*
 \widetilde h_n^\otimes(X,\cdot)
 \longrightarrow(b-l)(X,\cdot)
 &\mbox{pointwise on}&
 \dom_X b\cap K_{\bar\theta,\{\psi=\infty\}}(X),~~\mu-\mbox{a.s.}
 \e*
$(b-l)(X,\cdot)=\liminf_n\widetilde h_n^\otimes(X,\cdot)$, pointwise on $K_{\bar\theta,\{\psi=\infty\}}(X)$ (where $l$ is a limit of $l_n$), $\mu-$a.s. Clearly, on the last convergence set, $(b-l)(X,\cdot)>-\infty$ on $K_{\bar\theta,\{\psi=\infty\}}(X)$, and we now argue that $(b-l)(X,\cdot)<\infty$ on $K_{\bar\theta,\{\psi=\infty\}}(X)$, therefore $K_{\bar\theta,\{\psi=\infty\}}(X)\subset\dom_X b$, so that we deduce from the structure of $\widetilde h_n^\otimes$ that the last convergence holds also on $\aff K_{\bar\theta,\{\psi=\infty\}}(X)$:
 \be\label{hAff}
 \widetilde h_n^\otimes(X,\cdot)
 \longrightarrow(b-l)(X,\cdot)=:h^\otimes(X,\cdot)
 &\mbox{pointwise on}&
 K_{\bar\theta,\{\psi=\infty\}}(X),~~\mu-\mbox{a.s.}
 \ee
Indeed, let $x$ be an arbitrary point of the last convergence set, and consider an arbitrary $y\in K_{\bar\theta,\{\psi=\infty\}}(x)$. By the definition of $K_{\bar\theta,\{\psi=\infty\}}$, we have $x\in\mbox{ri} K_{\bar\theta,\{\psi=\infty\}}(x)$, and we may therefore find $y'\in K_{\bar\theta,\{\psi=\infty\}}(x)$ with $x=py+(1-p)y'$ for some $p\in(0,1)$. Then, $p\,\widetilde h_n^\otimes(x,y)+(1-p)\widetilde h_n^\otimes(x,y')=0$. Sending $n\to\infty$, by concavity of the $\liminf$, this provides $p(b-l)(x,y)+(1-p)(b-l)(x,y')\le0$, so that $(b-l)(x,y')>-\infty$ implies that $(b-l)(x,y)<\infty$.

\no \underline{\it Step 4.} Notice that by dual reflexivity of finite dimensional vector spaces, \eqref{hAff} defines a unique $h(X)$ in the vector space $\aff K_{\bar\theta,\{\psi=\infty\}}(X)-X$, such that $(b-l)(X,\cdot)=h^\otimes(X,\cdot)$ on $\aff K_{\bar\theta,\{\psi=\infty\}}(X)$. At this point, we have proceeded to a finite number of convex combinations which induce a final convex combination with coefficients $(\bar\lambda_n^k)_{k\ge n\ge 1}$. Denote $\bar\xi_n:=\sum_{k\ge n}\bar\lambda_n^k\xi_k$, and set $\theta:=\underline{\bar\theta}_\infty$. Then, applying this convex combination to the robust superhedging inequality \eqref{eq:superhedging-n}, we obtain by sending $n\to\infty$ that $(\varphi\oplus\psi+h^\otimes+\theta)(X,\cdot)\ge c(X,\cdot)$ on $\aff K_{\bar\theta,\{\psi=\infty\}}(X)$, $\mu-$a.s. and $\varphi\oplus\psi+h^\otimes+\theta=\infty$ on the complement $\mu$ null-set. As $\theta$ is the liminf of a convex extraction of $(\widehat\theta_n)$, we have $\theta\ge\underline{\widehat\theta}_\infty=\bar\theta$, and therefore $\aff K_{\theta,\{\psi=\infty\}}\subset\aff K_{\bar\theta,\{\psi=\infty\}}$. This shows that the limit point $\xi:=(\varphi,\psi,h,\theta)$ satisfies the pointwise robust superhedging inequality
 \be\label{lim-superhedging}
 \varphi\oplus\psi+\theta+ h^\otimes
 \ge 
 c,
 &\text{on }\left\{Y\in \aff K_{\theta,\{\psi=\infty\}}(X)\right\}.
 \ee
\no \underline{\it Step 5.} By Fatou's Lemma and Lemma \ref{lemma:KomlosLiminf}, we have
 \begin{equation}\label{geSmunu}
 \mu[\varphi]+\nu[\psi]+\nu\widehat{\ominus}\mu[\theta]
 \le 
 \liminf_n \mu[\varphi_n]+\nu[\psi_n]+\nu\widehat{\ominus}\mu[\theta_n]
 =
 \Sbf_{\mu,\nu}(c).
 \end{equation}
By \eqref{lim-superhedging}, we have $\mu[\varphi]+\nu[\psi]+\P[\theta]\ge \P[c]$ for all $\P\in\Mc(\mu,\nu)$. Then, $\mu[\varphi]+\nu[\psi]+\Sbf_{\mu,\nu}[\theta]\ge \Sbf_{\mu,\nu}[c]$. By Proposition \ref{prop:numoinsmu} (i), we have $\Sbf_{\mu,\nu}[\theta]\leq \nu\widehat{\ominus}\mu[\theta]$, and therefore
 $$
 \Sbf_{\mu,\nu}[c]
 \le 
 \mu[\varphi]+\nu[\psi]+\Sbf_{\mu,\nu}[\theta] 
 \le 
 \mu[\varphi]+\nu[\psi]+\nu\widehat{\ominus}\mu[\theta]
 \le
 \Sbf_{\mu,\nu}(c),
 $$
by \eqref{geSmunu}.
Then we have $\mbox{\rm Val}(\xi)=\mu[\varphi]+\nu[\psi]+\nu\widehat{\ominus}\mu[\theta] = \Sbf_{\mu,\nu}(c)$ and $\Sbf_{\mu,\nu}[\theta]=\nu\widehat{\ominus}\mu[\theta]$, so that $\xi\in\Dc^{mod}_{\mu,\nu}(c)$.
\ep

\subsection{Duality result}

We first prove the duality in the lattice $\mbox{USC}_b$ of bounded upper semicontinuous fonctions $\Omega\longrightarrow\overline{\R_+}$. This is a classical result using the Hahn-Banach Theorem, the proof is reported for completeness.

\begin{Lemma}\label{lem:dualUSC}
Let $f\in\mbox{\rm USC}_b$, then $\Sbf_{\mu,\nu}(f) = \Ibf^{mod}_{\mu,\nu}(f)$
\end{Lemma}
\proof
We have $\Sbf_{\mu,\nu}(f) \leq \Ibf^{mod}_{\mu,\nu}(f)$ by weak duality \eqref{eq:weakduality}, let us now show the converse inequality $\Sbf_{\mu,\nu}(f)\geq \Ibf^{mod}_{\mu,\nu}(f)$. By standard approximation technique, it suffices to prove the result for bounded continuous $f$. We denote by $\Ctn_l(\R^d)$ the set of continuous mappings $\R^d\to\R$ with linear growth at infinity, and by $\Ctn_b(\R^d,\R^d)$ the set of continuous bounded mappings $\R^d\longrightarrow \R^d$. Define 
 $$
 \Dc(f) 
 :=  
 \left\{ (\bar\varphi,\bar\psi,\bar h)
         \in \Ctn_l(\R^d)\x\Ctn_l(\R^d)\x\Ctn_b(\R^d,\R^d):~
         \bar\varphi\oplus\bar\psi+\bar h^\otimes\ge f 
 \right\},
 $$
and the associated $
 \Ibf_{\mu,\nu}(f) :=\inf_{(\bar\varphi,\bar\psi,\bar h)\in\Dc(f)}\mu(\bar\varphi)+\nu(\bar\psi)
$. By Theorem 2.1 in Zaev \cite{zaev2015monge}, and Lemma \ref{lemma:inclusionD} below, we have
 \b*
 \Sbf_{\mu,\nu}(f) 
 &=&
 \Ibf_{\mu,\nu}(f)
 \;=\;
 \inf_{(\bar \varphi,\bar \psi,\bar h)\in \Dc(f)}\mu(\bar \varphi)+\nu(\bar \psi)
 \;\ge\; 
 \Ibf^{mod}_{\mu,\nu}(f),
 \e*
which provides the required result.
\ep
\\

\no {\bf Proof of Theorem \ref{thm:probduality}} The existence of a dual optimizer follows from a direct application of the compactness Lemma \ref{lemma:compactnessprob} to a minimizing sequence of robust superhedging strategies.

As for the extension of duality result of Lemma \ref{lem:dualUSC} to non-negative upper semi-analytic functions, we shall use the capacitability theorem of Choquet, similar to \cite{kellerer1984duality} and \cite{beiglbock2015complete}. Let $[0,\infty]^{\Omega}$ denote the set of all nonnegative functions ${\Omega}\to [0,\infty]$, and $\mbox{USA}_+$ the sublattice of upper semianalytic functions. Note that $\mbox{USC}_b$ is stable by infimum.

Recall that a $\mbox{USC}_b$-capacity is a monotone map $\mathbf{C}:[0,\infty]^{\Omega}\longrightarrow [0,\infty]$, sequentially continuous upwards on $[0,\infty]^{\Omega}$, and sequentially continuous downwards on $\mbox{USC}_b$. The Choquet capacitability theorem states that a $\mbox{USC}_b-$capacity $\mathbf{C}$ extends to $\mbox{USA}_+$ by:
 \b*
 {\bf C}(f)
 \;=\;
 \sup\big\{{\bf C}(g): g\in\mbox{USC}_b~\mbox{and}~g\leq f\big\}
 &\mbox{for all}&
 f\in\mbox{USA}_+.
 \e*
In order to prove the required result, it suffices to verify that $\Sbf_{\mu,\nu}$ and $\Ibf^{mod}_{\mu,\nu}$ are $\mbox{USC}_b$-capacities. 
As $\Mc(\mu,\nu)$ is weakly compact, it follows from similar argument as in Prosposition 1.21, and Proposition 1.26 in Kellerer \cite{kellerer1984duality} that $\Sbf_{\mu,\nu}$ is a $\mbox{USC}_b$-capacity.
We next verify that $\Ibf^{mod}_{\mu,\nu}$ is a $\mbox{USC}_b$-capacity. Indeed, the upwards continuity is inherited from $\Sbf_{\mu,\nu}$ together with the compactness lemma \ref{lemma:compactnessprob}, and the downwards continuity follows from the downwards continuity of $\Sbf_{\mu,\nu}$ together with the duality result on $\mbox{USC}_b$ of Lemma \ref{lem:dualUSC}.
\ep
\\

\begin{Lemma}\label{lemma:inclusionD}
Let $c:\Omega\to\overline\R_+$, and $(\bar\varphi,\bar\psi,\bar h)\in\Dc(c)$. Then, we may find $\xi\in\Dc^{mod}_{\mu,\nu}(c)$ such that $\mbox{\rm Val}(\xi) = \mu[\bar\varphi]+\nu[\bar\psi]$.
\end{Lemma}
\proof
Let us consider $(\bar\varphi,\bar\psi,\bar h)\in \Dc(c)$. Then $\bar\varphi\oplus\bar\psi+\bar h^\otimes\geq c \geq 0$, and therefore
$$\bar\psi(y)\geq f(y) := \underset{x\in\R^d}{\sup}-\bar\varphi(x)-\bar h(x)\cdot(y-x).$$
Clearly, $f$ is convex, and $f(x)\geq -\bar\varphi(x)$ by taking value $x=y$ in the supremum. Hence $\bar\psi-f\geq 0$ and $\bar\varphi+f\geq 0$, implying in particular that $f$ is finite on $\R^d$. 
As $\bar\varphi$ and $\bar\psi$ have linear growth at infinity, $f$ is in $\Leb^1(\nu)\cap \Leb^1(\mu)$. We have $f\in\Cfrak_a$ for $a=\nu[f]-\mu[f]\ge 0$. Then we consider $p\in\partial f$ and denote $\theta:=\Tbf_pf$. $\theta\in\Tbf\big(\Cfrak_a\big)\subset\widehat{\Tc}(\mu,\nu)$. Then denoting $\varphi:=\bar\varphi+f$, $\psi:=\bar\psi-f$, and $ h:=\bar h+p$, we have $\xi:=(\varphi,\psi, h,\theta)\in \Dc_{\mu,\nu}^{mod}(c)$ and
$$\mu[\bar\varphi]+\nu[\bar\psi] = \mu[\varphi]+\nu[\psi] + (\nu-\mu)[f]=  \mu[\varphi]+\nu[\psi] + \nu\widehat{\ominus}\mu[\theta] = \mbox{\rm Val}(\xi).$$
\ep

\section{Polar sets and maximum support martingale plan}
\label{sect:polar}
\setcounter{equation}{0}

\subsection{Boundary of the dual paving}\label{subect:J}

Consider the optimization problems:
 \begin{equation}\label{eq:irreducible-optim1}
 \underset{(\theta,N_\nu)\in\widehat{\Tc}(\mu,\nu)\x\Nc_\nu}{\inf}\mu\big[G(R_{\theta,N_\nu})\big],
 ~\mbox{with}~
 R_{\theta,N_\nu}:= \cl\,\conv\big(\dom\theta(X,\cdot)\cap\partial \widehat{K}(X)\cap N_\nu^c\big),~~~
 \end{equation}
and for all $y\in\R^d$ we consider
 \be\label{eq:irreducible-optim3}
 \underset{(\theta,N_\nu)\in\widehat{\Tc}(\mu,\nu)\x\Nc_\nu}
 \inf\;\mu\big[y\in \partial \widehat{K}(X)
               \cap\dom\theta(X,\cdot)\cap N_\nu^c
          \big].
 \ee
These problems are well defined by the following measurability result, whose proof is reported in Subsection \ref{subsect:measurabilityI}.

\begin{Lemma}\label{lemma:measurabilityFrontier}
Let $F:\R^d\longrightarrow \Kc$, $\gamma-$measurable. Then we may find $N_\gamma\in\Nc_\gamma$ such that $\mathbf{1}_{Y\in F(X)}\mathbf{1}_{X\notin N_\gamma}$ is Borel measurable, and if $X\in \ri F(X)$ convex, $\gamma-$a.s., then $\mathbf{1}_{Y\in \partial F(X)}\mathbf{1}_{X\notin N_\gamma}$ is Borel measurable as well.
\end{Lemma}

By the same argument than that of the proof of existence and uniqueness in Proposition \ref{prop:icpdual}, we see that the problem \eqref{eq:irreducible-optim1}, (resp. \eqref{eq:irreducible-optim3} for $y\in\R^d$) has an optimizer $(\theta^*,N^*_\nu)\in\widehat{\Tc}(\mu,\nu)\x\Nc_\nu$, (resp. $(\theta^*_y,N^*_{\nu,y})\in\widehat{\Tc}(\mu,\nu)\x\Nc_\nu$). Furthermore, $D:=R_{\theta^*,N_\nu^*}$, (resp $D_y(x):=\{y\}$ if $y\in \partial \widehat{K}(x)\cap\dom\theta^*_y(x,\cdot)\cap N_{\nu,y}^*$, and $\emptyset$ otherwise, for $x\in\R^d$) does not depend on the choice of $(\theta^*,N^*_\nu)$, (resp. $\theta_y^*$) up to a $\mu-$negligible modification. 
 
We define $\bar K:=D\cup \widehat{K}$, and $K_\theta(X) := \dom\theta(X,\cdot)\cap \bar K(X)$ for $\theta\in\widehat{\Tc}(\mu,\nu)$. Notice that if $y\in\R^d$ is not an atom of $\nu$, we may chose $N_{\nu,y}$ containing $y$, which means that Problem \eqref{eq:irreducible-optim3} is non-trivial only if $y$ is an atom of $\nu$. We denote ${\rm atom}(\nu)$, the (at most countable) atoms of $\nu$, and define the mapping
$\underline{K}:= (\cup_{y\in{\rm atom}(\nu)}D_y)\cup \widehat{K}$, 

\begin{Proposition}\label{prop:Jqs}
Let $\theta\in\widehat{\Tc}(\mu,\nu)$. Up to a modification on a $\mu-$null set, we have\\
{\rm (i)} $\bar K$ is convex valued, moreover $Y\in \bar K(X)$, and $Y\in K_\theta(X)$, $\Mc(\mu,\nu)-$q.s.\\
{\rm (ii)} $\widehat{K}\subset\underline{K}\subset K_\theta\subset \bar K\subset \cl \widehat{K}$,\\
{\rm (iii)} $\underline{K}$, $K_\theta$, and $\bar K$ are constant on $\widehat{K}(x)$, for all $x\in\R^d$.
\end{Proposition}
\proof
\no\underline{\rm (i)}
For $x\in\R^d$, $\bar K(x)=D(x)\cup \widehat{K}(x)$. Let $y_1,y_2\in \bar K(x)$, $\lambda\in (0,1)$, and set $y:=\lambda y_1+(1-\lambda)y_2$. If $y_1,y_2\in \widehat{K}(x)$, or $y_1,y_2\in D(x)$, we get $y\in \bar K(x)$ by convexity of $\widehat{K}(x)$, or $D(x)$. Now, up to switching the indices, we may assume that $y_1\in \widehat{K}(x)$, and $y_2\in D(x)\setminus \widehat{K}(x)$. As $D(x)\setminus \widehat{K}(x)\subset \partial \widehat{K}(x)$, $y\in \widehat{K}(x)$, as $\lambda>0$. Then $y\in \bar K(x)$. Hence, $\bar K$ is convex valued.

Since $\dom\theta^*(X,\cdot)\setminus N_\nu^*\cap\cl \widehat{K}\setminus \widehat{K}\subset R_{\theta^*,N_\nu^*}$, we have $\dom\theta^*(X,\cdot)\setminus N_\nu^*\cap\cl \widehat{K}\subset R_{\theta^*,N_\nu^*}\cup \widehat{K} = \bar K$. Then, as $Y\in \dom\theta^*(X,\cdot)\setminus N_\nu^*$, and $Y\in\cl \widehat{K}(X)$, $Y\in \bar K(X)$, $\Mc(\mu,\nu)-$q.s.

Let $\theta\in\widehat{\Tc}(\mu,\nu)$, then $Y\in\dom\theta(X,\cdot)$, $\Mc(\mu,\nu)-$q.s. Finally we get $Y\in\dom\theta(X,\cdot)\cap \bar K(X)=K_\theta(X)$, $\Mc(\mu,\nu)-$q.s.

\no\underline{\rm (ii)}
As $R_{\theta,N_\nu}(X)\subset\cl\conv\partial \widehat{K}(X) = \cl \widehat{K}(X)$, $\bar K\subset \cl \widehat{K}$. By definition, $K_\theta\subset \bar K$, and $\widehat{K}\subset \underline{K}$. For $y\in{\rm atom}(\nu)$, and $\theta_0\in\widehat{\Tc}(\mu,\nu)$, by minimality,
\be\label{eq:trivialinclusion}
D_y(X)\subset\dom\theta_0(X,\cdot)\cap \partial \widehat{K}(X),&\mu-\mbox{a.s.}
\ee

Applying \eqref{eq:trivialinclusion} for $\theta_0 = \theta$, we get $D_y\subset \dom\theta(X,\cdot)$, and for $\theta_0 = \theta^*$, $D_y(X)\subset \bar K(X)$, $\mu-$a.s. Taking the countable union: $\underline{K}\subset K_\theta$, $\mu-$a.s. (This is the only inclusion that is not pointwise). Then we change $\underline{K}$ to $\widehat{K}$ on this set to get this inclusion pointwise.

\no\underline{\rm (iii)} For $\theta_0\in\widehat{\Tc}(\mu,\nu)$, let $N_\mu\in\Nc_\mu$ from Proposition \ref{prop:ConvexityInvariance}. Let $x\in N_\mu^c$, $y\in\partial \widehat{K}(x)$, and $y':= \frac{x+y}{2}\in \widehat{K}(x)$. Then for any other $x'\in \widehat{K}(x)\cap N_\mu^c$, $\frac12\theta_0(x,y)-\theta_0(x,y') = \frac12\theta_0(x',x)+\frac12\theta_0(x',y)-\theta_0(x',y')$, in particular, $y\in \dom\theta(x,\cdot)$ if and only if $y\in \dom\theta(x',\cdot)$. Applying this result to $\theta$, $\theta^*$, and $\theta^*_y$ for all $y\in{\rm atom}(\nu)$, we get $N_\mu$ such that for any $x\in\R^d$, $\bar K$, $K_\theta$, and $\underline{K}$ are constant on $\widehat{K}(x)\cap N_\mu^c$. To get it pointwise, we redefine these mappings to this constant value on $\widehat{K}(x)\cap N_\mu$, or to $\widehat{K}(x)$, if $\widehat{K}(x)\cap N_\mu^c = \emptyset$. The previous properties are preserved.
\ep

\subsection{Structure of polar sets}

\begin{Proposition}\label{prop:polar}
A Borel set $N\in\Bc(\Omega)$ is $\Mc(\mu,\nu)-$polar if and only if for some $(N_\mu,N_\nu)\in\Nc_\mu\x\Nc_\nu$ and $\theta\in\widehat{\Tc}(\mu,\nu)$, we have
 $$
 N \subset
 \{X\in N_\mu\}\cup\{Y\in N_\nu\}\cup\{Y\notin K_\theta(X)\}
 .
 $$
\end{Proposition}
\proof
One implication is trivial as $Y\in K_\theta(X)$, $\Mc(\mu,\nu)-$q.s. for all $\theta\in\widehat{\Tc}(\mu,\nu)$, by Proposition \ref{prop:Jqs}. We only focus on the non-trivial implication. For an $\Mc(\mu,\nu)$-polar set $N$, we have $\Sbf_{\mu,\nu}(\infty\mathbf{1}_N) = 0$, and it follows from the dual formulation of Theorem \ref{thm:probduality} that $0 = \mbox{\rm Val}(\xi)$ for some $\xi=(\varphi,\psi,h,\theta)\in \Dc_{\mu,\nu}^{mod}(\infty\mathbf{1}_N)$. Then, 
 \b*
 \varphi <\infty,~\mu-\mbox{a.s.}, 
 &\psi <\infty,~\nu-\mbox{a.s.}&\mbox{and }
 \theta\in\widehat{\Tc}(\mu,\nu),
 \e*
 
As $h$ is finite valued, and $\varphi,\psi$ are non-negative functions, the superhedging inequality $\varphi\oplus\psi+\theta+h^\otimes\ge\infty\mathbf{1}_N$ on $\{Y\in \aff K_{\theta,\{\psi = \infty\}}(X)\}$ implies that
\be\label{eq:superhedging-qs}
\mathbf{1}_{\{\varphi=\infty\}}\oplus\mathbf{1}_{\{\psi=\infty\}}+\mathbf{1}_{\{(\dom\theta)^c\}}\ge\mathbf{1}_N&\mbox{on}&\{Y\in \aff K_{\theta,\{\psi = \infty\}}(X)\}
\ee

By Proposition \ref{prop:icpdual} (ii), we have $\widehat{K}(X)\subset K_{\theta,\{\psi=\infty\}}(X)$, $\mu-$a.s. Then $\bar K(X)\subset \aff \widehat{K}(X)\subset \aff K_{\theta,\{\psi=\infty\}}(X)$, which implies that
\begin{equation}\label{eq:Jinclus}
K_\theta(X)
:=
\dom\theta(X,\cdot)\cap \bar K(X) \subset \dom\theta(X,\cdot)\cap\aff K_{\theta,\{\psi=\infty\}}(X),~\mu-\mbox{a.s.}
\end{equation}
We denote $N_\mu := \{\varphi=\infty\}\cup\{K_\theta(X)\not\subset \dom\theta(X,\cdot)\cap\aff K_{\theta,\{\psi=\infty\}}(X)\}\in\Nc_\mu$, and $N_\nu:=\{\psi = \infty\}\in\Nc_\nu$. Then by \eqref{eq:superhedging-qs}, $\mathbf{1}_N =0$ on $(\{\varphi =\infty\}^c\x \{\psi = \infty\}^c)\cap\{Y\in \dom\theta(X,\cdot)\cap \aff K_{\theta,\{\psi=\infty\}}(X)\}$, and therefore by \eqref{eq:Jinclus},
$N\subset \{X\in N_\mu\}\cup \{Y\in N_\nu\} \cup \{Y\notin K_\theta(X)\}$.
\ep

\subsection{The maximal support probability}\label{subsect:probamax}

In order to prove the existence of a maximum support martingale transport plan, we introduce the maximization problem.
 \be\label{supportmax}
 M:=\underset{\P\in\Mc(\mu,\nu)}{\sup}\mu[G(\wideparen{\supp}\P_X)].
 \ee
where we rely on the following measurability result whose proof is reported in Subsection \ref{subsect:measurabilityI}.

\begin{Lemma}\label{lemma:mesurabilityS}
For $\P\in \Pc(\Omega)$, the map $\wideparen{\supp}\P_X$ is analytically measurable, and the map $\csupp\left(\P_X|_{\partial \widehat{K}(X)}\right)$ is $\mu-$measurable.
\end{Lemma}

Now we prove a first Lemma about the existence of a maximal support probability.

\begin{Lemma}\label{lemma:existence_proba}
There exists $\widehat\P\in\Mc(\mu,\nu)$ such that for all $\P\in\Mc(\mu,\nu)$ we have the inclusion $\csupp\,\P_X
 \subset
 \csupp\,\widehat\P_X$, $\mu-$a.s.
\end{Lemma}

\proof
We proceed in two steps:

\no \underline{\it Step 1:} We first prove existence for the problem \ref{supportmax}. Let $(\P^n)_{n\ge 1}\subset\Mc(\mu,\nu)$ be a maximizing sequence. Then the measure $\widehat\P := \sum_{n\ge 1}2^{-n}\,\P^n\in\Mc(\mu,\nu)$, and satisfies $\wideparen{\supp}\P^n_X\subset \wideparen{\supp}\widehat\P_X$ for all $n\geq 1$. Consequently $\mu[G(\wideparen{\supp}_X\P^n_X)]\le\mu[G(\wideparen{\supp}\widehat\P_X)]$, and therefore $M=\mu[G(\wideparen{\supp}\widehat\P_X)]$.
\\
\no \underline{\it Step 2:} We next prove that $\wideparen{\supp}\P_X\subset \wideparen{\supp}\widehat\P_X$, $\mu$-a.s. for all $\P\in\Mc(\mu,\nu)$. Indeed, the measure $\Pb:=\frac{\widehat\P+\P}{2}\in\Mc(\mu,\nu)$ satisfies $M\geq\mu[G (\wideparen{\supp}\Pb_X)]\geq\mu[G (\wideparen{\supp}\widehat\P_X)] = M$, implying that $G(\wideparen{\supp}\Pb_X)=G(\wideparen{\supp}\widehat\P_X)$, $\mu-$a.s. The required result now follows from the inclusion $\wideparen{\supp}\widehat\P_X\subset \wideparen{\supp}\Pb_X$.
\ep\\
\\
\no{\bf Proof of Proposition \ref{prop:icpdual} (iii)}
Let $\widehat{\P}\in\Mc(\mu,\nu)$ from Lemma \ref{lemma:existence_proba}, if we denote $S(X):=\csupp\widehat\P_X$, then we have $\supp(\P_X)\subset S(X)$, $\mu-$a.s. Then $\{Y\notin S(X)\}$ is $\Mc(\mu,\nu)-$polar. By Lemma \ref{lemma:measurabilityFrontier}, $\{Y\notin S(X)\}\cup \{X\notin N_\mu'\}$ is Borel for some $N_\mu'\in \Nc_\mu$. By Theorem \ref{thm:CharacPolar}, we see that $\{Y\notin S(X)\}\subset\{Y\notin S(X)\}\cup \{X\notin N_\mu'\}\subset\{X\in N_\mu\}\cup\{Y\in N_\nu\}\cup\{Y\notin K_\theta(X)\}$, and therefore
 \b*
 \{Y\in S(X)\}
 &\supset& 
 \{X\notin N_\mu\}\cap\{Y\in K_\theta(X)\setminus N_\nu\},
 \e*
for some $N_\mu\in\Nc_\mu$, $N_\nu\in\Nc_\nu$, and $\theta\in\widehat{\Tc}(\mu,\nu)$. The last inclusion implies that $K_\theta(X)\setminus N_\nu\subset S(X)$, $\mu$-a.s. However, by Proposition \ref{prop:icpdual} (ii), $\widehat{K}(X)\subset\conv\big(\dom\theta(X,\cdot)\setminus N_\nu\big)$, $\mu-$a.s. Then, since $S(X)$ is closed and convex, we see that $\cl \widehat{K}(X)\subset S(X)$. 

To obtain the reverse inclusion, we recall from Proposition \ref{prop:icpdual} (i) that $\{Y\in \cl \widehat{K}(X)\}$, $\Mc(\mu,\nu)-$q.s. In particular $\widehat\P[Y\in\cl \widehat{K}(X)]=1$, implying that $S(X) \subset \cl \widehat{K}(X)$, $\mu$-a.s. as $\cl \widehat{K}(X)$ is closed convex. Finally, recall that by definition $I:=\ri S$ and therefore $\widehat{K}(X) = \cl I(X)$, $\mu-$a.s.
\ep

\begin{Lemma}\label{lemma:JisK}
We may choose $\widehat\P\in\Mc(\mu,\nu)$ in Theorem \ref{thm:icp} so that for all $\P\in\Mc(\mu,\nu)$ and $y\in\R^d$,
 $$
\mu\big[ \P_X[\{y\}]>0\big]
 \le 
 \mu\big[ \widehat \P_X[\{y\}]>0\big],
 ~\mbox{and}~
 \supp\,\P_X|_{\partial I(X)}
 \subset
 \csupp\,\widehat\P_X|_{\partial I(X)},~\mu-\mbox{a.s.}
 $$
 In this case the set-valued maps $\underline{J}(X):= I(X)\cup \big\{y\in\R^d:\nu[{y}]>0\mbox{ and }\widehat\P_X\big[\{y\}\big]>0\big\}$, and $\bar{J}(X):=I(X)\cup \csupp\,\widehat\P_X |_{\partial I(X)}$ are unique $\mu-$a.s. Furthermore $\underline{J}(X)=\underline{K}(X)$, $\bar{J}(X)=\bar{K}(X)$, and $J_\theta(X) = K_\theta(X)$, $\mu-$a.s. for all $\theta\in\widehat{\Tc}(\mu,\nu)$.
\end{Lemma}
\proof
\no\underline{\it Step 1:}
By the same argument as in the proof of Lemma \ref{lemma:existence_proba}, we may find $\widehat\P'\in\Mc(\mu,\nu)$ such that 
 \be\label{supportmax2}
 M':=
  \underset{\P\in\Mc(\mu,\nu)}{\sup}\mu\Big[G\Big(\csupp\left(\P_X|_{\partial \widehat{K}(X)}\right)\Big)\Big]=\mu\Big[G\Big(\csupp\left(\widehat\P_X'|_{\partial \widehat{K}(X)}\right)\Big)\Big].
 \ee
We also have similarly that $\csupp\left(\P_X|_{\partial \widehat{K}(X)}\right)\subset \csupp\left(\widehat\P_X'|_{\partial \widehat{K}(X)}\right)$, $\mu$-a.s. for all $\P\in\Mc(\mu,\nu)$. Then we prove similarly that $S'(X):=\csupp\left(\widehat\P_X'|_{\partial \widehat K(X)}\right)=D(X)$, $\mu-$a.s., where recall that $D$ is the optimizer for \eqref{eq:irreducible-optim1}. Indeed, by the previous step, we have $\csupp(\P_X|_{\partial \widehat K(X)})\subset S'(X)$, $\mu-$a.s. Then $\{Y\notin S'(X)\cup \widehat K(X)\}$ is $\Mc(\mu,\nu)-$polar. By Theorem \ref{thm:CharacPolar}, we see that $\{Y\notin S'(X)\cup \widehat K(X)\}\subset\{X\in N_\mu\}\cup\{Y\in N_\nu\}\cup\{Y\notin K_\theta(X)\cup \widehat K(X)\}$, or equivalently,
 \be
 \{Y\in S'(X)\cup \widehat K(X)\}
 &\supset& 
 \{X\notin N_\mu\}\cap\{Y\in K_\theta(X)\setminus N_\nu\},
 \ee
for some $N_\mu\in\Nc_\mu$, $N_\nu\in\Nc_\nu$, and $\theta\in\widehat{\Tc}(\mu,\nu)$. Similar to the previous analysis, we have $K_\theta(X)\setminus N_\nu\setminus \widehat K(X)\subset S'(X)$, $\mu$-a.s. Then, since $S'(X)$ is closed and convex, we see that $D(X)\subset S'(X)$. 

To obtain the reverse inclusion, we recall from Proposition \ref{prop:Jqs} that $\{Y\in \bar K(X)\}$, $\Mc(\mu,\nu)-$q.s. In particular $\widehat\P'[Y\in \widehat K(X)\cup D(X)]=1$, implying that $S'(X) \subset D(X)$, $\mu$-a.s. By Proposition \ref{prop:icpdual} (iii), we have $\bar{J}(X) = (I\cup S')(X)=(\widehat{K}\cup D)(X) = \bar{K}(X)$, $\mu-$a.s.

Finally, $\frac{\widehat\P+\widehat\P'}{2}$ is optimal for both problems \eqref{supportmax}, and \eqref{supportmax2}. By definition, the equality $J_\theta(X) = K_\theta(X)$, $\mu-$a.s. for $\theta\in\widehat{\Tc}(\mu,\nu)$ immediately follows.

\no\underline{\it Step 2:} Let $y\in{\rm atom}(\nu)$, if $y$ is an atom of $\gamma_1\in\Pc(\R^d)$ and $\gamma_2\in\Pc(\R^d)$, then $y$ in an atom of $\lambda\gamma_1+(1-\lambda)\gamma_2$ for all $0<\lambda<1$. By the same argument as in Step 1, we may find $\widehat\P^y\in\Mc(\mu,\nu)$ such that
 \be\label{supportmax3}
 M_y:=
  \underset{\P\in\Mc(\mu,\nu)}{\sup}\mu\Big[\P_X\big[\{y\}\cap \cl \widehat K(X)\big]>0\Big] = \mu\Big[\widehat\P_{X}^y\big[\{y\}\cap \cl \widehat K(X)\big]>0\Big].
 \ee

We denote $S_y(X):=\supp\widehat\P_X^y|_{\aff \widehat K(X)\cap\{y\}}$. Recall that $D_y$ is the notation for the optimizer of problem \eqref{eq:irreducible-optim3}. We consider the set $N:=\big\{Y\notin(\cl \widehat K(X)\setminus \{y\})\cup S_y(X)\big\}$. $N$ is polar as $Y\in\cl \widehat K(X)$, q.s., and by definition of $S_y$. Then $N\subset\{X\in N_\mu\}\cup\{Y\in N_\nu\}\cup\{Y\notin K_\theta(X)\}$, or equivalently,
 \be
 \big\{Y\notin(\cl \widehat K(X)\setminus \{y\})\cup S_y(X)\big\}
 &\supset& 
 \{X\notin N_\mu\}\cap\{Y\in K_\theta(X)\setminus N_\nu\},
 \ee
for some $N_\mu\in\Nc_\mu$, $N_\nu\in\Nc_\nu$, and $\theta\in\widehat{\Tc}(\mu,\nu)$. Then $D_y(X)\subset K_\theta(X)\setminus N_\nu\subset\cl \widehat K(X)\setminus \{y\}\cup S_y(X)$, $\mu-$a.s. Finally $D_y(X)\subset S_y(X)$, $\mu-$a.s.

On the other hand, $S_y\subset D_y$, $\mu-$a.s., as if $\widehat\P_{X}^y[\{y\}]>0$, we have $\theta(X,y)<\infty$, $\mu-$a.s. at the corresponding points. Hence, $D_y(X) = S_y(X)$, $\mu-$a.s. Now if we sum up the countable optimizers for $y\in{\rm atom}(\nu)$, with the previous optimizers, then the probability $\widehat\P$ we get is an optimizer for \eqref{supportmax}, \eqref{supportmax2}, and \eqref{supportmax3}, for all $y\in\R^d$ (the optimum is $0$ if it is not an atom of $\nu$). Furthermore, the $\mu-$a.e. equality of the maps $S_y$ and $D_y$ for these countable $y\in{\rm atom}(\nu)$ is preserved by this countable union, then together with Proposition \ref{prop:icpdual} (iii), we get $\underline{J} = \underline{K}$, $\mu-$a.s.
\ep\\

As a preparation to prove the main Theorem \ref{thm:icp}, we need the following lemma, which will be proved in Subsection \ref{subsect:measurabilityI}.

\begin{Lemma}\label{lemma:anal}
Let $F:\R^d\longrightarrow\ri\,\Kcirc$ be a $\gamma-$measurable function for some $\gamma\in\Pc(\R^d)$, such that $x\in F(x)$ for all $x\in\R^d$, and $\{F(x):x\in\R^d\}$ is a partition of $\R^d$. Then up to a modification on a $\gamma-$null set, $F$ can be chosen in addition to be analytically measurable.
\end{Lemma}

\no {\bf Proof of Theorem \ref{thm:icp}}
Existence holds by Lemma \ref{lemma:existence_proba} above, (i) is a consequence of Lemma \ref{lemma:mesurabilityS}, and (ii) directly stems from Lemma \ref{lemma:measurabilityI} (iii) together with Proposition \ref{prop:icpdual} (iii). Now we need to deal with the measurability issue. Lemma \ref{lemma:anal} allows to modify $\ri\,\csupp \widehat\P_X$ to get (ii) while preserving its analytic measurability, we denote $I$ its modification. However, we need to modify $\widehat\P_X$ to get the result. As $\csupp\widehat\P_X$ is analytically measurable by Lemma \ref{lemma:mesurabilityS}, the set of modification $N_\mu:=\{\csupp\widehat\P_X\neq \cl I(X)\}\in\Nc_\mu$ is analytically measurable. Then we may redefine $\widehat{\P}_X$ on $N_\mu$, so as to preserve a kernel for $\widehat{\P}$. By the same arguments than the proof of Lemma \ref{lemma:measurabilityI} (ii), the measure-valued map $\kappa_X:= g_{I(X)}$ is a kernel thanks to the analytic measurability of $I$, recall the definition of $g_K$ given by \eqref{eq:defG}. Furthermore, $\csupp\,\kappa_X = I(X)$ pointwise by definition. Then a suitable kernel modification from which the result follows is given by
$$\widehat\P_X':= \mathbf{1}_{\{X\in N_\mu\}}\kappa_X + \mathbf{1}_{\{X\notin N_\mu\}}\widehat{\P}_X.$$
\ep\\

\no {\bf Proof of Proposition \ref{prop:CharacProb}}
The existence and the uniqueness are given by Lemma \ref{lemma:JisK} and the other properties follow from the identity between the $J$ maps and the $K$ maps, also given by the Lemma, together with Proposition \ref{prop:Jqs}.
\ep\\

\no {\bf Proof of Theorem \ref{thm:CharacPolar}}
We simply apply Lemma \ref{lemma:JisK} to replace $K_\theta$ by $J_\theta$ in Proposition \ref{prop:polar}.
\ep

\section{Measurability of the irreducible components}
\label{sect:measurability}
\setcounter{equation}{0}

\subsection{Measurability of G}\label{subsect:measurabilityG}

\no {\bf Proof of Lemma \ref{lemma:measurabilityI} (ii)}
As $\R^d$ is locally compact, the Wijsman topology is locally equivalent to the Hausdorff topology\footnote{The Haussdorff distance on the collection of all compact subsets of a compact metric space $(\Xc,{\rm d})$ is defined by $d_H(K_1,K_2)=\sup_{x\in\Xc}\left|\dist(x,K_1)-\dist(x,K_2)\right|,$ for $K_1,K_2\subset\Xc$, compact subsets.}, i.e. as $n\to\infty$, $K_n\longrightarrow K$ for the Wijsman topology if and only if $K_n\cap B_M\longrightarrow K\cap B_M$ for the Hausdorff topology, for all $M\ge 0$.

We first prove that $K\longmapsto\dim\aff K$ is a lower semi-continuous map $\Kc\to\R$. Let $(K_n)_{n\ge 1}\subset\Kc$ with dimension $d_n\le d'\le d$ converging to $K$. We consider $A_n:=\aff\,K_n$. As $A_n$ is a sequence of affine spaces, it is homeomorphic to a $d+1$-uplet. Observe that the convergence of $K_n$ allow us to chose this $d+1$-uplet to be bounded. Then up to taking a subsequence, we may suppose that $A_n$ converges to an affine subspace $A$ of dimension less than $d'$. By continuity of the inclusion under the Wijsman topology, $K\subset A$ and $\dim K\le\dim A \le d'$.

We next prove that the mapping $K\mapsto g_K(K)$ is continuous on $\{\dim K = d'\}$ for $0\le d'\le d$, which implies the required measurability. Let $(K_n)_{n\ge 1}\subset\Kc$ be a sequence with constant dimension $d'$, converging to a $d'-$dimensional subset, $K$ in $\Kc$. Define $A_n:=\aff K_n$ and $A:=\aff K$, $A_n$ converges to $A$ as for any accumulation set $A'$ of $A_n$, $K\subset A'$ and $\dim A' = \dim A$, implying that $A'=A$. Now we consider the map $\phi_n:A_n\to A$, $x\mapsto proj_A(x)$. For all $M>0$, it follows from the compactness of the closed ball $B_M$ that $\phi_n$ converges uniformly to identity as $n\to\infty$ on $B_M$. Then, $\phi_n(K_n)\cap B_M\longrightarrow K\cap B_M$ as $n\to \infty$, and therefore $\leb_A[\phi_n(K_n\cap B_M)\setminus K]+\leb_A[K\setminus\phi_n(K_n)\cap B_M]\longrightarrow 0$. As the Gaussian density is bounded, we also have
$$ g_A[\phi_n(K_n\cap B_M)]\longrightarrow  g_A[K\cap B_M].$$
We next compare $ g_A[\phi_n(K_n\cap B_M)]$ to $ g_{K_n}(K_n\cap B_M)$. As $(\phi_n)$ is a sequence of linear functions that converges uniformly to identity, we may assume that $\phi_n$ is a $\Ctn^1-$diffeomorphism. Furthermore, its constant Jacobian $J_n$ converges to $1$ as $n\to\infty$. Then, 
\b*
\int_{K_n\cap B_M}\frac{e^{-|\phi_n(x)|^2/2}}
                       {(2\pi)^{d'/2}}
                  \leb_{K_n}(dx) 
&=&
\int_{\phi_n(K_n\cap B_M)}\frac{e^{-|y|^2/2}J_n^{-1}}
                               {(2\pi)^{d'/2}}
                          \leb_{A}(dy)
\\
&=&
J_n^{-1} g_A[\phi_n(K_n\cap B_M)].
\e*
As the Gaussian distribution function is $1$-Lipschitz, we have
 $$
 \left|\int_{K_n\cap B_M}\frac{e^{-|\phi_n(x)|^2/2}}
                              {(2\pi)^{d'/2}}
                         \leb_{K_n}(dx)
        - g_{K_n}(K_n\cap B_M)
 \right|
 \;\le\;
 \leb_{K_n}[K_n\cap B_M]|\phi_n-Id_A|_\infty,
 $$
where $|\cdot|_\infty$ is taken on $K_n\cap B_M$. Now for arbitrary $\epsilon>0$, by choosing $M$ sufficiently large so that $ g_V[V\setminus B_M]\le\epsilon$ for any $d'-$dimensional subspace $V$, we have
\b*
\left| g_{K_n}[K_n]- g_{K}[K]\right|&\le&\left| g_{K_n}[K_n\cap B_M]- g_{A}[K\cap B_M]\right|+2\epsilon\\
&\hspace{-20mm}\le &
\hspace{-10mm}
\left| g_{K_n}[K_n\cap B_M]-\int_{K_n\cap B_M}C\exp\left(\frac{-|\phi_n(x)|^2}{2}\right)\leb_{K_n}(dx)\right|\\
&&
+\left|J_n^{-1} g_A[\phi_n(K_n\cap B_M)]- g_{A}[K\cap B_M]\right|+2\epsilon
\\
&\hspace{-20mm}\le &
\hspace{-10mm}
4\epsilon,
\e*
for $n$ sufficiently large, by the previously proved convergence. Hence $G_{d'}:=G\big|_{\dim^{-1}\{d'\}}$ is continuous, implying that $G :K\longmapsto \sum_{d'=0}^d\mathbf{1}_{\dim^{-1}\{d'\}}(K)G_{d'}(K)$ is Borel-measurable.
\ep

\subsection{Further measurability of set-valued maps}\label{subsect:measurabilityI}

This subsection is dedicated to the proof of Lemmas \ref{lemma:measurabilityI} (i), \ref{lemma:measurabilityFrontier}, and \ref{lemma:mesurabilityS}. In preparation for the proofs, we start by giving some lemmas on measurability of set-valued maps. Let $\Ac$ be a $\sigma-$algebra of $\R^d$. In practice we will always consider either the $\sigma-$algebra of Borel sets, the $\sigma-$algebra of analytically measurable sets, or the $\sigma-$algebra of universally measurable sets.

\begin{Lemma}\label{lemma:unionbim}
Let $(F_n)_{n\ge 1}\subset\L^\Ac(\R^d,\Kc)$. Then $\cl\cup_{n\ge 1}F_n$ and $\cap_{n\ge 1}F_n$ are $\Ac-$measurable.
\end{Lemma}
\proof
The measurability of the union is a consequence of Propositions 2.3 and 2.6 in Himmelberg \cite{himmelberg1975measurable}. The measurability of the intersection follows from the fact that $\R^d$ is $\sigma$-compact, together with Corollary 4.2 in \cite{himmelberg1975measurable}.
\ep

\begin{Lemma}\label{lemma:bimesconv}
Let $F\in\L^\Ac(\R^d,\Kc)$. Then, $\cl\conv F$, $\aff F$, and $\cl\rf_X \cl\conv F$ are $\Ac-$measurable.
\end{Lemma}

\proof
The measurability of $\cl\conv F$ is a direct application of Theorem 9.1 in \cite{himmelberg1975measurable}.

We next verify that $\aff F$ is measurable. Since the values of $F$ are closed, we deduce from Theorem 4.1 in Wagner \cite{wagner1977survey}, that we may find a measurable $x\longmapsto y(x)$, such that $y(x)\in F(x)$ if $F(x)\neq\emptyset$, for all $x\in\R^d$. Then we may write $\aff F(x) = \cl\,\conv\,\cl\,\cup_{q\in\Q}\big(y(x)+q\left(F(x)-y(x)\right)\big)$ for all $x\in\R^d$. The measurability follows from Lemmas \ref{lemma:unionbim}, together with the first step of the present proof.

We finally justify that $\cl\rf_X \cl\conv F$ is measurable. We may assume that $F$ takes convex values. By convexity, we may reduce the definition of $\rf_x$ to a sequential form:
\b*
\cl\rf_x F(x) 
&\hspace{-5mm}=&
\hspace{-5mm}
\cl\cup_{n\ge 1}\left\{y\in \R^d , y+\frac{1}{n}(y-x)\in F(x) \text{ and } x-\frac{1}{n}(y-x)\in F(x)\right\}\\
&\hspace{-30mm}=&
\hspace{-15mm}
 \cl\cup_{n\ge 1} \left[\left\{y\in \R^d , y+\frac{1}{n}(y-x)\in F(x)\right\} \cap \left\{y\in \R^d , x-\frac{1}{n}(y-x)\in F(x)\right\}\right]\\
&\hspace{-30mm}=&
\hspace{-15mm}
\cl\cup_{n\ge 1} \left[\left(\frac{1}{n+1}x+\frac{n}{n+1}F(x)\right) \cap \left(-(n+1)x - nF(x)\right)\right],\\
\e*
so that the required measurability follows from Lemma \ref{lemma:unionbim}.
\ep\\

We denote $\Sc$ the set of finite sequences of positive integers, and $\Sigma$ the set of infinite sequences of positive integers. Let $s\in\Sc$, and $\sigma\in\Sigma$. We shall denote $s<\sigma$ whenever $s$ is a prefix of $\sigma$. 

\begin{Lemma}\label{lemma:suslinmeas}
Let $(F_s)_{s\in \Sc}$ be a family of universally measurable functions $\R^d\longrightarrow \Kc$ with convex image. Then the mapping $\cl\conv\big(\cup_{\sigma\in\Sigma}\cap_{s<\sigma}F_s\big)$ is universally measurable.
\end{Lemma}

\proof
Let $\Uc$ the collection of universally measurable maps from $\R^d$ to $\Kc$ with convex image. For an arbitrary $\gamma\in\Pc(\R^d)$, and $F:\R^d\longrightarrow\Kc$, we introduce the map
 \b*
 \gamma G^*[F]:=\inf_{F\subset F'\in\Uc}\gamma G[F'],
 &\mbox{where}&
 \gamma G[F']:= \gamma\big[G\big(F'(X)\big)\big]~\mbox{for all}~F'\in\Uc.
 \e*
Clearly, $\gamma G$ and $\gamma G^*$ are non-decreasing, and it follows from the dominated convergence theorem that $\gamma G$, and thus $\gamma G^*$, are upward continuous.
\\
\no \underline{\it Step 1:} In this step we follow closely the line of argument in the proof of Proposition 7.42 of Bertsekas and Shreve \cite{bertsekas1978stochastic}. Set $F:=\cl\conv\big(\cup_{\sigma\in\Sigma}\cap_{s<\sigma}F_s\big)$, and let $(\bar F_n)_n$ a minimizing sequence for $\gamma G^*[F]$. Notice that $F\subset\bar F := \cap_{n\ge 1}\bar F_n\in\Uc$, by Lemma \ref{lemma:unionbim}. Then $\bar F$ is a minimizer of $\gamma G^*[F]$.

For $s,s'\in S$, we denote $s\le s'$ if they have the same length $|s|=|s'|$, and $s_i\le s_i'$ for $1\le i\le |s|$. For $s\in S$, let
 \b*
 R(s):=\cl\,\conv\cup_{s'\le s}\cup_{\sigma>s'}\cap_{s''<\sigma}F_{s''}
 &\mbox{and}&
 K(s):=\cl\,\conv\cup_{s'\le s}\cap_{j=1}^{|s'|}F_{s_1',...,s_j'}.
 \e*
Notice that $K(s)$ is universally measurable, by Lemmas \ref{lemma:unionbim} and \ref{lemma:bimesconv}, and
 $$
 R(s)\subset K(s),
 ~\cl\cup_{s_1\ge 1}R(s_1) = F,
 ~\mbox{and}~
 \cl\cup_{s_k\ge 1}R(s_1,...,s_k) = R(s_1,...,s_{k-1}).
 $$
By the upwards continuity of $\gamma G^*$, we may find for all $\epsilon>0$ a sequence $\sigma^\epsilon\in\Sigma$ s.t.
 \b*
 \gamma G^*[F]\le \gamma G^*[R(\sigma^\epsilon_1)]+2^{-1}\epsilon,
 &\mbox{and}&
 \gamma G^*[R(\underline{\sigma}_{k-1})]
 \le 
 \gamma G^*[R(\underline{\sigma}_k)]+2^{-k}\epsilon,
 \e*
for all $k\ge 1$, with the notation $\underline{\sigma}^\eps_k:=(\sigma_1^\epsilon,\ldots,\sigma_k^\eps)$. Recall that the minimizer $\overline{F}$ and $K(s)$ are in $\Uc$ for all $s\in\Sc$. We then define the sequence $K_k^\epsilon:=\overline{F}\cap K(\underline{\sigma}^\epsilon_k)\in\Uc$, $k\ge 1$, and we observe that
 \begin{equation}\label{step2:BS-hyp}
 (K_k^\epsilon)_{k\ge 1}~~\mbox{decreasing,}
 ~\underline{F}^\epsilon:=\cap_{k\ge 1}K_k^\epsilon\subset F,
 ~\mbox{and}~
 \gamma G[K_k^\epsilon]\ge \gamma G^*[F]-\epsilon = \gamma G[\overline{F}]-\epsilon,
 \end{equation}
by the fact that $R(\underline{\sigma}^\epsilon_k)\subset K_k^\epsilon$. We shall prove in Step 2 that, for an arbitrary $\alpha>0$, we may find $\eps=\eps(\alpha)\le\alpha$ such that \eqref{step2:BS-hyp} implies that 
 \be\label{step2:BS-conc}
 \gamma G[\underline{F}^{\epsilon}]
 \ge 
 \inf_{k\ge 1}\gamma G[K_k^{\epsilon}] - \alpha
 \ge 
 \gamma G[\overline{F}]-\epsilon-\alpha.
 \ee
Now let $\alpha=\alpha_n:=n^{-1}$, $\eps_n:=\epsilon(\alpha_n)$, and notice that $\underline{F}:= \cl\conv\cup_{n\ge 1}\underline{F}^{\epsilon_n}\in\Uc$, with $\underline{F}^{\epsilon_n}\subset\underline{F}\subset F \subset\overline{F}$, for all $n\ge 1$. Then, it follows from \eqref{step2:BS-conc} that $\gamma G[\underline{F}]= \gamma G[\overline{F}]$, and therefore $\underline{F}= F =\overline{F}$, $\gamma -$a.s. In particular, $F$ is $\gamma-$measurable, and we conclude that $F\in\Uc$ by the arbirariness of $\gamma \in\Pc(\R^d)$.

\no \underline{\it Step 2:} It remains to prove that, for an arbitrary $\alpha>0$, we may find $\eps=\eps(\alpha)\le\alpha$ such that \eqref{step2:BS-hyp} implies \eqref{step2:BS-conc}. This is the point where we have to deviate from the argument of \cite{bertsekas1978stochastic} because $\gamma G$ is not downwards continuous, as the dimension can jump down. 

Set $A_n:=\{G\big(\overline{F}(X)\big)-\dim \overline{F}(X)\le 1/n\}$, and notice that $\cap_{n\ge 1}A_n=\emptyset$. Let $n_0\ge 1$ such that $\gamma[A_{n_0}]\le \frac12\frac{\alpha}{d+1}$, and set $\epsilon := \frac12\frac{1}{n_0}\frac{\alpha}{d+1}>0$. Then, it follows from \eqref{step2:BS-hyp} that
 \be
 \gamma\big[\inf_n G(K^\epsilon_n)-\dim\overline{F}\le 0\big]
 &\le&
 \gamma\big[\inf_n G(K^\epsilon_n)-G(\overline{F})\le n_0^{-1}\big]
 \\
 &&
 + \gamma\big[G(\overline{F})-\dim\overline{F}\le -n_0^{-1}\big]
 \nonumber\\
 &\le&
 n_0\big(\gamma\big[G(\overline{F})\big]
         -\gamma\big[\inf_n G(K^\epsilon_n)\big]
    \big)
 +\gamma\big[A_{n_0}\big]
 \nonumber\\
 &=&
 n_0\big(\gamma\big[G(\overline{F})\big]
         -\inf_n \gamma\big[G(K^\epsilon_n)\big]
    \big)
 +\gamma\big[A_{n_0}\big]
 \nonumber
 \\
 &\le&
 n_0\epsilon+\frac12\frac{\alpha}{d+1} \;=\; \frac{\alpha}{d+1},
 \label{last-step-meas}
 \ee
where we used the Markov inequality and the monotone convergence theorem. Then:
 \b*
 \gamma\big[\inf_n G(K^\epsilon_n)-G\big(\underline{F}^\epsilon\big)\big]
 &\le&
 \gamma\Big[\1_{\{\inf_n G(K^\epsilon_n)-\dim\overline{F}\le 0 \}}
            \big(\inf_n G(K^\epsilon_n)-G\big(\underline{F}^\epsilon\big)\big)
 \\
 &&
 \hspace{5mm} + \1_{\{\inf_n G(K^\epsilon_n)-\dim\overline{F}>0 \}}
            \big(\inf_n G(K^\epsilon_n)-G\big(\underline{F}^\epsilon\big)\big)
       \Big]
 \\
 &\le&
 \gamma\Big[(d+1)\1_{\{\inf_n G(K^\epsilon_n)-\dim\overline{F}\le 0 \}}
 \\
 &&
 \hspace{5mm} + \1_{\{\inf_n G(K^\epsilon_n)-\dim\overline{F}>0 \}}
            \big(\inf_n G(K^\epsilon_n)-G\big(\underline{F}^\epsilon\big)\big)
       \Big].
 \e*
We finally note that $\inf_n G(K^\epsilon_n)-G\big(\underline{F}^\epsilon\big)=0$ on $\{\inf_n G(K^\epsilon_n)-\dim\overline{F}>0 \}$. Then \eqref{step2:BS-conc} follows by substituting the estimate in \eqref{last-step-meas}.
 \ep\\

\noindent {\bf Proof of Lemma \ref{lemma:measurabilityI} (i)}
We consider the mappings $\theta:\Omega\to\bar\R_+$ such that $\theta = \sum_{k=1}^n \lambda_k 1_{C^1_k\x C^2_k}$ where $n\in\N$, the $\lambda_k$ are non-negative numbers, and the $C^1_k,C^2_k$ are closed convex subsets of $\R^d$. We denote the collection of all these mappings $\Fc$. Notice that $\cl\Fc$ for the pointwise limit topology contains all $\Leb^0_+(\Omega)$. Then for any $\theta\in\Leb^0_+(\Omega)$, we may find a family $(\theta_s)_{s\in\S}\subset \Fc$, such that $\theta = \inf_{\sigma\in\Sigma}\sup_{s<\sigma}\theta_s$. For $\theta\in\Leb_+^0(\Omega)$, and $n\ge 0$, we denote $F_\theta :x\longmapsto \cl\,\conv\,\dom\theta(x,\cdot)$, and $F_{\theta,n} :x\longmapsto \cl\,\conv\,\theta(x,\cdot)^{-1}([0,n])$. Notice that $F_\theta = \cl\cup_{n\ge 1}F_{\theta,n}$. Notice as well that $F_{\theta,n}$ is Borel measurable for $\theta\in\Fc$, and $n\ge 0$, as it takes values in a finite set, from a finite number of measurable sets. Let $\theta\in\Leb_+^0(\Omega)$, we consider the associated family $(\theta_s)_{s\in\S}\subset \Fc$, such that $\theta = \inf_{\sigma\in\Sigma}\sup_{s<\sigma}\theta_s$. Notice that $F_{\theta,n} = \cl \conv\big(\cup_{\sigma\in\Sigma}\cap_{s<\sigma}F_{\theta_s,n}\big)$ is universally measurable by Lemma \ref{lemma:suslinmeas}, thus implying the universal measurability of $F_\theta =\cl\,\dom\theta(X,\cdot)$ by Lemma \ref{lemma:unionbim}.

In order to justify the measurability of $\dom_X\theta$, we now define
 \b*
 F_{\theta}^0 := F_{\theta}
 &\mbox{and}&
 F_{\theta}^k 
 := 
 \cl\conv(\dom\theta(X,\cdot)\cap \aff\,\rf_XF_{\theta}^{k-1}),
 ~k \ge 1.
 \e*
Note that $F_{\theta}^{k} = \cl\cup_{n\ge 1}\big(\cl \conv\cup_{\sigma\in\Sigma}\cap_{s<\sigma}F_{\theta_s,n}\cap\aff\,\rf_xF_{\theta}^{k-1}\big)$. Then, as $F_{\theta}^0$ is universally measurable, we deduce that $\big(F_{\theta}^k\big)_{k\ge 1}$ are universally measurable, by Lemmas \ref{lemma:bimesconv} and \ref{lemma:suslinmeas}.

As $\dom_X\theta$ is convex and relatively open, the required measurability follows from the claim:
 \b*
 F_{\theta}^{d} 
 &=& 
 \cl\dom_X\theta.
 \e*
To prove this identity, we start by observing that $F_{\theta}^k(x) \supset \cl\dom_x\theta$. Since the dimension cannot decrease more than $d$ times, we have $\aff\,\rf_x F^d_{\theta}(x) = \aff F^{d}_{\theta}(x)$ and
 \b*
 F^{d+1}_{\theta}(x)
 &=&
 \cl\conv\big(\dom\theta(x,\cdot)\cap \aff\,\rf_xF_{\theta}^{d}(x)\big)
 \\
 &=&
 \cl\conv\big(\dom\theta(x,\cdot)\cap \aff \rf_x F_{\theta}^{d-1}(x)\big)
 \;=\;
 F_{\theta}^{d}(x).
 \e*
i.e. $(F^{d+1}_{\theta})_k$ is constant for $k\ge d$. Consequently,
\b*
\dim\rf_x\conv(\dom\theta(x,\cdot)\cap \aff\,\rf_xF_{\theta}^{d}(x)) &=& \dim F_{\theta}^{d}(x)\\
&\ge& \dim\conv(\dom\theta(x,\cdot)\cap \aff\,\rf_xF_{\theta}^{d}(x)).
\e*
As $\dim\conv\big(\dom\theta(x,\cdot)\cap \aff\,\rf_xF_{\theta}^{d}(x)\big)\ge \dim\rf_x\conv\big(\dom\theta(x,\cdot)\cap \aff\,\rf_xF_{\theta}^{d}(x)\big)$, we have equality of the dimension of $\conv\big(\dom\theta(x,\cdot)\cap \aff\,\rf_xF_{\theta}^{d}(x)\big)$ with its $\rf_x$. Then it follows from Proposition \ref{prop:ri} (ii) that $x\in \ri\,\conv\big(\dom\theta(x,\cdot)\cap \aff\,\rf_xF_{\theta}^{d}(x)\big),$ and therefore:
 \b*
 F_{\theta}^{d}(x)
 &=&  
 \cl\conv\big(\dom\theta(x,\cdot)\cap  \aff\,\rf_xF_{\theta}^{d}(x)\big)
 \\
 &=&
 \cl\ri\,\conv\big(\dom\theta(x,\cdot)\cap \aff\,\rf_xF_{\theta}^{d}(x)\big)
 \\
 &=&
 \cl\rf_x\conv\big(\dom\theta(x,\cdot)\cap \aff\,\rf_xF_{\theta}^{d}(x)\big)
 \;\subset\; 
 \cl\dom_x\theta.
 \e*
Hence $F_{\theta}^{d}(x)= \cl\dom_x\theta$.

Finally, $K_{\theta,A} = \dom_X(\theta+\infty 1_{\R^d\x A})$ is universally measurable by the universal measurability of $\dom_X$.
\ep\\

\no{\bf Proof of Lemma \ref{lemma:measurabilityFrontier}}
We may find $(F_n)_{n\ge 1}$, Borel-measurable with finite image, converging $\gamma-$a.s. to $F$. We denote $N_\gamma\in\Nc_\gamma$, the set on which this convergence does not hold. For $\epsilon>0$, we denote $F_k^\epsilon(X) := \{y\in\R^d:dist\big(y,F_k(X)\big)\le\epsilon\}$, so that
\b*
F(x) = \cap_{i\ge1}\liminf_{n\to\infty}F_n^{1/i}(x),&\mbox{for all}&x\notin N_\gamma.
\e*
Then, as $\mathbf{1}_{Y\in F(X)}\mathbf{1}_{X\notin N_\gamma} = \inf_{i\ge1}\liminf_{n\to\infty}\mathbf{1}_{Y\in F_n^{1/i}(X)}\mathbf{1}_{X\notin N_\gamma}$, the Borel-measurability of this function follows from the Borel-measurability of each $\mathbf{1}_{Y\in F_n^{1/i}(X)}$.

Now we suppose that $X\in \ri F(X)$ convex, $\gamma-$a.s. Up to redefining $N_\gamma$, we may suppose that this property holds on $N_\gamma^c$, then $\partial F(x) = \cap_{n\ge 1}F(x)\setminus \big(x+\frac{n}{n+1}(F(x)-x)\big)$, for $x\notin N_\gamma$. We denote $a := \mathbf{1}_{Y\in F(X)}\mathbf{1}_{X\notin N_\gamma}$. The result follows from the identity $\mathbf{1}_{Y\in \partial F(X)}\mathbf{1}_{X\notin N_\gamma} = a - \sup_{n\ge 1}a\big(X,X+\frac{n}{n+1}(Y-X)\big)$.
\ep\\

\no{\bf Proof of Lemma \ref{lemma:mesurabilityS}}
Let $\Kc_\Q:=\{K=\conv(x_1,\ldots,x_n): n\in\N,(x_i)_{i\le n}\subset\Q^d\}$. Then
 $$
 \csupp\P_x
 =
  \cl\cup_{N\ge 1}\cap\{K\in\Kc_\Q:\csupp\P_x\cap B_N\subset K\}
 =
 \cl\cup_{N\ge 1}\cap_{K\in\Kc_\Q}F_K^N(x),
 $$
where $F_K^N(x) := K$ if $\P_x[B_N\cap K] =\P_x[B_N]$, and $F_K^N(x) :=\R^d$ otherwise.
As for any $K\in\Kc_\Q$ and $N\ge 1$, the map $\P_X[B_N\cap K] -\P_X[B_N]$ is analytically measurable, then $F_K^N$ is analytically measurable. The required measurability result follows from lemma \ref{lemma:unionbim}.

Now, in order to get the measurability of $\csupp(\P_X|_{\partial I(X)})$, we have in the same way
 \b*
 \csupp(\P_X|_{\partial I(X)})
 &=&
 \cl\cup_{n\ge 1}\cap_{K\in\Kc_\Q}F_K'^N(x),
 \e*
where $F_K'^N(x) := K$ if $\P_x[\partial I(x)\cap B_N\cap K] =\P_x[\partial I(x)\cap B_N]$, and $F_K'^N(x) :=\R^d$ otherwise. As $\P_X[\partial I(X)\cap B_N\cap K] = \P_X[\mathbf{1}_{Y\in\partial I(X)}\mathbf{1}_{X\notin N_\mu}\mathbf{1}_{Y\notin B_N\cap K}]$, $\mu-$a.s., where $N_\mu\in\Nc_\mu$ is taken from Lemma \ref{lemma:measurabilityFrontier}, $\P_X[\partial I(X)\cap B_N\cap K]$ is $\mu-$measurable, as equal $\mu-$a.s. to a Borel function. Then similarly, $\P_X[\partial I(X)\cap B_N\cap K] -\P_X[\partial I(X)\cap B_N]$ is $\mu-$measurable, and therefore $\csupp(\P_X|_{\partial I(X)})$ is $\mu-$measurable.
\ep\\

\no{\bf Proof of Lemma \ref{lemma:anal}}
By $\gamma-$measurability of $F$, we may find a Borel function $F_B:\R^d\longrightarrow\ri\,\Kcirc$ such that $F = F_B$, $\gamma-$a.s. Let a Borel $N_\gamma\in\Nc_\gamma$ such that $F = F_B$ on $N_\gamma^c$. By the fact that $\ri\Kcirc$ is Polish, we may find a sequence $(F_n)_{n\ge 1}$ of Borel functions with finite image converging pointwise towards $F_B$ when $n\longrightarrow\infty$. We will give an explicit expression for $F_n$ that will be useful later in the proof. Let $(K_n)_{n\ge 1}\subset\ri\,\Kcirc$ a dense sequence,
\be\label{eq:defFn}
F_n(x) &:=&\argmin_{K\in(K_i)_{i\le n}}dist\big(F_B(x),K\big),
\ee

Where $dist$ is the distance on $\ri\,\Kcirc$ that makes it Polish, and we chose the $K$ with the smallest index in case of equality.

We fix $n\ge 1$, let $K\in F_n(N_\gamma^c)$, the image of $F_n$ outside of $N_\gamma$, and $A_K:= F_n^{-1}\big(\{K\}\big)$. We will modify the image of $F_n$ so that it is the same for all $x'\in F_B(x)= F(x)$, for all $x\in N_\gamma^c\cap A_K$. Then we consider the set $A'_K:= \cup_{x\in N_\gamma^c\cap A_K} F_B(x)$, we now prove that this set in analytic. By Theorem 4.2 (b) in \cite{wagner1977survey}, $Gr F_B:=\{Y\in \cl F_B(X)\}$ is a Borel set. Let $\lambda >0$, we define the affine deformation $f_\lambda:\Omega\longrightarrow \Omega$ by $f_\lambda(X,Y) := \big(X,X+\lambda(Y-X)\big)$. By the fact that for $k\ge 1$, $f_{1-1/k}(Gr F_B)$ is Borel together with the fact that $x\in f_B(x)$ for $x\notin N_\gamma$, we have
$$\{Y\in F_B(X)\}\cap \{X\notin N_\gamma\} = \cup_{k\ge 1}f_{1-1/k}(Gr F_B)\cap \{X\notin N_\gamma\}.$$
Therefore, $\{Y\in F_B(X)\}\cap \{X\notin N_\gamma\}$ is Borel, and so is $\{Y\in F_B(X)\}\cap \{X\in N_\gamma^c\cap A_K\}$. Finally,
$$A'_K = Y\big(\{Y\in F_B(X)\}\cap \{X\in N_\gamma^c\cap A_K\}\big),$$
therefore, $A'_K$ is the projection of a Borel set, which is one of the definitions of an analytic set (see Proposition 7.41 in \cite{bertsekas1978stochastic}). Now we define a suitable modification of $F_n$ by $F_n'(x) := K$ for all $x\in A'_K$, we do this redefinition for all $K\in F_B(N_\gamma^c)$. Notice that thanks to the definition \eqref{eq:defFn} and the fact that $F_B(x) = F_B(x')$ if $x,x'\notin N_\gamma$ and $x'\in F_B(x) = F(x)$, we have the inclusion $A'_K\subset A_K\cup N_\gamma$. Then the redefinitions of $F_n$ only hold outside of $N_\gamma$, furthermore for different $K_1, K_2\in F_n(N_\gamma^c)$, $A'_{K_1}\cap A'_{K_2} = \emptyset$ as the value of $F_n(x)$ only depends on the value of $F_B(x)$ by \eqref{eq:defFn}. Notice that
\be\label{eq:Ngammaprime}
N_\gamma' := \left(\cup_{K\in F_n(N_\gamma^c)}A'_K\right)^c =  \left(\cup_{x\notin N_\gamma}F_B(x)\right)^c\subset N_\gamma 
\ee
is analytically measurable, as the complement of an analytic set, and does not depend on $n$. For $x\in N_\gamma'$, we define $F_n'(x):= \{x\}$. Notice that $F_n'$ is analytically measurable as the modification of a Borel Function on analytically measurable sets.

Now we prove that $F_n'$ converges pointwise when $n\longrightarrow \infty$. For $x\in N_\gamma'$, $F_n'(x)$ is constant equal to $\{x\}$, if $x\notin N_\gamma'$, by \eqref{eq:Ngammaprime} $x\in \cup_{x\notin N_\gamma}F_B(x)$, and therefore $F_n'(x) = F_B(x')=F(x')$ for some $x\in N_\gamma^c$, for all $n\ge 1$. Then as $F_n'(x')$ converges to $F(x')$, $F_n'(x)$ converges to $F(x)$. Let $F'$ be the pointwise limit of $F_n'$. the maps $F_n'$ are analytically measurable, and therefore, so does $F'$. For all $n\ge 1$, $F_n' = F_n$, $\gamma-$a.e. and therefore $F' = F_B = F$, $\gamma-$a.e. Finally, $F'(N_\gamma^c) = F(N_\gamma^c)$, and $\cup F(N_\gamma^c) = (N_\gamma')^c$. By property of $F$, $F'(N_\gamma^c)$ is a partition of $(N_\gamma')^c$ such that $x\in F'(x)$ for all $x\notin N_\gamma'$. On $N_\gamma'$, this property is trivial as $F'(x) = \{x\}$ for all $x\in N_\gamma'$.
\ep

\section{Properties of tangent convex functions}
\setcounter{equation}{0}

\subsection{x-invariance of the y-convexity}
\label{subsect:convexityinvariance}

We first report a convex analysis lemma.

\begin{Lemma}\label{lemma:convexlsc}
Let $f:\R^d\to\bar\R$ be convex finite on some convex open subset $U\subset\R^d$. We denote $f_*:\R^d\to\bar\R$ the lower-semicontinuous envelop of $f$ on $U$, then
\b*
f_*(y)=\lim_{\epsilon\searrow 0}f\big(\epsilon x+(1-\epsilon)y\big),&\mbox{for all}&(x,y)\in U\x\cl U.
\e*
\end{Lemma}
\proof
$f_*$ is the lower semi-continuous envelop of $f$ on $U$, i.e. the lower semi-continuous envelop of $f':=f+\infty\mathbf{1}_{U^c}$. Notice that $f'$ is convex $\R^d\longrightarrow\R\cup\{\infty\}$. Then by Proposition 1.2.5 in Chapter IV of \cite{hiriart2013convex}, we get the result as $f=f'$ on $U$.
\ep\\

\no {\bf Proof of Proposition \ref{prop:ConvexityInvariance}}
The result is obvious in $\Tbf(\Cfrak_1)$, as the affine part depending on $x$ vanishes. We may use $N_\nu = \emptyset$. Now we denote $\Tc$ the set of mappings in $\Theta_\mu$ such that the result from the proposition holds. Then we have $\Tbf(\Cfrak_1)\subset \Tc$.

We prove that $\Tc$ is $\muxpw-$Fatou closed. Let $(\theta_n)_n$ be a sequence in $\Tc$ converging $\muxpw$ to $\theta\in\Theta_\mu$. Let $n\ge 1$, we denote $N_\mu$, the set in $\Nc_\mu$ from the proposition applied to $\theta_n$, and let $N_\mu^0\in\Nc_\mu$ corresponding to the $\muxpw$ convergence of $\theta_n$ to $\theta$. We denote $N_\mu:=\cup_{n\in\N}N_\mu^n\in\Nc_\mu$. Let $x_1,x_2\notin N_\mu$, and $\yb\in\dom_{x_1}\theta\cap\dom_{x_2}\theta$. Let $y_1,y_2\in\dom_{x_1}\theta$, such that we have the convex combination $\yb = \lambda y_1 + (1-\lambda) y_2$, and $0\le\lambda\le 1$. Then for $i=1,2$, $\theta_n(x_1,y_i)\longrightarrow\theta(x_1,y_i)$, and $\theta_n(x_1,\yb)\longrightarrow\theta(x_1,\yb)$, as $n\to\infty$. Using the fact that $\theta_n\in\Tc$, for all $n$, we have
\begin{equation}\label{eq:convinvn}
 \Delta_n
 :=
 \lambda \theta_n(x_i,y_1)+(1-\lambda)\theta_n(x_i,y_2)-\theta_n(x_i,\yb)
 \ge
 0,
 ~\mbox{and independent of}~
 i=1,2.
 \end{equation}
Taking the limit $n\to\infty$ gives that $\underline\theta_\infty(x_2,y_i)<\infty$, and $y_i\in \dom\underline\theta_\infty(x_2,\cdot)$. $\yb$ is interior to $\dom_{x_1}\theta$, then for any $y\in\dom_{x_1}\theta$, $y':= \yb + \frac{\epsilon}{1-\epsilon}(\yb-y)\in \dom_{x_1}\theta$ for $0<\epsilon<1$ small enough. Then $\yb = \epsilon y+(1-\epsilon)y'$. As we may chose any $y\in\dom_{x_1}\theta$, we have $\dom_{x_1}\theta\subset \dom\underline\theta_\infty(x_2,\cdot)$. Then, we have
 \begin{equation}\label{eq:inclusion}
 \rf_{x_2}\conv(\dom_{x_1}\theta\cup\dom_{x_2}\theta)
 \subset
 \rf_{x_2}\conv\,\dom\big(\underline\theta_\infty(x_2,\cdot)
                     \big) 
 = 
 \dom_{x_2}\theta.
\end{equation}
By Lemma \ref{lemma:intersectinterior}, as $\dom_{x_1}\theta\cap\dom_{x_2}\theta\neq \emptyset$, $\conv(\dom_{x_1}\theta\cup\dom_{x_2}\theta)=\ri\,\conv(\dom_{x_1}\theta\cup\dom_{x_2}\theta)$. In particular, $\conv(\dom_{x_1}\theta\cup\dom_{x_2}\theta)$ is relatively open and contains $x_2$, and therefore $\rf_{x_2}\conv(\dom_{x_1}\theta\cup\dom_{x_2}\theta) = \conv(\dom_{x_1}\theta\cup\dom_{x_2}\theta)$. Finally, by $\eqref{eq:inclusion}$, $\dom_{x_1}\theta\subset\dom_{x_2}\theta$. As there is a symmetry between $x_1$, and $x_2$, we have $\dom_{x_1}\theta = \dom_{x_2}\theta$. Then we may go to the limit in equation \eqref{eq:convinvn}:
\begin{equation}\label{eq:convinv}
 \Delta_\infty
 :=
 \lambda \theta(x_i,y_1)+(1-\lambda)\theta(x_i,y_2)-\theta(x_i,\yb)
 \ge
 0,
 ~\mbox{and independent of}~
 i=1,2
 \end{equation}
Now, let $y_1,y_2\in\R^d$, such that we have the convex combination $\yb = \lambda y_1 + (1-\lambda) y_2$, and $0\le\lambda\le 1$. we have three cases to study.
\\
\no\underline{\it Case 1:} $y_i\notin\cl\dom_{x_1}\theta$ for some $i=1,2$.
Then, as the average $\yb$ of the $y_i$ is in $\dom_{x_1}\theta$, by Proposition \ref{prop:ri} (ii), me may find $i'=1,2$ such that $y_{i'}\notin\conv\,\dom\theta(x_1,\cdot)$, thus implying that $\theta(x_1,y_i)=\infty$. Then $\lambda\theta(x_1,y_1)+(1-\lambda) \theta(x_1,y_2)-\theta(x_1,\yb) = \infty\ge 0$. As $\dom_{x_1}\theta = \dom_{x_2}\theta$, we may apply the same reasoning to $x_2$, we get $\lambda\theta(x_1,y_1)+(1-\lambda) \theta(x_2,y_2)-\theta(x_2,\yb) = \infty\ge 0$. We get the result.
\\
\no\underline{\it Case 2:} $y_1,y_2\in\dom_{x_1}\theta$.
This case is \eqref{eq:convinv}.
\\
\no\underline{\it Case 3:} $y_1,y_2\in\cl\dom_{x_1}\theta$.
The problem arises here if some $y_i$ is in the boundary $\partial \dom_{x_1}\theta$. Let $x\notin N_\mu$, we denote the lower semi-continuous envelop of $\theta(x,\cdot)$ in $\cl\dom_x\theta$, by $\theta_*(x,y) := \lim_{\epsilon\searrow 0}\theta(x,\epsilon x+(1-\epsilon)y')$, for $y\in\cl\dom_{x}\theta$, where the latest equality follows from Lemma \ref{lemma:convexlsc} together with that fact that $\theta(x,\cdot)$ is convex on $\dom_{x}\theta$. Let $y\in\cl\dom_{x_1}\theta$, for $1\ge\epsilon>0$, $y^\epsilon :=\epsilon x_1 + (1-\epsilon )y\in\dom_{x_1}\theta$. By \eqref{eq:convinvn}, $(1-\epsilon )\theta_n(x_1,y)-\theta_n(x_1,y^\epsilon) = (1-\epsilon )\theta_n(x_2,y)-\theta_n(x_2,y^\epsilon)$. Taking the $\liminf$, we have $(1-\epsilon )\theta(x_1,y)-\theta(x_1,y^\epsilon) = (1-\epsilon )\theta(x_2,y)-\theta(x_2,y^\epsilon)$. Now taking $\epsilon\searrow 0$, we have $\theta(x_1,y)-\theta_*(x_1,y) = \theta(x_2,y)-\theta_*(x_2,y)$. Then the jump of $\theta(x,\cdot)$ in $y$ is independent of $x=x_1$ or $x_2$. Now for $1\ge\epsilon>0$, by \eqref{eq:convinv}
$$
 \lambda \theta(x_1,y_1^\epsilon)+(1-\lambda)\theta(x_1,y_2^\epsilon)-\theta(x_1,\yb^\epsilon)
 =
\lambda \theta(x_2,y_1^\epsilon)+(1-\lambda)\theta(x_2,y_2^\epsilon)-\theta(x_2,\yb^\epsilon)
 \ge 0.
 $$
 By going to the limit $\epsilon\searrow 0$, we get
$$
 \lambda \theta_*(x_1,y_1)+(1-\lambda)\theta_*(x_1,y_2)-\theta_*(x_1,\yb)
 =
\lambda \theta_*(x_2,y_1)+(1-\lambda)\theta_*(x_2,y_2)-\theta_*(x_2,\yb)
 \ge 0.
$$
As the (nonnegative) jumps do not depend on $x= x_1$ or $x_2$, we finally get
$$
 \lambda \theta(x_1,y_1)+(1-\lambda)\theta(x_1,y_2)-\theta(x_1,\yb)
 =
\lambda \theta(x_2,y_1)+(1-\lambda)\theta(x_2,y_2)-\theta(x_2,\yb)
 \ge 0.
$$
 Finally, $\Tc$ is $\muxpw-$Fatou closed, and convex. $\widehat{\Tc}_1\subset\Tc$. As the result is clearly invariant when the function is multiplied by a scalar, the Result is proved on $\widehat{\Tc}(\mu,\nu)$.
\ep

\subsection{Compactness}\label{subsect:Komlos}
~~~
\\
~~
\\
\no {\bf Proof of Proposition \ref{prop:Komlos}}
We first prove the result for $\theta=(\theta_n)_{n\ge 1}\subset \Theta$. Denote $\conv(\theta) := \{\theta'\in\Theta^\N: \theta_n'\in \conv(\theta_k,k\ge n), n\in\N\}$. Consider the minimization problem:
 \be\label{min:muGdom}
 m
 &:=&
 \inf_{\theta'\in \conv(\theta)}
 \mu[G(\dom_X\underline\theta'_\infty)],
 \ee
where the measurability of $G(\dom_X\underline\theta'_\infty)$ follows from Lemma \ref{lemma:measurabilityI}.
\\
\no \underline{\it Step 1:} We first prove the existence of a minimizer. Let $(\theta'^k)_{k\in\N}\in \conv(\theta)^\N$ be a minimizing sequence, and define the sequence $\widehat\theta\in \conv(\theta)$ by:
 \b*
 &\widehat\theta_n 
 := 
 (1-2^{-n})^{-1}\sum_{k=1}^n 2^{-k}\theta'^k_n,
 ~n\ge 1.
 \e*
Then, $\dom(\underline{\widehat\theta}_\infty)\subset \bigcap_{k\ge 1}\dom(\underline{\theta}'^k_\infty)$ by the non-negativity of $\theta'$, and we have the inclusion $\big\{\widehat\theta_n\underset{n\to\infty}{\longrightarrow}\infty\big\}\subset\big\{\theta_n'^k\underset{n\to\infty}{\longrightarrow}\infty~\mbox{for some}~k\ge 1\big\}$. Consequently,
 \b*
 &\dom_x\underline{\widehat\theta}_\infty
 \;\subset\; 
 \conv\big(\bigcap_{k\ge 1}\dom\underline{\theta}'^k_\infty(x,\cdot)\big)
 \;\subset\;
 \bigcap_{k\ge 1}\dom_x\underline\theta'^k_\infty
 &
 \mbox{for all}~~
 x\in\R^d.
 \e*
Since $(\theta'^k)_k$ is a minimizing sequence, and $\widehat\theta\in \conv(\theta)$, this implies that $\mu[G(\dom_X\underline{\widehat\theta}_\infty)]= m$.
\\
\no \underline{\it Step 2:} We next prove that we may find a sequence $(y_i)_{i\ge 1}\subset \Leb^0(\R^d,\R^d)$ such that 
 \be\label{step2compactness}
 y_i(X)\in \aff(\dom_X\underline{\widehat\theta}_\infty)
 &\mbox{and}& 
 (y_i(X))_{i\ge 1}~\mbox{dense in}~\aff\dom_X\underline{\widehat\theta}_\infty,
 ~~\mu-\mbox{a.s.}
 \ee
Indeed, it follows from Lemmas \ref{lemma:measurabilityI}, and \ref{lemma:bimesconv} that the map $x\mapsto \aff(\dom_x\underline{\widehat\theta}_\infty)$ is universally measurable, and therefore Borel-measurable up to a modification on a $\mu-$null set. Since its values are closed and nonempty, we deduce from the implication $(ii) \implies (ix)$ in Theorem 4.2 of the survey on measurable selection \cite{wagner1977survey} the existence of a sequence $(y_i)_{i\ge 1}$ satisfying \eqref{step2compactness}.
\\
\no \underline{\it Step 3:} Let $m(dx,dy):=\mu(dx)\otimes\sum_{i\ge 0}2^{-i}\delta_{\{y_i(x)\}}(dy)$. By the Koml\`os lemma (in the form of Lemma A1.1 in \cite{delbaen1994general}, similar to the one used in the proof of Proposition 5.2 in \cite{beiglbock2015complete}), we may find $\widetilde\theta\in\conv(\widehat\theta)$ such that $\widetilde\theta_n\longrightarrow\widetilde\theta_\infty\in\L^0(\Omega)$, $m-$a.s. Clearly, $\dom_x\underline{\widetilde\theta}_\infty\subset\dom_x\underline{\widehat\theta}_\infty$, and therefore $\mu\big[G(\dom_X\underline{\widetilde\theta}_\infty)\big]\le \mu\big[G(\dom_x\underline{\widehat\theta}_\infty)\big]$, for all $x\in\R^d$. This shows that 
 \be\label{thetatildeoptimal}
 G(\dom_X\underline{\widetilde\theta}_\infty)
 \;=\;
 G(\dom_X\underline{\widehat\theta}_\infty),~~\mu-\mbox{a.s.}
 \ee
so that $\widetilde\theta$ is also a solution of the minimization problem \eqref{min:muGdom}. Moreover, it follows from \eqref{eq:G} that
 \b*
 \ri\,\dom_X\underline{\widetilde\theta}_\infty
 =
 \ri\,\dom_X\underline{\widehat\theta}_\infty,
 &\mbox{and therefore}&
 \aff\,\dom_X\underline{\widetilde\theta}_\infty
 =
 \aff\,\dom_X\underline{\widehat\theta}_\infty,
 ~~\mu-a.s.
 \e*
\no \underline{\it Step 4:} Notice that the values taken by $\widetilde\theta_\infty$ are only fixed on an $m-$full measure set. By the convexity of elements of $\Theta$ in the $y-$variable, $\dom_X\widetilde\theta_n$ has a nonempty interior in $\aff(\dom_X\underline{\widetilde\theta}_\infty)$. Then as $\mu-$a.s., $\widetilde\theta_n(X,\cdot)$ is convex, the following definition extends $\widetilde\theta_\infty$ to $\Omega$:
 \b*
 \widetilde\theta_\infty(x,y)
 :=
 \sup\big\{ a\cdot y + b: (a,b)\in\R^d\times\R,~a\cdot y_n(x)+b\le \widetilde\theta_\infty(x,y_n(x))~\mbox{for all}~n\ge 0
     \big\}.
 \e*
This extension coincides with $\widetilde\theta_\infty$, in $(x,y_n(x))$ for $\mu-$a.e. $x\in\R^d$, and all $n\ge 1$ such that $y_n(x)\notin \partial\dom_X\widetilde\theta_k$ for some $k\ge 1$ such that $\dom_x\widetilde\theta_n$ has a nonempty interior in $\aff(\dom_x\underline{\widetilde\theta}_\infty)$. As for $k$ large enough, $\partial\dom_X\widetilde\theta_k$ is Lebesgue negligible in $\aff(\dom_x\underline{\widetilde\theta}_\infty)$, the remaining $y_n(x)$ are still dense in $\aff(\dom_x\underline{\widetilde\theta}_\infty)$. Then, for $\mu-$a.e. $x\in\R^d$, $\widetilde\theta_n(x,\cdot)$ converges to $\widetilde\theta_\infty(x,\cdot)$ on a dense subset of $\aff(\dom_x\underline{\widetilde\theta}_\infty)$. We shall prove in Step 6 below that 
 \be\label{step5:nonemptyinterior}
 \dom\,\widetilde\theta_\infty(X,\cdot)
 ~\mbox{has nonempty interior in}~
 \aff(\dom_X\underline{\widetilde\theta}_\infty),
 &\mu-\mbox{a.s.}
 \ee
Then, by Theorem \ref{thm:ConvexConv}, $\widetilde\theta_n(X,\cdot)\longrightarrow\widetilde\theta_\infty(X,\cdot)$ pointwise on $\aff(\dom_X\underline{\widetilde\theta}_\infty)\setminus\partial\dom\widetilde\theta_\infty(X,\cdot)$, $\mu-$a.s. Since $\dom_X\theta_\infty = \dom_X\underline{\theta}_\infty$, and $\widetilde\theta$ converges to $\theta_\infty$ on $\dom_X\theta_\infty$, $\mu-$a.s., $\widetilde\theta$ converges to $\theta_\infty\in\Theta$, $\muxpw$.
\\
\no \underline{\it Step 5:} Finally for general $(\theta_n)_{n\ge 1}\subset \Theta_\mu$, we consider $\theta_n'$, equal to $\theta_n$, $\muxpw$, such that $\theta_n'\le\theta_n$, for $n\ge 1$, from the definition of $\Theta_\mu$. Then we may find $\lambda_n^k$, coefficients such that $\widehat\theta_n':=\sum_{k\ge n}\lambda_n^k\theta_k'\in \conv(\theta')$ converges $\muxpw$ to $\widehat\theta_\infty\in\Theta$. We denote $\widehat\theta_n:=\sum_{k\ge n}\lambda_n^k\theta_k\in \conv(\theta)$, $\widehat\theta_n = \widehat\theta_n'$, $\muxpw$, and $\widehat\theta_n \ge \widehat\theta_n'$. By Proposition \ref{prop:convergenceprop} (iii), $\widehat\theta$ converges to $\widehat\theta_\infty$, $\muxpw$. The Proposition is proved.

\no \underline{\it Step 6:} In order to prove \eqref{step5:nonemptyinterior}, suppose to the contrary that there is a set $A$ such that $\mu[A]>0$ and $\dom\widetilde\theta_\infty(x,\cdot)$ has an empty interior in $\aff(\dom_x\underline{\widetilde\theta}_\infty)$ for all $x\in A$. Then, by the density of the sequence $(y_n(x))_{n\ge 1}$ stated in \eqref{step2compactness}, we may find for all $x\in A$ an index $i(x)\ge 0$ such that 
 \be\label{contra:step5}
 \widehat y(x):=y_{i(x)}(x)\in \ri\,\dom_x\underline{\widetilde\theta}_\infty,
 &\mbox{and}&
 \widetilde\theta_\infty(x,\widehat y(x))=\infty.
 \ee
Moreover, since $i(x)$ takes values in $\N$, we may reduce to the case where $i(x)$ is a constant integer, by possibly shrinking the set $A$, thus guaranteeing that $\widehat y$ is measurable. Define the measurable function on $\Omega$:
 \be\label{Lnx}
 \theta^0_n(x,y)
 :=
 \dist(y,L^n_x)
 &\mbox{with}&
 L^n_x:=\big\{y\in\R^d:\widetilde\theta_n(x,y)<\widetilde\theta_n(x,\widehat y(x))\big\}.
 \ee
Since $L^n_x$ is convex, and contains $x$ for $n$ sufficiently large by \eqref{contra:step5}, we see that 
 \be\label{theta0}
 \theta^0_n~~\mbox{is convex in $y$ and} 
 &\theta^0_n(x,y)\leq |x-y|,& 
 \mbox{for all}~~(x,y)\in\Omega.
 \ee
In particular, this shows that $\theta^0_n\in\Theta$. By Koml\`os Lemma, we may find
 \b*
 &\widehat\theta^0_n:=\sum_{k\ge n}\lambda^n_k \theta^0_k\in\conv(\theta^0)&
 \mbox{such that}~~~
 \widehat\theta^0_n\longrightarrow \widehat\theta^0_\infty,~~m-\mbox{a.s.}
 \e*
for some non-negative coefficients $(\lambda^n_k,k\ge n)_{n\ge 1}$ with $\sum_{k\ge n}\lambda^n_k=1$. By convenient extension of this limit, we may assume that $\widehat\theta^0_\infty\in\Theta$. We claim that 
 \begin{equation}\label{claim:step5}
 \widehat\theta^0_\infty>0
 ~\mbox{on}~
 H_x:=\{h(x)\cdot(y-\widehat y(x))>0\},
 ~\mbox{for some}~
 h(x)\in\R^d.
 \end{equation}
We defer the proof of this claim to Step 7 below and we continue in view of the required contradiction. By definition of $\theta^0_n$ together with \eqref{theta0}, we compute that
 \b*
 \theta^1_n(x,y)
 \;:=\;
 \sum_{k\ge n}\lambda^n_k\widetilde\theta_k(x,y)
 &\ge&
 \sum_{k\ge n}\lambda^n_k\widetilde\theta_k(x,\widehat y(x)) \1_{\{\theta^0_n>0\}}
 \\
 &\ge&
 \sum_{k\ge n}\lambda^n_k\widetilde\theta_k(x,\widehat y(x)) \frac{\theta^0_k(x,y)}{|x-y|}
 \\
 &\ge&
 \frac{\widehat\theta^0_n(x,y)}{|x-y|}\inf_{k\ge n}\widetilde\theta_k(x,\widehat y(x)).
 \e*
By \eqref{contra:step5} and \eqref{claim:step5}, this shows that the sequence $\theta^1\in\conv(\theta)$ satisfies
 \b*
 \theta^1_n(x,\cdot) \longrightarrow \infty,
 &\mbox{on}&
 H_x,
 ~~\mbox{for all}~~
x\in A.
 \e*
We finally consider the sequence $\widetilde\theta^1:=\frac12(\widetilde\theta+\theta^1)\in\conv(\theta)$. Clearly, $\dom\underline{\widetilde\theta}^1_\infty(X,\cdot)\subset\dom\underline{\widetilde\theta}_\infty(X,\cdot)$, and it follows from the last property of $\theta^1$ that $\dom\underline{\widetilde\theta}^1_\infty(x,\cdot)\subset H_x^c\cap\dom\underline{\widetilde\theta}_\infty(x,\cdot)$ for all $x\in A$. Notice that $\widehat y(x)$ lies on the boundary of the half space $H_x$ and, by \eqref{contra:step5}, $\widehat y(x)\in\ri\dom_x\underline{\widetilde\theta}_\infty$. Then $G(\dom_x\underline{\widetilde\theta}^1_\infty)<
G(\dom_x\underline{\widetilde\theta}_\infty)$ for all $x\in A$ and, since $\mu[A]>0$, we deduce that $\mu\big[G(\dom_X\underline{\widetilde\theta}^1_\infty)\big]<
 \mu\big[G(\dom_X\underline{\widetilde\theta}_\infty)\big]$, contradicting the optimality of $\widetilde\theta$, by \eqref{thetatildeoptimal}, for the minimization problem \eqref{min:muGdom}.
\\
\no \underline{\it Step 7:} It remains to justify \eqref{claim:step5}. Since $\widetilde\theta_n(x,\cdot)$ is convex, it follows from the Hahn-Banach separation theorem that:
 \b*
 \widetilde\theta_n(x,\cdot) \ge \widetilde\theta_n(x,\widehat y(x))
 ~\mbox{on}~
 H^n_x := \big\{y\in\R^d: h^n(x)\cdot(y-\widehat y(x))>0\big\},
 \e*
for some $h^n(x)\in\R^d$, so that it follows from \eqref{Lnx} that $L^n_x\subset (H^n_x)^c$, and
 \b*
 \theta^0_n(x,y)
 \ge
 \dist\big(y,(H^n_x)^c\big)
 =
 \big[\big(y-\widehat y(x)\big)\cdot h^n(x)\big]^+.
 \e*
Denote $g_x:=g_{\dom_x\underline{\widehat\theta}_\infty}$ the Gaussian kernel restricted to the affine span of $\dom_x\underline{\widehat\theta}_\infty$, and $B_r(x_0)$ the corresponding ball with radius $r$, centered at some point $x_0$. By \eqref{contra:step5}, we may find $r^x$ so that $B^x_r:=B_r(\widehat y(x))\subset\ri\,\dom_x\underline{\widetilde\theta}_\infty$ for all $r\le r^x$, and
 \b*
 \int_{B_r^x} \theta^0_n(x,y)g_x(y)dy
 &\ge &
 \int_{B_r^x} \big[\big(y-\widehat y(x)\big)\cdot h^n(x)\big]^+g_x(y)dy
 \\
 &\ge&
 \min_{B_r^x}g_x\int_{B_r(0)} (y\cdot e_1)^+ dy
 =:b^r_x >0,
 \e*
where $e_1$ is an arbitrary unit vector of the affine span of $\dom_x\underline{\widehat\theta}_\infty$. Then we have the inequality $\int_{B_r^x} \widehat\theta^0_n(x,y)g_x(y)dy
 \ge b^r_x$, and since $\widehat\theta^0_n$ has linear growth in $y$ by \eqref{theta0}, it follows from the dominated convergence theorem that $\int_{B^r_x} \widehat\theta^0_\infty(x,y)g(dy)\ge b^r_x>0$, and therefore $\widehat\theta^0_\infty(x,y^r_x)>0$ for some $y^r_x\in B^r_x$. From the arbitrariness of $r\in(0,r_x)$, We deduce \eqref{claim:step5} as a consequence of the convexity of $\widehat\theta^0(x,\cdot)$. 
 \ep\\

\no {\bf Proof of Proposition \ref{prop:convergenceprop}} (iii)
 We need to prove the existence of some
 \be\label{thetaprime}
 \theta'\in\Theta
 &\mbox{such that}&
 \underline{\theta}_\infty= \theta',~~\muxpw,~~
 \mbox{and}~~
 \underline{\theta}_\infty\ge \theta'.
 \ee
For simplicity, we denote $\theta:=\underline{\theta}_\infty$. Let 
 \b*
 &F^1 := \cl\,\conv\,\dom{{\theta}}(X,\cdot),
 ~~ 
 F^k 
 := 
 \cl\,\conv\big(\dom{{\theta}}(X,\cdot)\cap \aff\,\rf_XF^{k-1}\big),
 ~k \ge 2,&
 \\
& ~~\mbox{and}~~
 F := \cup_{n\ge 1}(F^n\setminus \cl\rf_X F^{n})\cup\cl\dom_X{{\theta}}.
 &
 \e*
Fix some sequence $\eps_n\searrow 0$, and denote $\theta_*:=\liminf_{n\to\infty}{{\theta}}\big(X,\eps_n X+(1-\eps_n)Y\big)$, and
 $$
 \theta':= \big[\infty\mathbf{1}_{Y\notin F(X)}+\mathbf{1}_{Y\in \cl\dom_X{{\theta}}}\theta_*\big]\mathbf{1}_{X\notin N_\mu},
 $$
where $N_\mu\in\Nc_\mu$ is chosen such that $\mathbf{1}_{Y\in F^{k}(X)}\mathbf{1}_{X\notin N_\mu}$ are Borel measurable for all $k$ from Lemma \ref{lemma:measurabilityFrontier}, and ${\theta}(x,\cdot)$ (resp. $\theta_n(x,\cdot)$) is convex finite on $\dom_x{\theta}$ (resp. $\dom_x\theta_n$), for $x\notin N_\mu$. Consequently, $\theta'$ is measurable. In the following steps, we verify that $\theta'$ satisfies \eqref{thetaprime}.

\no\underline{\it Step 1:} We prove that $\theta'\in\Theta$. Indeed, $\theta'\in \Leb^0_+(\Omega)$, and $\theta'(X,X) = 0$. Now we prove that $\theta'(x,\cdot)$ is convex for all $x\in\R^d$. For $x\in N_\mu$, $\theta'(x,\cdot)=0$. For $x\notin N_\mu$, ${\theta}(x,\cdot)$ is convex finite on $\dom_x{\theta}$, then by the fact that $\dom_x{\theta}$ is a convex relatively open set containing $x$, it follows from Lemma \ref{lemma:convexlsc} that $\theta_*(x,\cdot)=\lim_{n\to\infty}{{\theta}}\big(x,\eps_n x+(1-\eps_n)\cdot\big)$ is the lower semi-continuous envelop of ${\theta}(x,\cdot)$ on $\cl\dom_x{\theta}$. We now prove the convexity of $\theta'(x,\cdot)$ on all $\R^d$. We denote $\widehat F(x) := F(x)\setminus \cl\dom_x{\theta}$ so that $\R^d = F(x)^c\cup\widehat F(x) \cup \cl\dom_x{\theta}$. Now, let $y_1,y_2\in\R^d$, and $\lambda\in(0,1)$. If $y_1\in F(x)^c$, the convexity inequality is verified as $\theta'(x,y_1)=\infty$. Moreover, $\theta'(x,\cdot)$ is constant on $\widehat F(x)$, and convex on $\cl\dom_x{\theta}$. We shall prove in Steps 4 and 5 below that
 \be\label{Step45}
 F(x)
 ~\mbox{is convex, and}~
 \rf_xF(x) = \dom_x{\theta}.
 \ee
In view of Proposition \ref{prop:ri} (ii), this implies that the sets $\widehat F(x)$ and $\cl\dom_x{\theta}$ are convex. Then we only need to consider the case when $y_1\in\widehat{F}(x)$, and $y_2\in\cl\dom_x{\theta}$. By Proposition \ref{prop:ri} (ii) again, we have $[y_1,y_2)\subset \widehat{F}(x)$, and therefore $\lambda y_1 + (1-\lambda)y_2\in\widehat{F}(x)$, and $\theta'(x,\lambda y_1 + (1-\lambda)y_2) = 0$, which guarantees the convexity inequality.

\no\underline{\it Step 2:} We next prove that ${\theta}=\theta'$, $\muxpw$. By the second claim in \eqref{Step45}, it follows that $\theta_*(X,\cdot)$ is convex finite on $\dom_X{\theta}$, $\mu-$a.s. Then as a consequence of Proposition \ref{prop:Thetamu} (ii), we have $\dom_X\theta' = \dom_X(\infty\mathbf{1}_{Y\notin F(X)})\cap \dom_X(\theta_*\mathbf{1}_{Y\in \cl\dom_X{{\theta}}})$, $\mu-$a.s. The first term in this intersection is $\rf_XF(X) = \dom_X{{\theta}}$. The second contains $\dom_X{{\theta}}$, as it is the $\dom_X$ of a function which is finite on $\dom_X{{\theta}}$, which is convex relatively open, containing $X$. Finally, we proved that $\dom_X{{\theta}}=\dom_X\theta'$, $\mu-$a.s. Then $\theta'(X,\cdot)$ is equal to $\theta_*(X,\cdot)$ on $\dom_X{{\theta}}$, and therefore, equal to ${\theta}(X,\cdot)$, $\mu-$a.s. We proved that ${\theta}=\theta'$, $\muxpw$.

\no\underline{\it Step 3:} We finally prove that $\theta'\le{\theta}$ pointwise. We shall prove in Step 6 below that 
 \be\label{Step6}
 \dom{{\theta}}(X,\cdot)
 &\subset& F.
 \ee
Then, $\infty\mathbf{1}_{Y\notin F(X)}\mathbf{1}_{X\notin N_\mu}\le {{\theta}}$, and it remains to prove that
 \b*
 {\theta}(x,y)
 \ge  
 \theta_*(x,y)
 &\mbox{for all}&
 y\in\cl\dom_x\theta,~~x\notin N_\mu.
 \e*
To see this, let $x\notin N_\mu$. By definition of $N_\mu$, $\theta_n(x,\cdot)\longrightarrow\theta(x,\cdot)$ on $\dom_x\theta$. Notice that ${\theta}(x,\cdot)$ is convex on $\dom_x\theta$, and therefore as a consequence of Lemma \ref{lemma:convexlsc},
\b*
\theta_*(x,y)
=
\lim_{\epsilon\searrow 0}{\theta}\big(x,\epsilon x + (1-\epsilon)y\big),&\mbox{for all}&y\in\cl\dom_x\theta.
\e*
Then $y^\epsilon := (1-\epsilon)y + \epsilon x\in \dom_x\theta_n$, for $\eps\in(0,1]$, and $n$ sufficiently large by (i) of this Proposition, and therefore $(1-\epsilon)\theta_n(x,y)-\theta_n(x,y_\epsilon)\ge (1-\epsilon)\theta_n'(x,y)-\theta_n'(x,y_\epsilon)\ge 0$, for $\theta_n'\in\Theta$ such that $\theta_n'=\theta_n$, $\muxpw$, and $\theta_n\ge \theta_n'$. Taking the $\liminf$ as $n\to\infty$, we get $(1-\epsilon){\theta}(x,y)-\theta(x,y_\epsilon)\ge 0$, and finally $
 {\theta}(x,y)
 \ge 
 \lim_{\epsilon\searrow 0}{\theta}\big(x,\epsilon x + (1-\epsilon)y\big) 
 = 
 \theta'(x,y),
 $ by sending $\epsilon\searrow 0$.

\no\underline{\it Step 4:} (First claim in \eqref{Step45}) Let $x_0\in\R^d$, let us prove that $F(x_0)$ is convex. Indeed, let $x,y\in F(x_0)$, and $0<\lambda<1$. Since $\cl\dom_x{\theta}$ is convex, and $F^n(x_0)\setminus \cl\rf_X F^{n}(x_0)$ is convex by Proposition \ref{prop:ri} (ii), we only examine the following non-obvious cases:

$\bullet$ Suppose $x\in F^n(x_0)\setminus \cl\rf_{x_0} F^{n}(x_0)$, and $y\in F^p(x_0)\setminus \cl\rf_{x_0} F^{p}(x_0)$, with $n<p$. Then as $F^p(x_0)\setminus \cl\rf_{x_0} F^{p}(x_0)\subset \cl\rf_{x_0} F^{n}(x_0)$, we have $\lambda x + (1-\lambda)y\in F^n(x_0)\setminus \cl\rf_{x_0} F^{n}(x_0)$ by Proposition \ref{prop:ri} (ii).

$\bullet$ Suppose $x\in F^n(x_0)\setminus \cl\rf_{x_0} F^{n}(x_0)$, and $y\in \cl\dom_{x_0}{{\theta}}$, then as $\cl\dom_{x_0}{{\theta}}\subset\cl\rf_{x_0} F^{n}(x_0)$, this case is handled similar to previous case.

\no\underline{\it Step 5:} (Second claim in \eqref{Step45}). We have $\dom_X{{\theta}}\subset F(X)$, and therefore $\dom_X{{\theta}}\subset \rf_XF(X)$. Now we prove by induction on $k\ge 1$ that $\rf_XF(X)\subset \cup_{n\ge k}(F^n\setminus \cl\rf_X F^{n})\cup\cl\dom_X{{\theta}}$. The inclusion is trivially true for $k=1$. Let $k\ge 1$, we suppose that the inclusions holds for $k$, hence $\rf_XF(X)\subset \cup_{n\ge k}(F^n\setminus \cl\rf_X F^{n})\cup\cl\dom_X{{\theta}}$. As $\cup_{n\ge k}(F^n\setminus \cl\rf_X F^{n})\cup\cl\dom_X{{\theta}}\subset F^k$. Applying $\rf_X$ gives
\b*
\rf_XF(X)&\subset& \rf_X\Big[\cup_{n\ge k}(F^n\setminus \cl\rf_X F^{n})\cup\cl\dom_X{{\theta}}\Big]\\
&=&\rf_X\Big[F^k\cap\cup_{n\ge k}(F^n\setminus \cl\rf_X F^{n})\cup\cl\dom_X{{\theta}}\Big]\\
&=&\rf_X F^k\cap\rf_X\Big[\cup_{n\ge k}(F^n\setminus \cl\rf_X F^{n})\cup\cl\dom_X{{\theta}}\Big]\\
&\subset&\cl\rf_X F^k\cap\cup_{n\ge k}(F^n\setminus \cl\rf_X F^{n})\cup\cl\dom_X{{\theta}}\\
&\subset&\cup_{n\ge k+1}(F^n\setminus \cl\rf_X F^{n})\cup\cl\dom_X{{\theta}}.
\e*
Then the result is proved for all $k$. In particular we apply it for $k=d+1$. Recall from the proof of Lemma \ref{lemma:measurabilityI} that for $n\ge d+1$, $F^n$ is stationary at the value $\cl\dom_X{{\theta}}$. Then $\cup_{n\ge d+1}(F^n\setminus \cl\rf_X F^{n})=\emptyset$, and $\rf_XF(X)\subset\rf_X\cl\dom_X{{\theta}} = \dom_X{{\theta}}$. The result is proved.

\no\underline{\it Step 6:} We finally prove \eqref{Step6}. Indeed, $\dom{{\theta}}(X,\cdot)\subset F^1$ by definition. Then
\b*
\dom{{\theta}}(X,\cdot)&\subset&F^1\setminus\aff F^1\cup\big(\cup_{2\le k\le d+1}(\dom{{\theta}}(X,\cdot)\cap \aff\,\rf_XF^{k-1})\setminus\aff F^{k}\big)\cup F^{d+1}\\
&\subset&F^1\setminus\cl F^1\cup\big(\cup_{k\ge 2}\,\cl\,\conv(\dom{{\theta}}(X,\cdot)\cap \aff\,\rf_XF^{k-1})\setminus\cl F^{k}\big)\cup\,\cl\dom_X{\theta}\\
&&
=\cup_{k\ge 1}F^{k}\setminus\cl F^{k}\cup\cl\dom_X{\theta}
=F.
\e*
 \ep

\section{Some convex analysis results}
\label{sect:convexanalysis}
\setcounter{equation}{0}

As a preparation, we first report a result on the union of intersecting relative interiors of convex subsets which was used in the proof of Proposition \ref{prop:partition}. We shall use the following characterization of the relative interior of a convex subset $K$ of $\R^d$:
 \be
 \ri K 
 &=& 
 \big\{x\in\R^d : x-\epsilon(x'-x)\in K\mbox{ for some }\epsilon>0,\mbox{ for all }x'\in K\big\}
 \label{characri-i}
 \\
 &=&
 \big\{x\in\R^d : x\in (x',x_0],\mbox{ for some }x_0\in\ri K,\mbox{ and }x'\in K \big\}.
 \label{characri-ii}
 \ee
 
We start by proving the required properties of the notion of relative face.\\

\no{\bf Proof of Proposition \ref{prop:ri}}
\no\underline{\rm (i)} The proofs of the first properties raise no difficulties and are left as an exercise for the reader. We only prove that $\rf_a A  = \ri A \neq\emptyset$ iff $a\in\ri A$. We assume that $\rf_a A  = \ri A\neq\emptyset$. The non-emptiness implies that $a\in A$, and therefore $a\in\rf_aA = \ri A$. Now we suppose that $a\in\ri A$. Then for $x\in\ri A$, $\big[x,a-\epsilon (x-a)\big]\subset \ri A\subset A$, for some $\epsilon>0$, and therefore $\ri A\subset\rf_a A$. On the other hand, by \eqref{characri-ii}, $\ri A = \{x\in\R^d : x\in (x',x_0],\mbox{ for some }x_0\in\ri A,\mbox{ and }x'\in A\}$. Taking $x_0 := a\in\ri A$, we have the remaining inclusion $\rf_a A\subset\ri A$.

\no\underline{\rm (ii)} We now assume that $A$ is convex.

\no\quad\underline{\it Step 1:} We first prove that $\rf_aA$ is convex. Let $x,y\in\rf_aA$ and $\lambda\in[0,1]$. We consider $\epsilon >0$ such that $\big(a-\epsilon(x-a),x+\epsilon(x-a)\big)\subset A$ and $\big(a-\epsilon(y-a),y+\epsilon(y-a)\big)\subset A$. Then if we write $z=\lambda x + (1-\lambda)y$, $\big(a-\epsilon(z-a),z+\epsilon(z-a)\big)\subset A$ by convexity of $A$, because $a,x,y\in A$.

\no\quad\underline{\it Step 2:} In order to prove that $\rf_aA$ is relatively open, we consider $x,y\in\rf_aA$, and we verify that $\big(x-\epsilon(y-x),y+\epsilon(y-x)\big)\subset \rf_aA$ for some $\epsilon>0$. Consider the two alternatives:

\no\quad\quad\underline{\it Case 1:} If $a,x,y$ are on a line. If $a=x=y$, then the required result is obvious. Otherwise,
$$\big(a-\epsilon(x-a),x+\epsilon(x-a)\big)\cup\big(a-\epsilon(y-a),y+\epsilon(y-a)\big)\subset \rf_a A$$
This union is open in the line and $x$ and $y$ are interior to it. We can find $\epsilon'>0$ such that $\big(x-\epsilon'(y-x),y+\epsilon'(y-x)\big)\subset \rf_aA$.

\no\quad\quad\underline{\it Case 2:} If $a,x,y$ are not on a line. Let $\epsilon>0$ be such that $\big(a-2\epsilon(x-a),x+2\epsilon(x-a)\big)\subset A$ and $\big(a-2\epsilon(y-a),y+2\epsilon(y-a)\big)\subset A$. Then $x+\epsilon(x-a)\in \rf_a A$ and $a-\epsilon(y-a)\in \rf_a A$. Then, if we take $\lambda := \frac{\epsilon}{1+2\epsilon}$,
$$
\lambda (a-\epsilon(y-a))+(1-\lambda)(x+\epsilon(x-a))
=
 (1-\lambda)(1+\epsilon)x-\lambda\epsilon y
=
 x + \lambda\epsilon (x-y)
$$
Then $x + \lambda\epsilon (x-y)\in \rf_a A$ and symmetrically, $y + \lambda\epsilon (y-x)\in \rf_a A$ by convexity of $\rf_aA$. And still by convexity, $\big(x-\epsilon'(y-x),y+\epsilon'(y-x)\big)\subset \rf_aA$ for  $\epsilon':=\frac{\epsilon^2}{1+2\epsilon}>0$.

\no\quad\underline{\it Step 3:} Now we prove that $A\setminus \cl \rf_a A$ is convex, and that if $x_0\in A\setminus \cl \rf_a A$ and $y_0\in A$, then $[x_0,y_0)\subset A\setminus \cl\rf_a A$. We will prove these two results by an induction on the dimension of the space $d$. First if $d = 0$ the results are trivial. Now we suppose that the result is proved for any $d'<d$, let us prove it for dimension $d$.

\no\quad\quad\underline{\it Case 1:} $a\in\ri A$. This case is trivial as $\rf_aA = \ri A$ and $A\subset \cl\ri A = \cl\rf_aA$ because of the convexity of $A$. Finally $A\setminus\cl\rf_aA = \emptyset$ which makes it trivial.

\no\quad\quad\underline{\it Case 2:} $a\notin\ri A$. Then $a\in \partial A$ and there exists a hyperplan support $H$ to $A$ in $a$ because of the convexity of $A$. We will write the equation of $E$, the corresponding half-space containing $A$, $E:c\cdot x \leq b$ with $c\in\R^d$ and $b\in\R$. As $x\in\rf_aA$ implies that $[a-\epsilon(x-a),x+\epsilon(x-a)]\subset A$ for some $\epsilon >0$, we have $(a-\epsilon(x-a))\cdot c \leq b$ and $(x+\epsilon(x-a))\cdot c \leq b$. These equations are equivalent using that $a\in H$ and thus $a\cdot c = b$ to $-\epsilon(x-a)\cdot c\leq 0$ and $(1+\epsilon)(x-a)\cdot c \leq 0$. We finally have $(x-a)\cdot c = 0$ and $x\in H$. We proved that $\rf_aA\subset H$.

Now using (i) together with the fact that $\rf_aA\subset H$ and $a\in H$ affine, we have
$$
\rf_a(A\cap H) = \rf_a A\cap \rf_a H = \rf_a A\cap H = \rf_aA.$$
Then we can now have the induction hypothesis on $A\cap H$ because $\dim H = d-1$ and $A\cap H\subset H$ is convex. Then we have $A\cap H \setminus \cl\rf_aA$ which is convex and if $x_0\in A\cap H\setminus \cl \rf_a (A\cap H)$, $y_0\in A\cap H$ and if $\lambda\in (0,1]$ then $\lambda x_0 + (1-\lambda)y_0\in A\setminus \cl\rf_a (A\cap H)$.

First $A\setminus \cl\rf_aA = (A\setminus H) \cup \big((A\cap H) \setminus \cl\rf_aA\big)$, let us show that this set is convex. The two sets in the union are convex ($A\setminus H = A\cap (E\setminus H)$), so we need to show that a non trivial convex combination of elements coming from both sets is still in the union. We consider $x\in A\setminus H$, $y\in A\cap H \setminus \cl\rf_aA$ and $\lambda >0$, let us show that $z:= \lambda x + (1-\lambda)y\in (A\setminus H) \cup (A\cap H \setminus \cl\rf_aA)$. As $x,y\in A$ ($\cl\rf_a A\subset A$ because $A$ is closed), $z\in A$ by convexity of $A$. We now prove $z\notin H$,
$$
z\cdot c
=
\lambda x\cdot c + (1-\lambda)y\cdot c
=
\lambda x\cdot c + (1-\lambda)b
<
\lambda b  + (1-\lambda)b = b.
$$
Then $z$ is in the strict half space: $z\in E\setminus H$. Finally $z\in A\setminus H$ and $A\setminus \cl\rf_aA$ is convex.\\
Let us now prove the second part: we consider $x_0\in A\setminus \cl \rf_a A$, $y_0\in\cl\rf_a A$ and $\lambda\in (0,1]$ and write $z_0 := \lambda x_0 + (1-\lambda)y_0$.

\no\quad\quad\underline{\it Case 2.1:} $x_0,y_0\in H$. We apply the induction hypothesis.

\no\quad\quad\underline{\it Case 2.2:} $x_0,y_0\in A\setminus H$. Impossible because $\rf_aA\subset H$ and $\cl\rf_aA\subset \cl H = H$. $y_0\in H$.

\no\quad\quad\underline{\it Case 2.3:} $x_0 \in A\setminus H$ and $y_0\in H$. Then by the same computation than in Step 1,
$$z_0\in A\setminus H\subset A\setminus\cl\rf_aA.$$

\no\quad\underline{\it Step 4:} Now we prove that if $a\in A$, then $\dim(\rf_a\cl A) =\dim(A)$ if and only if $a\in \ri A$, and that in this case, we have $\cl\rf_a\cl A = \cl\ri\,\cl A = \cl A = \cl\rf_aA$. We first assume that $a\in\ri A$. As by the convexity of $A$, $\ri A = \ri\,\cl A$, $\rf_a\cl A = \ri\,\cl A$, and therefore $\cl\,\rf_a\cl A = \cl A$. Finally, taking the dimension, we have $\dim(\cl\rf_a\cl A) =\dim(A)$. In this case we proved as well that $\cl\rf_a\cl A = \cl\ri\,\cl A = \cl A = \cl\rf_aA$, the last equality coming from the fact that $\ri A=\rf_aA$ as $a\in\ri A$.

Now we assume that $a\notin\ri A$. Then $a\in\partial \cl A$, and $\rf_a\cl A\subset\partial \cl A$. Taking the dimension (in the local sense this time), and by the fact that $\dim\partial\cl A=\dim\partial A<\dim A$, we have $\dim(\cl\rf_a\cl A) <\dim(A)$ (as $\cl\rf_a\cl A$ is convex, the two notions of dimension coincide).
\ep

\begin{Lemma}\label{lemma:intersectinterior}
Let $K_1,K_2\subset \R^d$ be convex with $\ri K_1\cap\ri K_2\neq\emptyset$. Then $\conv(\ri K_1\cup\ri K_2) = \ri\,\conv(K_1\cup K_2)$.
\end{Lemma}

\proof We fix $y\in \ri K_1\cap\ri K_2$.

Let $x\in\conv(\ri K_1\cup\ri K_2)$, we may write $x=\lambda x_1+(1-\lambda)x_2$, with $x_1\in\ri K_1$, $x_2\in\ri K_2$, and $0\le\lambda\le 1$. If $\lambda$ is $0$ or $1$, we have trivially that $x\in\ri\,\conv(K_1\cap K_2)$. Let us now treat the case $0<\lambda<1$. Then for $x'\in \conv(K_1\cup K_2)$, we may write $x'=\lambda' x_1'+(1-\lambda')x_2'$, with $x_1'\in K_1$, $x_2'\in K_2$, and $0\le\lambda'\le 1$. We will use $y$ as a center as it is in both the sets. For all the variables, we add a bar on it when we subtract $y$, for example $\xb := x-y$. The geometric problem is the same when translated with $y$,
\begin{equation}\label{eq:decompositionconvex}
\xb-\epsilon(\xb'-\xb)
\;=\;
\lambda \left(\xb_1-\epsilon\left(\frac{\lambda'}{\lambda} \xb_1'- \xb_1\right)\right) + (1-\lambda)\left(\xb_2-\epsilon\left(\frac{1-\lambda'}{1-\lambda} \xb_2'-\xb_2\right)\right).
\end{equation}
However, as $\xb_1$ and $\xb_1'$ are in $K_1-y$, as $0$ is an interior point, $\epsilon(\frac{\lambda'}{\lambda} \xb_1'- \xb_1)\in K_1-y$ for $\epsilon$ small enough. Then as $\xb_1$ is interior to $K_1-y$ as well, $\xb_1-\epsilon(\frac{\lambda'}{\lambda} \xb_1'- \xb_1)\in K_1-y$ as well. By the same reasoning, $\xb_2-\epsilon(\frac{1-\lambda'}{1-\lambda} \xb_2'-\xb_2)\in K_2-y$. Finally, by \eqref{eq:decompositionconvex}, for $\epsilon$ small enough, $x-\epsilon(x'-x)\in \conv(K_1\cup K_2)$. By \eqref{characri-i}, $x\in\ri\,\conv(K_1\cup K_2)$.

Now let $x\in \ri\,\conv(K_1\cup K_2)$. We use again $y$ as an origin with the notation $\xb := x-y$. As $\xb$ is interior, we may find $\epsilon >0$ such that $(1+\epsilon)\xb\in \conv(K_1\cup K_2)$. We may write $(1+\epsilon)\xb=\lambda \xb_1+(1-\lambda)\xb_2$, with $\xb_1\in K_1-y$, $\xb_2\in K_2-y$, and $0\le\lambda\le 1$. Then $\xb = \lambda\frac{1}{1+\epsilon} \xb_1+(1-\lambda)\frac{1}{1+\epsilon}\xb_2$. By \eqref{characri-ii}, $\frac{1}{1+\epsilon} \xb_1\in\ri K_1$, and $\frac{1}{1+\epsilon} \xb_2\in\ri K_2$. $\xb\in \conv\big(\ri (K_1-y)\cup\ri (K_2-y)\big)$, and therefore $x\in \conv(\ri K_1\cup\ri K_2)$.
\ep\\

Now we use the measurable selection theory to establish the non-emptiness of $\partial f$.

\begin{Lemma}\label{lemma:partialnonempty}
For all $f\in\Cfrak$, we have $\partial f\neq\emptyset$.
\end{Lemma}
\proof By the fact that $f$ is continuous, we may write $\partial f(x) = \cap_{n\ge 1}F_n(x)$ for all $x\in\R^d$, with $F_n(x):=\{p\in\R^d:f(y_n)-f(x)\ge p\cdot (y_n-x)\}$ where $(y_n)_{n\ge 1}\subset\R^d$ is some fixed dense sequence. All $F_n$ are measurable by the continuity of $(x,p)\longmapsto f(y_n)-f(x)-p\cdot (y_n-x)$ together with Theorem 6.4 in \cite{himmelberg1975measurable}. Therefore the mapping $x\longmapsto\partial f(x)$ is measurable by Lemma \ref{lemma:unionbim}. Moreover, the fact that this mapping is closed nonempty-valued is a well-known property of the subgradient of finite convex functions in finite dimension. Then the result holds by Theorem 4.1 in \cite{wagner1977survey}. 
\ep\\

We conclude this section with the following result which has been used in our proof of Proposition \ref{prop:Komlos}. We believe that this is a standard convex analysis result, but we could not find precise references. For this reason, we report the proof for completeness.

\begin{Theorem}\label{thm:ConvexConv}
Let $f_n,f:\R^d\to\bar\R$ be convex functions with $\interior\,\dom f \neq \emptyset$. Then $f_n \longrightarrow f$ pointwise on $\R^d\setminus \partial \dom f$ if and only if $f_n \longrightarrow f$ pointwise on some dense subset $A\subset \R^d\setminus \partial \dom f$.
\end{Theorem}

\proof
We prove the non-trivial implication "if". We first prove the convergence on $\interior\,\dom f$. $f_n$ converges to $f$ on a dense set. The reasoning will consist in proving that the $f_n$ are Lipschitz, it will give a uniform convergence and then a pointwise convergence. First we consider $K\subset \interior\,\dom f$ compact convex with nonempty interior. We can find $N\in\N$ and $x_1,...x_N\in A\cap(\interior\,\dom f \setminus K)$ such that $K\subset \interior\, \conv(x_1,...,x_N)$. We use the pointwise convergence on $A$ to get that for $n$ large enough, $f_n(x)\leq M$ for $x\in \conv(x_1,...,x_N)$, $M>0$ (take $M=\max_{1\leq k\leq N}f(x_k)+1$). Then we will prove that $f_n$ is bounded from below on $K$. We consider $a\in A\cap K$ and $\delta_0:=\sup_{x\in K}|x-a|$. For $n$ large enough, $f_n(a)\geq m$ for any $a\in A$ (take for example $m=f(a)-1$). We write $\delta_1 := \min_{(x,y)\in K\x\partial \conv(x_1,...,x_N)}|x-y|$. Finally we write $\delta_2:= \sup_{x,y\in \conv(x_1,...,x_N)}|x-y|$. Now, for $x\in K$, we consider the half line $x+\R_+(a-x)$, it will cut $\partial \conv(x_1,...,x_N)$ in one only point $y\in \partial \conv(x_1,...,x_N)$. Then $a\in[x,y]$, and therefore $a= \frac{|a-y|}{|x-y|}x+\frac{|a-x|}{|x-y|}y$. By the convex inequality, $f_n(a)\leq \frac{|a-y|}{|x-y|} f_n(x)+\frac{|a-x|}{|x-y|}f_n(y)$. Then $f_n(x)
\geq -\frac{|a-x|}{|a-y|} M + \frac{|x-y|}{|a-y|} m
\ge -\frac{\delta_0}{\delta_1} M + \frac{\delta_2}{\delta_1} m$. Finally, if we write $m_0:=-\frac{\delta_0}{\delta_1} M + \frac{\delta_2}{\delta_1} m$,
\b*
M\geq f_n\geq m_0,&\mbox{on }K.
\e*
This will prove that $f_n$ is $\frac{M-m0}{\delta_1}$-Lipschitz. We consider $x\in K$ and a unit direction $u\in\Sc^{d-1}$ and $f_n'\in\partial f_n(x)$. For a unique $\lambda > 0$, $y:=x+\lambda u\in \partial \conv(x_1,...,x_N)$. As $u$ is a unit vector, $\lambda = |y-x|\geq \delta_1$. By the convex inequality, $f_n(y)\geq f_n(x)+f_n'(x)\cdot (y-x)$.
Then $M-m_0\geq \delta_0 |f_n'\cdot u|$ and finally $|f_n'\cdot u|\leq \frac{M-m0}{\delta_1}$ as this bound does not depend on $u$, $|f_n'|\leq \frac{M-m0}{\delta_1}$ for any such subgradient. For $n$ large enough, the $f_n$ are uniformly Lipschitz on $K$, and so in $f$. The convergence is uniform on $K$, it is then pointwise on $K$. As this is true for any such $K$, the convergence is pointwise on $\interior\,\dom f$.\\
Now let us consider $x\in (\cl\dom f)^c$. The set $\conv(x,\interior\,\dom f)\setminus \dom f$ has a nonempty interior because $\dist(x,\dom f)>0$ and $\interior\,\dom f\neq \emptyset$. As $A$ is dense, we can consider $a\in A\cap \conv(x,\interior\,\dom f)\setminus \dom f$. By definition of $\conv(x,\interior\,\dom f)$, we can find $y\in\interior\,\dom f$ such that $a=\lambda y + (1-\lambda) x$. We have $\lambda <1$ because $a\notin \dom f$. If $\lambda = 0$, $f_n(x)=f_n(a)\underset{n\to\infty}{\longrightarrow}\infty$. Otherwise, by the convexity inequality, $f_n(a)\leq \lambda f_n(y) + (1-\lambda) f_n(x)$. Then, as $f_n(a)\underset{n\to\infty}{\longrightarrow}\infty$, and $f_n(y)\underset{n\to\infty}{\longrightarrow}f(y)<\infty$, we have $f_n(x)\underset{n\to\infty}{\longrightarrow}\infty$.
\ep

\section*{Acknowledgements}

We are grateful to the two anonymous referees, whose fruitful remarks and comments contributed to enhance deeply this paper.

\bibliographystyle{plain}
\bibliography{mabib}

\begin{thebibliography}{10}

\bibitem{beer1991polish}
Gerald Beer.
\newblock A polish topology for the closed subsets of a polish space.
\newblock {\em Proceedings of the American Mathematical Society},
  113(4):1123--1133, 1991.

\bibitem{beiglbock2013model}
Mathias Beiglb{\"o}ck, Pierre Henry-Labord{\`e}re, and Friedrich Penkner.
\newblock Model-independent bounds for option prices: a mass transport
  approach.
\newblock {\em Finance and Stochastics}, 17(3):477--501, 2013.

\bibitem{beiglboeck2016problem}
Mathias Beiglb{\"o}ck and Nicolas Juillet.
\newblock On a problem of optimal transport under marginal martingale
  constraints.
\newblock {\em The Annals of Probability}, 44(1):42--106, 2016.

\bibitem{beiglbock2015complete}
Mathias Beiglb{\"o}ck, Marcel Nutz, and Nizar Touzi.
\newblock Complete duality for martingale optimal transport on the line.
\newblock {\em arXiv preprint arXiv:1507.00671}, 2015.

\bibitem{bertsekas1978stochastic}
Dimitri~P Bertsekas and Steven~E Shreve.
\newblock {\em Stochastic optimal control: The discrete time case}, volume~23.
\newblock Academic Press New York, 1978.

\bibitem{cox2011robust}
Alexander~MG Cox and Jan Ob{\l}{\'o}j.
\newblock Robust pricing and hedging of double no-touch options.
\newblock {\em Finance and Stochastics}, 15(3):573--605, 2011.

\bibitem{delbaen1994general}
Freddy Delbaen and Walter Schachermayer.
\newblock A general version of the fundamental theorem of asset pricing.
\newblock {\em Mathematische annalen}, 300(1):463--520, 1994.

\bibitem{ekren2016constrained}
Ibrahim Ekren and H~Mete Soner.
\newblock Constrained optimal transport.
\newblock {\em arXiv preprint arXiv:1610.02940}, 2016.

\bibitem{galichon2014stochastic}
Alfred Galichon, Pierre Henry-Labordere, and Nizar Touzi.
\newblock A stochastic control approach to no-arbitrage bounds given marginals,
  with an application to lookback options.
\newblock {\em The Annals of Applied Probability}, 24(1):312--336, 2014.

\bibitem{ghoussoub2015structure}
Nassif Ghoussoub, Young-Heon Kim, and Tongseok Lim.
\newblock Structure of optimal martingale transport plans in general
  dimensions.
\newblock {\em arXiv preprint arXiv:1508.01806}, 2015.

\bibitem{hess1986contribution}
Christian Hess.
\newblock {\em Contribution {\`a} l'{\'e}tude de la mesurabilit{\'e}, de la loi
  de probabilit{\'e} et de la convergence des multifonctions}.
\newblock PhD thesis, 1986.

\bibitem{himmelberg1975measurable}
Charles~J Himmelberg.
\newblock Measurable relations.
\newblock {\em Fundamenta mathematicae}, 87(1):53--72, 1975.

\bibitem{hiriart2013convex}
Jean-Baptiste Hiriart-Urruty and Claude Lemar{\'e}chal.
\newblock {\em Convex analysis and minimization algorithms I: Fundamentals},
  volume 305.
\newblock Springer science \& business media, 2013.

\bibitem{hobson2011skorokhod}
David Hobson.
\newblock The skorokhod embedding problem and model-independent bounds for
  option prices.
\newblock In {\em Paris-Princeton Lectures on Mathematical Finance 2010}, pages
  267--318. Springer, 2011.

\bibitem{hobson2015robust}
David Hobson and Martin Klimmek.
\newblock Robust price bounds for the forward starting straddle.
\newblock {\em Finance and Stochastics}, 19(1):189--214, 2015.

\bibitem{hobson2012robust}
David Hobson and Anthony Neuberger.
\newblock Robust bounds for forward start options.
\newblock {\em Mathematical Finance}, 22(1):31--56, 2012.

\bibitem{kellerer1984duality}
Hans~G Kellerer.
\newblock Duality theorems for marginal problems.
\newblock {\em Zeitschrift f{\"u}r Wahrscheinlichkeitstheorie und verwandte
  Gebiete}, 67(4):399--432, 1984.

\bibitem{OblojSiopraes2017}
Jan Obl{\'o}j and Pietro Siorpaes.
\newblock Structure of martingale transports in finite dimensions.
\newblock 2017.

\bibitem{rockafellar2015convex}
Ralph~Tyrell Rockafellar.
\newblock {\em Convex analysis}.
\newblock Princeton university press, 2015.

\bibitem{strassen1965existence}
Volker Strassen.
\newblock The existence of probability measures with given marginals.
\newblock {\em The Annals of Mathematical Statistics}, pages 423--439, 1965.

\bibitem{wagner1977survey}
Daniel~H Wagner.
\newblock Survey of measurable selection theorems.
\newblock {\em SIAM Journal on Control and Optimization}, 15(5):859--903, 1977.

\bibitem{zaev2015monge}
Danila~A Zaev.
\newblock On the monge--kantorovich problem with additional linear constraints.
\newblock {\em Mathematical Notes}, 98(5-6):725--741, 2015.

\end{thebibliography}

\end{document}